\documentclass[11pt,a4,fleqn]{article}
\usepackage{graphicx}
\input{amssym}
\newtheorem{theorem}{\bf Theorem}[section]

\newtheorem{lemma}[theorem]{\bf Lemma}

\newtheorem{emp}[theorem]{\bf Claim }
\newtheorem{remark}[theorem]{\bf Remark}

\newtheorem{conjecture}[theorem]{\bf Conjecture}
\newtheorem{question}[theorem]{\bf Question}

\newtheorem{preproof}{{\bf Proof.}}

\begin{document}
\title{{\Large Diagonal Ramsey numbers of   loose cycles in uniform hypergraphs}}
\author{\small  G.R. Omidi$^{\textrm{a},\textrm{b},1}$, M. Shahsiah$^{\textrm{b},2}$\\
\small  $^{\textrm{a}}$Department of Mathematical Sciences,
Isfahan University
of Technology,\\ \small Isfahan, 84156-83111, Iran\\
\small  $^{\textrm{b}}$School of Mathematics, Institute for
Research in Fundamental Sciences (IPM),\\
\small  P.O.Box: 19395-5746, Tehran,
Iran\\
\small \texttt{E-mails: romidi@cc.iut.ac.ir, shahsiah@ipm.ir}}
\date {}
\maketitle \footnotetext[1] {\tt This research is partially carried out in the IPM-Isfahan Branch and in part
 supported by a grant from IPM (No. 91050416).} \vspace*{-0.5cm}

\begin{abstract}
\hspace{-0.25 cm}A  $k$-uniform  loose cycle $\mathcal{C}_n^k$  is
a hypergraph with vertex set $\{v_1,v_2,\ldots,v_{n(k-1)}\}$ and
with the set of $n$ edges
$e_i=\{v_{(i-1)(k-1)+1},v_{(i-1)(k-1)+2},\ldots,
v_{(i-1)(k-1)+k}\}$, $1\leq i\leq n$, where we use mod $n(k-1)$
arithmetic. The  Ramsey number
$R(\mathcal{C}^k_n,\mathcal{C}^k_n)$ is asymptotically
$\frac{1}{2}(2k-1)n$ as has been proved by
 Gy\'{a}rf\'{a}s, S\'{a}rk\"{o}zy and   Szemer\'{e}di.
 In this paper, we investigate to determining the exact value of
 diagonal Ramsey number of $\mathcal{C}^k_n$ and we show that for $n\geq 2$ and  $k\geq 8$
$$R(\mathcal{C}^k_n,\mathcal{C}^k_n)=(k-1)n+\lfloor\frac{n-1}{2}\rfloor.$$

\noindent{\small {\bf Keywords:} Ramsey number, Uniform hypergraph, Loose path, Loose cycle.}\\
{\small AMS subject classification: 05C65, 05C55, 05D10.}

\end{abstract}

\section{\normalsize Introduction}

\bigskip
For given $k$-uniform hypergraphs $\mathcal{H}_1$ and
$\mathcal{H}_2$,  the \textit{Ramsey number}
$R(\mathcal{H}_1,\mathcal{H}_2)$ is the smallest number $N$ such
that in every red-blue coloring of the
  edges of the complete $k$-uniform hypergraph $\mathcal{K}^k_N$ there is  a red copy of $\mathcal{H}_1$ or a blue copy of $\mathcal{H}_2.$
 A {\it $k$-uniform loose path} $\mathcal{P}_n^k$ (shortly, a {\it
path of length $n$}) is a hypergraph with vertex set
$\{v_1,v_2,\ldots,v_{n(k-1)+1}\}$  and with the set of $n$ edges
$e_i=\{v_1,v_2,\ldots, v_k\}+(i-1)(k-1)$, $i=1,2,\ldots, n$.
 For an edge $e$ of a given loose path (also a given loose cycle)
$\mathcal{K}$, the first
 vertex and the last vertex are denoted by $f_{\mathcal{K},e}$ and
 $l_{\mathcal{K},e}$, respectively.  For $k=2$ we get the usual definitions of a cycle $C_n$ and a path $P_n$ with $n$
 edges.\\

One of the most important problems in combinatorics and graph
theory is determining or estimating the Ramsey numbers. In
contrast to the graph case, there are  few known  results about
the Ramsey numbers of hypergraphs. Recently, several interesting
results were obtained on the Ramsey numbers of loose cycles in
uniform hypergraphs. Haxell et al. \cite{Ramsy number of loose
cycle} showed that
the  Ramsey number of $3$-uniform loose
cycles  is asymptotically $\frac{5n}{2}$.
 More precisely, they proved that for
all $\eta>0$ there exists $n_0=n_0(\eta)$ such that for every $n>
n_0$, every $2$-coloring of $\mathcal{K}^{3}_{ 5(1+\eta)n/2}$
contains a monochromatic copy of $\mathcal{C}_n^3$.
 As
extension of this result proved by  Gy\'{a}rf\'{a}s,
S\'{a}rk\"{o}zy and Szemer\'{e}di \cite{Ramsy number of loose
cycle for k-uniform}, shows that the  Ramsey number of
$k$-uniform loose cycles is asymptotically $\frac{1}{2}(2k-1)n$.
However, those proofs are all based on  the method of the
Regularity Lemma. Concerning to  exact values of  Ramsey numbers
of loose paths and cycles with at most $4$ edges we have the
following result.
\begin{theorem}\label{R(Pk3,Pk3)}
For every $k\geq 3$,
\begin{itemize}
\item[\rm(i)] $R(\mathcal{C}^k_2,\mathcal{C}^k_2)+1=2k-2$ {\rm (to
see a proof see Appendix $A$)}. \item[\rm(ii)]{\rm (\cite{subm})}
$R(\mathcal{P}^k_3,\mathcal{P}^k_3)=R(\mathcal{C}^k_3,\mathcal{P}^k_3)=R(\mathcal{C}^k_3,\mathcal{C}^k_3)+1=3k-1$.
\item[\rm(iii)]{\rm
(\cite{subm})}$R(\mathcal{P}^k_4,\mathcal{P}^k_4)=R(\mathcal{C}^k_4,\mathcal{P}^k_4)=R(\mathcal{C}^k_4,\mathcal{C}^k_4)+1=4k-2$.
\end{itemize}
\end{theorem}
Regarding to Ramsey numbers of $3$-uniform loose paths and cycles with arbitrary number of edges,
in  \cite{subm}, the authors  posed the following question.
\begin{question}\label{question} For every $n\geq m\geq 3$, is it
true that
\begin{eqnarray*}
R(\mathcal{P}^3_n,\mathcal{P}^3_m)=R(\mathcal{P}^3_n,\mathcal{C}^3_m)=R(\mathcal{C}^3_n,\mathcal{C}^3_m)+1=2n+\Big\lfloor\frac{m+1}{2}\Big\rfloor?
\end{eqnarray*}
In particular, is it true that
\begin{eqnarray*}
R(\mathcal{P}^3_n,\mathcal{P}^3_n)=R(\mathcal{C}^3_n,\mathcal{C}^3_n)+1=\Big\lceil\frac{5n}{2}\Big\rceil?
\end{eqnarray*}
\end{question}
Recently, this question is answered positively (see  \cite{The Ramsey number of loose paths  and loose cycles
in 3-uniform hypergraphs} and \cite{Ramsey numbers of loose cycles
in uniform hypergraphs}).
 In \cite{Ramsey numbers of loose cycles in
uniform hypergraphs} the authors   posed the following
conjecture on the  values of  Ramsey numbers of $k$-uniform  loose paths and
cycles for $k\geq 3$.

\begin{conjecture}\label{our conjecture}
Let $k\geq 3$ be an integer number. For every $n\geq m \geq 3$,
\begin{eqnarray}\label{2}
R(\mathcal{P}^k_n,\mathcal{P}^k_m)=R(\mathcal{P}^k_n,\mathcal{C}^k_m)=R(\mathcal{C}^k_n,\mathcal{C}^k_m)+1=(k-1)n+\lfloor\frac{m+1}{2}\rfloor.
\end{eqnarray}
\end{conjecture}
In \cite{Ramsey numbers of loose cycles in uniform hypergraphs},
the authors  also demonstrated  that

\begin{theorem}\label{connection}
For $n\geq m,$ if
$R(\mathcal{C}^k_{n},\mathcal{C}^k_{m})=(k-1)n+\lfloor\frac{m-1}{2}\rfloor,$
then
\begin{eqnarray*}\label{222}
R(\mathcal{P}^k_n,\mathcal{P}^k_m)=R(\mathcal{P}^k_n,\mathcal{C}^k_m)=(k-1)n+\lfloor\frac{m+1}{2}\rfloor.
\end{eqnarray*}
\end{theorem}
So they digested Conjecture \ref{our conjecture} to the following
conjecture.

\begin{conjecture}\label{our conjecture2}
Let $k\geq 3$ be an integer number. For every $n\geq m \geq 3$,
\begin{eqnarray*}\label{3}
R(\mathcal{C}^k_n,\mathcal{C}^k_m)=(k-1)n+\lfloor\frac{m-1}{2}\rfloor.
\end{eqnarray*}
\end{conjecture}

In this paper, we investigate Conjecture \ref{our conjecture2} for
the case $n=m$. More precisely, we prove the following theorem.
\begin{theorem}\label{main theorem}
For every integers  $k\geq 8$ and $n\geq 2$
\begin{eqnarray*}\label{3}
R(\mathcal{C}^k_n,\mathcal{C}^k_n)=(k-1)n+\lfloor\frac{n-1}{2}\rfloor.
\end{eqnarray*}
\end{theorem}

Using Theorems \ref{connection} and \ref{main theorem} we conclude
the following result.

\begin{theorem}\label{conclusion}
For every integers  $k\geq 8$ and $n\geq 3$
\begin{eqnarray*}
R(\mathcal{P}^k_n,\mathcal{P}^k_n)=R(\mathcal{P}^k_n,\mathcal{C}^k_n)=R(\mathcal{C}^k_n,\mathcal{C}^k_n)+1=(k-1)n+\lfloor\frac{n+1}{2}\rfloor.
\end{eqnarray*}
\end{theorem}

Using Lemma 1 of \cite{subm}, we have
$R(\mathcal{C}^k_n,\mathcal{C}^k_m)\geq
(k-1)n+\lfloor\frac{m-1}{2}\rfloor$.  Therefore, in order to  prove Theorem \ref{main theorem}, it
 suffices to verify  that the known lower bound are also   upper bound.

Throughout the paper, we denote by $\mathcal{H}_{\rm red}$ and
$\mathcal{H}_{\rm blue}$ the induced $k$-uniform hypergraphs on
the edges of $\mathcal{H}$ with color red and blue, respectively.
Also we denote by $|\mathcal{H}|$ and $\|\mathcal{H}\|$ the number
of vertices and edges of $\mathcal{H}$, respectively.\\

The rest of this paper is organized as follows. In the next section, we give an outline of the proof of Theorem \ref{main theorem}.
In Section 3, we
state some definitions  and  key technical lemmas
 required to prove the main theorem. For the sake of clarity of presentation, we  omit the proofs
  in Section 3 and we refer the reader to Appendix $A$ to see the complete proofs.
  We will present a complete proof of
 Theorem \ref{main theorem} in Section 4. In section $5$, we  conclude with some further remarks and open problems.

\section{Outline of the proof of Theorem \ref{main theorem}}
In this section, we sketch the main idea of our proof for Theorem \ref{main theorem}. We give a
proof by induction on the number of vertices. By Theorem
\ref{R(Pk3,Pk3)} we may assume that $n\geq 5$. Let
$f(n)=(k-1)n+\lfloor\frac{n-1}{2}\rfloor$ and suppose to the
contrary that $\mathcal{H}=\mathcal{K}^k_{f(n)}$ is 2-edge colored
red and blue with  no monochromatic  copy of $\mathcal{C}^k_n$.
 We consider three following
cases.\\


\noindent{\bf Case 1: $n\equiv 1,2,3$ (mod 4)}\\

First
we show that there are two monochromatic copies of
$\mathcal{C}^k_{\lfloor\frac{n}{2}\rfloor+1}$ of colors red and blue (see Claims
\ref{red-blue C} and \ref{red-blue C_{(n/2)+1}}). Then
 among all red-blue copies of
$\mathcal{C}^k_{\lfloor\frac{n}{2}\rfloor+1}$'s choose red-blue copies with
maximum intersection, say $\mathcal{C}_1$ and $\mathcal{C}_2$.
It is not
difficult to see that
\begin{eqnarray*}
|V(\mathcal{C}_1\cup \mathcal{C}_2)|\leq
R(\mathcal{C}^k_{\lfloor\frac{n}{2}\rfloor+1},\mathcal{C}^k_{\lfloor\frac{n}{2}\rfloor+1})+1.
\end{eqnarray*}
Clearly $|V(\mathcal{H})\setminus V(\mathcal{C}_1\cup
\mathcal{C}_2)|\geq f(n-1-\lfloor\frac{n}{2}\rfloor)$.
 So, by  induction hypothesis, there is a monochromatic  $\mathcal{C}^k_{n-1-\lfloor\frac{n}{2}\rfloor}$,
 say $\mathcal{C}_3=e_1e_2\ldots e_{n-1-\lfloor\frac{n}{2}\rfloor}$,
 disjoint from $\mathcal{C}_1$ and $\mathcal{C}_2$.
  By symmetry we may assume that $\mathcal{C}_2$ and
  $\mathcal{C}_3$ are both in $\mathcal{H}_{\rm
  blue}$ and  $\mathcal{C}_1\subseteq \mathcal{H}_{\rm
  red}$.\\

Since there is no red copy of
  $\mathcal{C}^k_n$, we can find a red copy of $\mathcal{C}^k_{n-1-\lfloor\frac{n}{2}\rfloor}$, say $\mathcal{C}_4=h_1h_2 \ldots
  h_{n-1-\lfloor\frac{n}{2}\rfloor}$, so that for each $1\leq i \leq
  n-1-\lfloor\frac{n}{2}\rfloor$, $k-2 \leq|h_i\cap e_i|\leq k-1$. For even $n$, this follows from Claims  \ref{28} and \ref{32}. For odd $n$ we use Claims \ref{12} and \ref{13}.
In the rest  of the   proof, we find a red copy of
$\mathcal{C}^k_{n}$ by joining $\mathcal{C}_1$ and
$\mathcal{C}_4$. If $n=6$, we do this  by applying Claims
 \ref{21}, \ref{32} and \ref{26}. Otherwise, we  apply Claims
\ref{12}, \ref{13} and \ref{14} when $n\equiv^4 1,3$  and   Claims \ref{21}, \ref{28} and \ref{26} when $n\equiv^4 2$.
 This is a contradiction to our
assumption. So we are done.\\\\

\noindent{\bf Case 2: $n\equiv 0$ (mod 4)}\\

 In this case, first  we show that
 there are two disjoint isochromatic
$\mathcal{C}^k_{\frac{n}{2}}$. By symmetry we can suppose that
$\mathcal{C}_1=g_1g_2\ldots g_{\frac{n}{2}}$ and
$\mathcal{C}_2=h_1h_2\ldots
  h_{\frac{n}{2}}$
are such cycles in $\mathcal{H}_{\rm blue}$. As we assumed that
there is no blue copy of $\mathcal{C}^k_n$, we make a copy of
$\mathcal{C}^k_n$ in $\mathcal{H}_{\rm red}$ as follows.
Let $W=V(\mathcal{H})\setminus
V(\mathcal{C}_1\cup \mathcal{C}_2)$. Since $n\geq 8$, we have  $|W|\geq 3$.
Use Lemma \ref{there is a P2:3} for $e_i=g_1$, $f_j=h_1$ and
$B=\{w_1,w_2,w_3\}\subseteq W$ to obtain two red paths $E_1$ and
$F_1$ with desired properties.  Assume that
$$g_i=\{x_1,x_2,\ldots,x_k\}+(k-1)(i-1) ({\rm mod}\ (k-1)\frac{n}{2}),\hspace{0.5 cm} i=1,2,\ldots,\frac{n}{2},$$
$$h_i=\{y_1,y_2,\ldots,y_k\}+(k-1)(i-1) ({\rm mod}\  (k-1)\frac{n}{2}),\hspace{0.5 cm} i=1,2,\ldots,\frac{n}{2}.$$
Now, use Lemma \ref{there is a P:3} for $e_i=g_2$, $f_j=h_2$ and
$\mathcal{E}_1=E_1$ (resp. $\mathcal{E}_1=F_1$) to obtain two red
paths $\mathcal{E}_2$ and $\mathcal{F}_2$ (resp.
$\overline{\mathcal{E}}_2$ and $\overline{\mathcal{F}}_2$) of
length $4$ with mentioned properties of Lemma \ref{there is a
P:3} (this can be done by a suitable renaming of the edges of $\mathcal{C}_1$ and $\mathcal{C}_2$). In the next step, for $2\leq l \leq \frac{n}{2}-2$, use Lemma
\ref{there is a P}  for red paths $\mathcal{E}_l$, $\mathcal{F}_l$ (resp.
 $\overline{\mathcal{E}}_l$ and $\overline{\mathcal{F}}_l$),
$e_i=g_{l+1}$ and $f_j=h_{l+1}$ to obtain two red paths
$\mathcal{E}_{l+1}$ and $\mathcal{F}_{l+1}$ (resp.
$\overline{\mathcal{E}}_{l+1}$ and $\overline{\mathcal{F}}_{l+1}$)
with properties of Lemma \ref{there is a P}. Note that each of
$\mathcal{E}_{\frac{n}{2}-1}$, $\mathcal{F}_{\frac{n}{2}-1}$,
$\overline{\mathcal{E}}_{\frac{n}{2}-1}$ and
$\overline{\mathcal{F}}_{\frac{n}{2}-1}$ has $n-2$ edges.

Let $\mathcal{E}_{\frac{n}{2}-1}=p_1p_2\ldots p_{n-2}$, $x\in
g_{\frac{n}{2}-1}\setminus(\mathcal{E}_{\frac{n}{2}-1}\cup\{f_{\mathcal{C}_1,g_{\frac{n}{2}-1}}\})$,
$y'\in (h_1\setminus\{y_k\})\cap(p_1\setminus p_2)$ and  $y''\in
 (h_{\frac{n}{2}-1}\setminus\{f_{\mathcal{C}_2,h_{\frac{n}{2}-1}}\})\cap(p_{n-2}\setminus
 p_{n-3})$. Now,  we define two edges $q$ and $q'$ as
 follows. Let
\begin{eqnarray*}
q=\{x,x_{(k-1)(\frac{n}{2}-1)+2},\ldots,x_{(k-1)(\frac{n}{2}-1)+\lfloor\frac{k}{2}\rfloor}\}\cup
\{y_{(k-1)\frac{n}{2}-\lceil\frac{k}{2}\rceil+2},\ldots,y_{(k-1)\frac{n}{2}},y'\},\\
\end{eqnarray*}
and
\begin{eqnarray*}
q'&=&\{x_{(k-1)(\frac{n}{2}-1)+\lfloor\frac{k}{2}\rfloor},x_{(k-1)(\frac{n}{2}-1)+\lfloor\frac{k}{2}\rfloor+2},x_{(k-1)(\frac{n}{2}-1)+\lfloor\frac{k}{2}\rfloor+3}
 ,\ldots,x_{(k-1)\frac{n}{2}},v\}\\
&&\cup\{y'',y_{(k-1)(\frac{n}{2}-1)+2},\ldots,y_{(k-1)(\frac{n}{2}-1)+\lfloor\frac{k}{2}\rfloor}\}.
\end{eqnarray*}
If at least one of $q$ or $q'$ is blue, then we can find either a red copy of $\mathcal{C}^k_{n}$ (using
the red paths $\mathcal{E}_{\frac{n}{2}-1}$,
$\mathcal{F}_{\frac{n}{2}-1}$,
$\overline{\mathcal{E}}_{\frac{n}{2}-1}$ and
$\overline{\mathcal{F}}_{\frac{n}{2}-1}$) or a blue copy of $\mathcal{C}^k_{n}$ (using the blue cycles $\mathcal{C}_1$ and $\mathcal{C}_2$). Hence,
  $\mathcal{E}_{\frac{n}{2}-1}q'q$ is
 a red copy of $\mathcal{C}^k_{n}$. This contradicts our
 assumption.\\\\

\section{\normalsize Preliminaries}

\medskip

In this section, we present some results that will be used later
on. For the clarity of presentation, except Lemmas \ref{there is a P2:2} and \ref{there is a P:2} (their proofs are completely similar to the proofs of Lemmas \ref{there is a P2} and \ref{there is a P}, respectively), the proofs of all results in this section will be given in
 Appendix $A$.

\medskip
\begin{lemma}\label{No  red Cn-1 implies blue Cn}
Let $n\geq 3$, $k\geq 6$, $f\geq (k-1)n$ and $\mathcal{H}=\mathcal{K}^k_{f}$
be $2$-edge colored red and blue.
 If there is a $\mathcal{C}_{n}^k\subseteq
\mathcal{H}_{\rm red}$ and there is
  no red copy of
$\mathcal{C}_{n-1}^k$,  then $\mathcal{C}_{n}^k\subseteq
\mathcal{H}_{\rm blue}.$
\end{lemma}


\medskip
\begin{lemma}\label{No  red Cn-2 implies blue Cn}
Let $n\geq 3$, $k\geq 6$, $f\geq (k-1)n$ and $\mathcal{H}=\mathcal{K}^k_{f}$
be $2$-edge colored red and blue. If there is a $\mathcal{C}_{n}^k\subseteq
\mathcal{H}_{\rm red}$ and there is  no red copy of
$\mathcal{C}_{n-2}^k$,  then
$\mathcal{C}_{n}^k\subseteq\mathcal{H}_{\rm blue}$.
\end{lemma}


 \begin{remark}\label{assumption}
  In the sequel of this section  assume that $n\geq 5$, $k\geq 8$, $l_2\geq l_1\geq 2$, $l_1+l_2\leq  n$ and
  $\mathcal{H}=\mathcal{K}^k_{(k-1)n+\lfloor\frac{n-1}{2}\rfloor}$ is $2$-edge
colored red and blue. Also, let $\mathcal{C}_1=e_1e_2\ldots
e_{l_1}$ and $\mathcal{C}_2=f_1f_2\ldots f_{l_2}$ be two disjoint
cycles in $\mathcal{H}_{\rm blue}$ with edges

\begin{eqnarray*}
e_i=\{v_1,v_2,\ldots,v_k\}+(k-1)(i-1) \ \ {\rm (mod} (k-1)l_1{\rm
)}, \hspace{0.5 cm}i=1,\ldots,l_1
\end{eqnarray*}
 and
\begin{eqnarray*}
f_i=\{u_1,u_2,\ldots,u_k\}+(k-1)(i-1) \ \ {\rm (mod} (k-1)l_2{\rm
)}, \hspace{0.5 cm}i=1,\ldots,l_2,
\end{eqnarray*}
respectively, and $W$ with $|W|\geq 2$ be the set of vertices
uncovered by $\mathcal{C}_1\cup \mathcal{C}_2$.
\end{remark}

\noindent Let $e_i$ and $f_j$ be two edges of $\mathcal{C}_1$ and
$\mathcal{C}_2$, respectively and  $g=E\dot{\cup}W'\dot{\cup}F$,
where

\begin{itemize}
\item[$\bullet$] $E=E'\cup\{v'\}$ for $v'\in
(e_{i-1}\setminus\{f_{\mathcal{C}_1,e_{i-1}}\})\cup(e_{i+1}\setminus\{l_{\mathcal{C}_1,e_{i+1}}\})$
and $E'\subseteq V(e_i)\setminus
\{f_{\mathcal{C}_1,e_{i}},l_{\mathcal{C}_1,e_{i}}\}$.

\item[$\bullet$] $F=F'\cup\{u'\}$ for $u'\in
(f_{j-1}\setminus\{f_{\mathcal{C}_2,f_{j-1}}\})\cup(f_{j+1}\setminus\{l_{\mathcal{C}_2,f_{j+1}}\})$
and $F'\subseteq
V(f_j)\setminus\{f_{\mathcal{C}_2,f_{j}},l_{\mathcal{C}_2,f_{j}}\}$.
\item[$\bullet$] $W'\subseteq W$.

 \item[$\bullet$] $|E|=p\geq
1,|W'|=q\geq 0,|F|=r\geq 1$ and $p+q+r=k$.

\end{itemize}
We say that  $g$ is of type $A$, $B$, $C$ or $D$ if $(v',u')\in
(e_{i-1}\setminus\{f_{\mathcal{C}_1,e_{i-1}}\})\times
(f_{j+1}\setminus\{l_{\mathcal{C}_2,f_{j+1}}\})$,
$(v',u')\in(e_{i+1}\setminus\{l_{\mathcal{C}_1,e_{i+1}}\})\times(f_{j-1}\setminus\{f_{\mathcal{C}_2,f_{j-1}}\})$,
$(v',u')\in
(e_{i-1}\setminus\{f_{\mathcal{C}_1,e_{i-1}}\})\times(f_{j-1}\setminus\{f_{\mathcal{C}_2,f_{j-1}}\})$
or $(v',u')\in
(e_{i+1}\setminus\{l_{\mathcal{C}_1,e_{i+1}}\})\times(f_{j+1}\setminus\{l_{\mathcal{C}_2,f_{j+1}}\})$,
respectively. By $\mathcal{A}_{ij}$, $\mathcal{B}_{ij}$,
$\mathcal{C}_{ij}$, $\mathcal{D}_{ij}$ we mean all edges of type
$A$, $B$, $C$, $D$ corresponding to the edges $e_i$ and $f_j$,
respectively.





\begin{remark}\label{edge complement}
Consider the edges $e_i$ and $f_j$. For  $g\in \mathcal{A}_{ij}$
{\rm(}resp. $g\in \mathcal{B}_{ij})$ and for every $v''\in
e_{i+1}\setminus(V(g)\cup\{l_{\mathcal{C}_1,e_{i+1}}\})$ {\rm(}resp.
$v''\in e_{i-1}\setminus(V(g)\cup\{f_{\mathcal{C}_1,e_{i-1}}\}))$ and
$u''\in f_{j-1}\setminus\{f_{\mathcal{C}_2,f_{j-1}}\}$ {\rm(}resp.
$u''\in f_{j+1}\setminus\{l_{\mathcal{C}_2,f_{j+1}}\})$,
there is an edge $g'\in \mathcal{B}_{ij}$ {\rm(}resp. $g'\in
\mathcal{A}_{ij})$ where $\{u'',v''\}\subseteq g'$ and $g\cap
g'=\emptyset$. The same result holds for the edges of types $C$
and $D$.

\end{remark}

\begin{remark}\label{complementary edges are not blue}
If there is no blue copy of $\mathcal{C}_{l_1+l_2}$ in
$\mathcal{H}$, then there are no two disjoint edges
$g\in\mathcal{A}_{ij}$ {\rm(}resp. $g\in\mathcal{C}_{ij})$ and
$g'\in\mathcal{B}_{ij}$ {\rm(}resp. $g'\in\mathcal{D}_{ij}${\rm)} of colors
blue in $\mathcal{H}$.
\end{remark}

\medskip
\begin{lemma}\label{there is a P2}
With the same assumptions in Remark {\rm \ref{assumption}}, let
$e_i$ and $f_j$ be two arbitrary edges of $\mathcal{C}_1$ and
$\mathcal{C}_2$, respectively, and $C\subseteq \{v\}$, where $v\in
e_i\setminus
\{f_{\mathcal{C}_1,e_i},l_{\mathcal{C}_1,e_i}\}$. Also, let
$B=\{w_1,w_2\}\subseteq W$  and let $e\in \{e_i,f_j\}$ be an edge
of $\mathcal{C}_r$ for some $r\in\{1,2\}$. Assume that
 $v'$, $v''$, $u'$ and $u''$ are distinct vertices so that
 $v'\in e_{i-1}\setminus\{f_{\mathcal{C}_1,e_{i-1}}\}$,
$v''\in e_{i+1}\setminus\{l_{\mathcal{C}_1,e_{i+1}}\}$,
$u'\in f_{j-1}\setminus\{f_{\mathcal{C}_2,f_{j-1}}\}$ and $u''\in
f_{j+1}\setminus\{l_{\mathcal{C}_2,f_{j+1}}\}$.  If there is
no blue copy of $\mathcal{C}_{l_1+l_2}$, then we can find a red
path $\mathcal{P}=g_1g_2$ so that for some vertex $w\in
e\setminus(\{f_{\mathcal{C}_r,e}\}\cup C)$, we have  $w\notin \mathcal{P}$
and the following conditions hold. \hspace*{-1 cm}\begin{itemize}
\item[$i)$] $V(\mathcal{P})\subseteq \Big(V(e_i)\setminus
(C\cup\{f_{\mathcal{C}_1,e_{i}},l_{\mathcal{C}_1,e_{i}}\})\Big)\cup
\Big(V(f_j)\setminus\{f_{\mathcal{C}_2,f_{j}},l_{\mathcal{C}_2,f_{j}}\}\Big)\cup
B'\cup\{v',v'',u',u''\}$ where $B'\subseteq B$ with $|B'|\leq
|C|+1$.
\item[$ii)$] $|g_m\cap
(e_i\setminus\{f_{\mathcal{C}_1,e_{i}},l_{\mathcal{C}_1,e_{i}}\})|\geq
2, |g_m\cap
(f_j\setminus\{f_{\mathcal{C}_2,f_{j}},l_{\mathcal{C}_2,f_{j}}\})|\geq
2,\ \ for \ \   m=1,2$.
 \item[$iii)$] If $|C|=1$, then either  $w_1\in g_1\setminus g_2$ and
 $w_2\in g_2\setminus g_1$ or $|B\cap g_1|=|B\cap g_2|+1=1$.
\end{itemize}
\end{lemma}

By an argument similar to the proof of  Lemma \ref{there is a P2} in Appendix $A$,
we have the following.
\medskip
\begin{lemma}\label{there is a P2:2}
With the same assumptions in Remark {\rm \ref{assumption}}, let
$e_i$ and $f_j$ be two arbitrary edges of $\mathcal{C}_1$ and
$\mathcal{C}_2$, respectively.  Assume that
 $v'$, $v''$, $u'$ and $u''$ are distinct vertices so that
 $v'\in e_{i-1}\setminus\{f_{\mathcal{C}_1,e_{i-1}}\}$,
$v''\in e_{i+1}\setminus\{l_{\mathcal{C}_1,e_{i+1}}\}$,
$u'\in f_{j-1}\setminus\{f_{\mathcal{C}_2,f_{j-1}}\}$ and $u''\in
f_{j+1}\setminus\{l_{\mathcal{C}_2,f_{j+1}}\}$.  If there is
no blue copy of $\mathcal{C}_{l_1+l_2}$, then we can find a red
path $\mathcal{P}=g_1g_2$ so that the following conditions hold.
\begin{itemize}
\item
 $V(\mathcal{P})\subseteq
(V(e_i)\setminus
\{f_{\mathcal{C}_1,e_{i}},l_{\mathcal{C}_1,e_{i}}\})\cup
(V(f_j)\setminus\{f_{\mathcal{C}_2,f_{j}},l_{\mathcal{C}_2,f_{j}}\})
\cup\{v',v'',u',u''\}.$ \item $|g_m\cap
(e_i\setminus\{f_{\mathcal{C}_1,e_{i}},l_{\mathcal{C}_1,e_{i}}\})|\geq
3, |g_m\cap
(f_j\setminus\{f_{\mathcal{C}_2,f_{j}},l_{\mathcal{C}_2,f_{j}}\})|\geq
3,\ \ for  \ \    m=1,2.$
\end{itemize}
\end{lemma}

 \medskip
\begin{lemma}\label{there is a P2:3}
With the same assumptions in Remark {\rm \ref{assumption}}, let
$e_i$ and $f_j$ be two arbitrary edges of $\mathcal{C}_1$ and
$\mathcal{C}_2$, respectively.
 Also, let
$B=\{w_1,w_2,w_3\}\subseteq W$.
If there is
no blue copy of $\mathcal{C}_{l_1+l_2}$, then we can find
two red paths $E_1=g_1g'_1$ and $F_1=g_1\overline{g'}_1$
 so that
 the following conditions hold. \hspace*{-1 cm}\begin{itemize}
\item $V(E_1),V(F_1)\subseteq\Big(V(e_i) \cup V(f_j)\cup
B'\Big)\setminus\{\overline{u},\overline{v}\}$ for some
$\overline{v}\in e_i$, $\overline{u}\in f_j$ and $B'\subseteq B$
with $|B'|=2$. \item
$e_{i}\setminus(V(E_{1})\cup\{f_{\mathcal{C}_1,e_{i}},l_{\mathcal{C}_1,e_{i}},\overline{v}\})\neq
\emptyset$ and
$f_{j}\setminus(V(F_{1})\cup\{f_{\mathcal{C}_2,f_{j}},l_{\mathcal{C}_2,f_{j}},\overline{u}\})\neq
\emptyset$. \item $|(g'_1\setminus g_1)\cap
(e_i\setminus\{f_{\mathcal{C}_1,e_{i}},l_{\mathcal{C}_1,e_{i}}\})|\geq
1$, $|(g'_1\setminus g_1)\cap
(f_j\setminus\{f_{\mathcal{C}_2,f_{j}},l_{\mathcal{C}_2,f_{j}}\})|\geq
1$. \item $|(g_1\setminus g'_1)\cap
(e_i\setminus\{f_{\mathcal{C}_1,e_{i}},l_{\mathcal{C}_1,e_{i}}\})|\geq
1$, $|(g_1\setminus g'_1)\cap
(f_j\setminus\{f_{\mathcal{C}_2,f_{j}},l_{\mathcal{C}_2,f_{j}}\})|\geq
1$. \item $|(\overline{g'}_1\setminus {g}_1)\cap
(e_i\setminus\{f_{\mathcal{C}_1,e_{i}},l_{\mathcal{C}_1,e_{i}}\})|\geq
1$, $|(\overline{g'}_1\setminus {g}_1)\cap
(f_j\setminus\{f_{\mathcal{C}_2,f_{j}},l_{\mathcal{C}_2,f_{j}}\})|\geq
1$. \item $|(g_1\setminus \overline{g'}_1)\cap
(e_i\setminus\{f_{\mathcal{C}_1,e_{i}},l_{\mathcal{C}_1,e_{i}}\})|\geq
1$, $|(g_1\setminus \overline{g'}_1)\cap
(f_j\setminus\{f_{\mathcal{C}_2,f_{j}},l_{\mathcal{C}_2,f_{j}}\})|\geq
1$.
 \item
 $|B'\cap (g_1\setminus g'_1)|=|B'\cap (g'_1\setminus g_1)|=|B'\cap (g_1\setminus \overline{g'}_1)|=|B'\cap (\overline{g'}_1\setminus g_1)|=1$.
\end{itemize}
\end{lemma}

\begin{lemma}\label{there is a P:3}
With the same assumptions in Remark {\rm \ref{assumption}}, let $l_1\geq 3$ and   $e_i$ and $f_j$ be two arbitrary edges of $\mathcal{C}_1$ and
$\mathcal{C}_2$, respectively. Also,  let
$\mathcal{E}_{1}=E_1=g_1g'_1$  be a red path of length $2$
 with the following properties.
\begin{itemize}
\item  $V(E_1)\subseteq (V(e_{i-1})\cup V(f_{j-1})\cup
B')\setminus\{\overline{u},\overline{v}\}$ where $\overline{v}\in
e_{i-1}\setminus\{f_{\mathcal{C}_1,e_{i-1}}\}$, $\overline{u}\in
f_{j-1}\setminus\{f_{\mathcal{C}_2,f_{j-1}}\}$ and $B'\subseteq W$
with $|B'|=2$.
 \item $|(g'_1\setminus g_1)\cap
(e_{i-1}\setminus\{f_{\mathcal{C}_1,e_{i-1}},l_{\mathcal{C}_1,e_{i-1}}\})|\geq
1$, $|(g'_1\setminus g_1)\cap
(f_{j-1}\setminus\{f_{\mathcal{C}_2,f_{j-1}},l_{\mathcal{C}_2,f_{j-1}}\})|\geq
1$. \item $({g'}_1\setminus {g}_1)\cap B'=\{w\}$.
\end{itemize}
 If there
is no blue copy of $\mathcal{C}_{l_1+l_2}$, then there are two red
paths $E_2=g_2g'_2$ and $F_2=\overline{g}_2\overline{g'}_2$ of
length $2$ so that the following conditions hold.
\begin{itemize}
\item  $V(E_2),V(F_2)\subseteq \Big((V(e_{i})\cup
V(f_{j}))\setminus\{f_{\mathcal{C}_1,e_{i}},f_{\mathcal{C}_2,f_{j}}\}\Big)\cup\{w,\hat{u},\hat{v}\}$
for some  $\hat{v}\in
(e_{i-1}\setminus\{f_{\mathcal{C}_1,e_{i-1}}\})\cap(V(E_1)\cup\{\overline{v}\})$,
$\hat{u}\in
(f_{j-1}\setminus\{f_{\mathcal{C}_2,f_{j-1}}\})\cap(V(E_1)\cup\{\overline{u}\})$.
\item $e_i\setminus(E_2\cup\{f_{\mathcal{C}_1,e_{i}}\})\neq
\emptyset$ and
$f_j\setminus(F_2\cup\{f_{\mathcal{C}_2,f_{j}}\})\neq\emptyset$.
\item $|(g'_2\setminus g_2)\cap
(e_{i}\setminus\{f_{\mathcal{C}_1,e_{i}},l_{\mathcal{C}_1,e_{i}}\})|\geq
1$, $|(\overline{g'}_2\setminus \overline{g}_2)\cap
(f_{j}\setminus\{f_{\mathcal{C}_2,f_{j}},l_{\mathcal{C}_2,f_{j}}\})|\geq
1$. \item $\mathcal{E}_{2}=\mathcal{E}_{1} E_2$ and
$\mathcal{F}_{2}=\mathcal{E}_{1}F_2$ are two red paths of length
$4$.
\end{itemize}
\end{lemma}

\medskip
\begin{lemma}\label{there is a P}
With the same assumptions in Remark {\rm \ref{assumption}}, let
$i\geq 2$ and
\begin{eqnarray*}\mathcal{E}_{i-1}=E_1E_2\ldots E_{i-1},
\mathcal{F}_{i-1}=F_1F_2\ldots F_{i-1}
\end{eqnarray*}
 be two red paths of
length $2i-2$ where for each $1\leq t\leq i-1$, $E_t=g_tg'_t$ and
$F_t=\overline{g}_t\overline{g'}_t$
 are red paths of length
$2$ with the following properties.
\begin{enumerate}
\item[$P_1$:] For each $1\leq t\leq i-1$,
\begin{eqnarray*}
V(E_t),V(F_t)\subseteq
(V(e_{t})\setminus \{f_{\mathcal{C}_1,e_{t}}\})\cup
(V(f_{t})\setminus\{f_{\mathcal{C}_2,f_{t}}\})\cup
B_t\cup\{\hat{v}_t,\hat{u}_t\}
\end{eqnarray*} where  $\hat{v}_t\in
e_{t-1}\setminus\{f_{\mathcal{C}_1,e_{t-1}}\}$,  $\hat{u}_t\in
f_{t-1}\setminus\{f_{\mathcal{C}_2,f_{t-1}}\}$, $|B_t|\leq 2$ for
$t=1$ and $|B_t\setminus \bigcup_{j=1}^{t-1}B_j|\leq 1$ for $t\neq
1$.
 \item[$P_2$:] For $t=i-1$
\begin{itemize} \item $e_{t}\setminus(E_{t}\cup\{f_{\mathcal{C}_1,e_{t}}\})\neq
\emptyset$ and
$f_{t}\setminus(F_{t}\cup\{f_{\mathcal{C}_2,f_{t}}\})\neq
\emptyset$. \item $|(g'_t\setminus g_t)\cap
(e_t\setminus\{f_{\mathcal{C}_1,e_{t}},l_{\mathcal{C}_1,e_{t}}\})|\geq
1$, $|(g'_t\setminus g_t)\cap
(f_t\setminus\{f_{\mathcal{C}_2,f_{t}},l_{\mathcal{C}_2,f_{t}}\})|\geq
1$. \item $|(\overline{g'}_t\setminus \overline{g}_t)\cap
(e_t\setminus\{f_{\mathcal{C}_1,e_{t}},l_{\mathcal{C}_1,e_{t}}\})|\geq
1$, $|(\overline{g'}_t\setminus \overline{g}_t)\cap
(f_t\setminus\{f_{\mathcal{C}_2,f_{t}},l_{\mathcal{C}_2,f_{t}}\})|\geq
1$.
\end{itemize}
\end{enumerate}
 If there
is no blue copy of $\mathcal{C}_{l_1+l_2}$ and $W\setminus \bigcup
_{t=1}^{i-1}B_t\neq \emptyset$, then there are two red paths
$E_i=g_ig'_i$ and $F_i=\overline{g}_i\overline{g'}_i$ of length
$2$ such that the properties $P_1$ and $P_2$ hold for $t=i$ and
for some $\mathcal{P}\in\{\mathcal{E}_{i-1},\mathcal{F}_{i-1}\}$,
$\mathcal{E}_i=\mathcal{P}E_i$ and $\mathcal{F}_i=\mathcal{P}F_i$
are two red paths of length $2i$.



\end{lemma}

The proof of the following Lemma is similar to the proof of Lemma
\ref{there is a P}.

\medskip
\begin{lemma}\label{there is a P:2}
Let  $\mathcal{E}_{i-1}=E_1E_2\ldots E_{i-1}$ and
$\mathcal{F}_{i-1}=F_1F_2\ldots F_{i-1}$ be two red paths of
length $2i-2$ with the same properties of Lemma {\rm\ref{there is a
P}}.
 Let
$v\in e_{i+1}\setminus \{l_{\mathcal{C}_1,e_{i+1}}\}$ and
$u\in f_{i+1}\setminus\{l_{\mathcal{C}_2,f_{i+1}}\}$ are two
distinct vertices so that $\{v,u\}\cap V(\mathcal{E}_{i-1}\cup
\mathcal{F}_{i-1})=\emptyset$. If there is no blue copy of
$\mathcal{C}_{l_1+l_2}$, then there is a red path
$E_i=g_ig'_i$  of length
$2$ such that the following conditions hold.
\begin{itemize}
\item
 $V(E_i)\subseteq
(V(e_{i})\setminus
\{f_{\mathcal{C}_1,e_{i}},l_{\mathcal{C}_1,e_{i}}\})\cup
(V(f_{i})\setminus\{f_{\mathcal{C}_2,f_{i}},l_{\mathcal{C}_2,f_{i}}\})\cup\{u,v,u',v'\}$
where  $v'\in e_{i-1}\setminus\{f_{\mathcal{C}_1,e_{i-1}}\}$ and
$u'\in f_{i-1}\setminus\{f_{\mathcal{C}_2,f_{i-1}}\}$.   \item
$|(g'_i\setminus g_i)\cap
(e_i\setminus\{f_{\mathcal{C}_1,e_{i}},l_{\mathcal{C}_1,e_{i}}\})|\geq
2$, $|(g'_i\setminus g_i)\cap
(f_i\setminus\{f_{\mathcal{C}_2,f_{i}},l_{\mathcal{C}_2,f_{i}}\})|\geq
2$. \item For some
$\mathcal{P}\in\{\mathcal{E}_{i-1},\mathcal{F}_{i-1}\}$,
$\mathcal{E}_i=\mathcal{P}E_i$ is a red path of length $2i$.

\end{itemize}

\end{lemma}

\medskip
\begin{lemma}\label{No  blue Cn implies red Cn+1/2}
Let $n\geq 5$ be odd, $k\geq 8$ and
$\mathcal{H}=\mathcal{K}^k_{(k-1)n+\lfloor\frac{n-1}{2}\rfloor}$
be $2$-edge colored red and blue. Let  $\mathcal{C}_1$ and
$\mathcal{C}_2$ be two disjoint blue cycles of length
$\frac{n-1}{2}$ and $\frac{n+1}{2}$, respectively. If there is no
blue copy of $\mathcal{C}_{n}^k$, then there is a copy of
$\mathcal{C}_{\frac{n+1}{2}}^k$ in $\mathcal{H}_{\rm red}.$
\end{lemma}

\medskip
\begin{lemma}\label{No  red Cn implies blue Cn/2+1:2}
Let $n\geq 6$ be even, $k\geq 8$ and
$\mathcal{H}=\mathcal{K}^k_{(k-1)n+\lfloor\frac{n-1}{2}\rfloor}$
be $2$-edge colored red and blue. Let  $\mathcal{C}_1$ and
$\mathcal{C}_2$ be two disjoint blue cycles of length
$\frac{n}{2}-1$ and $\frac{n}{2}+1$, respectively. If there is no
blue copy of $\mathcal{C}_{n}^k$, then there is a copy of
$\mathcal{C}_{\frac{n}{2}+1}^k$ in $\mathcal{H}_{\rm red}.$
\end{lemma}

\section{Proof of Theorem \ref{main theorem}}

 We prove the theorem by induction on
$n$.  By Theorem \ref{R(Pk3,Pk3)} the theorem holds when
$n\leq 4$. Let $f(n)=(k-1)n+\lfloor\frac{n-1}{2}\rfloor$ and suppose
to the contrary that $\mathcal{H}=\mathcal{K}^k_{f(n)}$ is 2-edge
colored red and blue with  no monochromatic copy of
$\mathcal{C}^k_n$. We consider the following cases.\\

\noindent{\bf Case 1: $n$ is odd}
\begin{emp}\label{red-blue C}
There are two monochromatic copies of
$\mathcal{C}^k_{\frac{n+1}{2}}$ of colors red and blue.
\end{emp}
\noindent{\bf Proof.} Since $f(\frac{n+1}{2})<f(n),$ using
induction hypothesis, there is a monochromatic copy of
$\mathcal{C}^k_{\frac{n+1}{2}}$, say $\mathcal{C}_1$. Because of the symmetry we may  assume that
$\mathcal{C}_1\subseteq \mathcal{H}_{\rm blue}.$ Since
$|V(\mathcal{C}_1)|=(k-1)(\frac{n+1}{2})$ and
\begin{eqnarray*}f(\frac{n-1}{2})=(k-1)(\frac{n-1}{2})+\lfloor\frac{n-3}{4}\rfloor<f(n)-(k-1)(\frac{n+1}{2}),
\end{eqnarray*}
there is a monochromatic
$\mathcal{C}_2=\mathcal{C}^k_{\frac{n-1}{2}}$ in
$V(\mathcal{H})\setminus V(\mathcal{C}_1).$ If
$\mathcal{C}_2\subseteq \mathcal{H}_{\rm blue}$,
 using Lemma
\ref{No  blue Cn implies red Cn+1/2}, there is  a copy of
$\mathcal{C}^k_{\frac{n+1}{2}}$ in $\mathcal{H}_{\rm red}$ and so
we are done. Now let $\mathcal{C}_2\subseteq \mathcal{H}_{\rm
red}$. If there is no $\mathcal{C}^k_{\frac{n-1}{2}}$ in
$\mathcal{H}_{\rm blue}$, then  using Lemma \ref{No  red Cn-1
implies blue Cn} there is a red copy of
$\mathcal{C}^k_{\frac{n+1}{2}}$ and again
 we are done. Therefore, we may assume that we
have a red copy of $\mathcal{C}^k_{\frac{n-1}{2}}$ and a blue copy
of $\mathcal{C}^k_{\frac{n-1}{2}}$.  Among all red-blue copies of
$\mathcal{C}^k_{\frac{n-1}{2}}$'s choose red-blue copies with
maximum intersection, say $\mathcal{C'}_1$ and $\mathcal{C'}_2$. Set $A=V(\mathcal{C'}_1)\cup V(\mathcal{C'}_2)$.
We can show that $|A|\leq
R(\mathcal{C}^k_{\frac{n-1}{2}},\mathcal{C}^k_{\frac{n-1}{2}})+1$.
To see that, suppose to the contrary that  $|A|\geq
R(\mathcal{C}^k_{\frac{n-1}{2}},\mathcal{C}^k_{\frac{n-1}{2}})+2$.
Since
\begin{eqnarray*}|A|\geq {\rm
max}\{|V(\mathcal{C'}_1)|,|V(\mathcal{C'}_2)|\}+2,
\end{eqnarray*}
 then
$B_1=V(\mathcal{C'}_1)\setminus V(\mathcal{C'}_2)$ and
$B_2=V(\mathcal{C'}_2)\setminus V(\mathcal{C'}_1)$ are non-empty.
Choose $v_1\in B_1$ and $v_2\in B_2$ and set
\begin{eqnarray*}
U=\Big(V(\mathcal{C'}_1)\cup V(\mathcal{C'}_2)\Big)\setminus
\{v_1,v_2\}.
\end{eqnarray*}
 Since $|U|\geq
R(\mathcal{C}_\frac{n-1}{2},\mathcal{C}_\frac{n-1}{2})$, we have a
monochromatic  copy of $\mathcal{C}_\frac{n-1}{2},$ say
$\mathcal{C}$. If $\mathcal{C}$ is red, $|\mathcal{C}\cap
\mathcal{C'}_2|>|\mathcal{C'}_1\cap \mathcal{C'}_2|$. If
$\mathcal{C}$ is blue, then $|\mathcal{C}\cap \mathcal{C'}_1|>
|\mathcal{C'}_1\cap \mathcal{C'}_2|$.
 Both cases contradict the choices of $\mathcal{C'}_1$ and $\mathcal{C'}_2$. So $|A|\leq
R(\mathcal{C}^k_{\frac{n-1}{2}},\mathcal{C}^k_{\frac{n-1}{2}})+1$.
 Clearly
\begin{eqnarray*}
f(\frac{n+1}{2})\leq f(n)-|A|.
 \end{eqnarray*}
  So  using
  induction hypothesis, we have a monochromatic $\mathcal{C}^k_{\frac{n+1}{2}}$,
 say $\mathcal{C}$, in $V(\mathcal{H})\setminus V(\mathcal{C'}_1\cup \mathcal{C'}_2)$. If $\mathcal{C}$
  is red, then
   $\mathcal{C}$ and  $\mathcal{C}_1$ are desired monochromatic copies of  $\mathcal{C}^k_{\frac{n+1}{2}}$.
   If not, since there is no blue $\mathcal{C}^k_n$, using Lemma
   \ref{No  blue Cn implies red Cn+1/2}, there is a copy of $\mathcal{C}^k_{\frac{n+1}{2}}$, say $\mathcal{C'}$, in $\mathcal{H}_{\rm
red}$. Clearly $\mathcal{C'}$ and $\mathcal{C}_1$ are our
favorable  cycles. $\square$\\\\

 Among all red-blue copies of
$\mathcal{C}^k_{\frac{n+1}{2}}$'s choose red-blue copies with
maximum intersection. Let $\mathcal{C}_1\subseteq \mathcal{H}_{\rm
red}$ and $\mathcal{C}_2\subseteq \mathcal{H}_{\rm blue}$ be such
copies. Assume that $\mathcal{C}_1=d_1d_2\ldots
d_{\frac{n+1}{2}}$. An argument similar to the proof of Claim \ref{red-blue C} yields
\begin{eqnarray*}|V(\mathcal{C}_1\cup \mathcal{C}_2)|\leq
R(\mathcal{C}^k_{\frac{n+1}{2}},\mathcal{C}^k_{\frac{n+1}{2}})+1.
\end{eqnarray*}
Let $B=V(\mathcal{H})\setminus V(\mathcal{C}_1\cup
\mathcal{C}_2)$.
 It is easy to show that
 $|B|\geq f(\frac{n-1}{2})$. So, by  induction hypothesis, there is a monochromatic  $\mathcal{C}^k_{\frac{n-1}{2}}$, say $\mathcal{C}_3$, in $B$.
  By symmetry we may assume that $\mathcal{C}_3\subseteq \mathcal{H}_{\rm
  blue}$. Let
    $\mathcal{C}_2=f_1f_2\ldots
  f_{\frac{n+1}{2}}$, $\mathcal{C}_3=e_1e_2\ldots e_{\frac{n-1}{2}}$ and
$W=V(\mathcal{H})\setminus V(\mathcal{C}_2\cup \mathcal{C}_3)$
where
$$f_i=\{u_1,u_2,\ldots,u_k\}+(k-1)(i-1) \Big({\rm mod}\  (k-1)(\frac{n+1}{2})\Big),\hspace{0.5 cm} i=1,2,\ldots,\frac{n+1}{2}$$
\noindent and
$$e_i=\{v_1,v_2,\ldots,v_k\}+(k-1)(i-1) \Big({\rm mod}\ (k-1)(\frac{n-1}{2})\Big),\hspace{0.5 cm} i=1,2,\ldots,\frac{n-1}{2}.$$
\begin{emp}\label{333}
Let  $f_j$ be an arbitrary edge of $\mathcal{C}_2$ and $z\in W$.
There are vertices $\overline{u}\in
f_{j-1}\setminus\{f_{\mathcal{C}_2,f_{j-1}}\}$, $\overline{u'}\in
f_{j+1}\setminus\{l_{\mathcal{C}_2,f_{j+1}}\}$ and
$\hat{u}\in
f_{j}\setminus\{f_{\mathcal{C}_2,f_{j}},l_{\mathcal{C}_2,f_{j}}\}$
so that  the edge
\begin{eqnarray*}
g=(f_{j}\setminus\{f_{\mathcal{C}_2,f_{j}},l_{\mathcal{C}_2,f_{j}},\hat{u}\})\cup\{\overline{u},\overline{u'},z\}
\end{eqnarray*}
is blue.
\end{emp}

\noindent{\bf Proof of Claim \ref{333}.} Suppose not.  We may assume that
$f_j=f_{\frac{n+1}{2}}$.
 We  find a red copy of
$\mathcal{C}^k_n$ as follows.\\

Use Lemma
 \ref{there is a P2} for $e=e_1$
(resp. $e=f_1$) (by putting
 $i=j=1$, $v'=v_1$, $v''=v_k$, $u'=u_{(k-1)(\frac{n+1}{2})}$,
$u''=u_k$, $C=\{v_{k-1}\}$ and  $B=\{z,w_1\}\subseteq W$)
 to obtain a red path
 $E_1=g_1g'_1$ (resp. $F_1=\overline{g}_1\overline{g'}_1$)   with the mentioned properties of Lemma
 \ref{there is a P2}. Let $\mathcal{E}_1=E_1$ and $\mathcal{F}_1=F_1$.  We may assume that $(g_1\setminus g'_1)\cap\{z\}=E_1\cap\{z\}$ (also  ($\overline{g}_1\setminus \overline{g'}_1)\cap\{z\}=F_1\cap\{z\}$).
 Use
Lemma \ref{there is a P}, $\frac{n-5}{2}$ times to obtain  two red
paths
$\mathcal{E}_{\frac{n-3}{2}}$ and
$\mathcal{F}_{\frac{n-3}{2}}$  of length $n-3$ with desired properties.\\

 Let $i=\frac{n-1}{2}$, $v=v_{k-1}$,
 $u=f_{\mathcal{C}_2,f_{\frac{n+1}{2}}}$ and use Lemma
 \ref{there is a P:2} to obtain a  red path
 $E_{\frac{n-1}{2}}=g_{\frac{n-1}{2}}g'_{\frac{n-1}{2}}$ so that for some
$\mathcal{P}\in\{\mathcal{E}_{\frac{n-3}{2}},\mathcal{F}_{\frac{n-3}{2}}\}$,
$\mathcal{E}_{\frac{n-1}{2}}=\mathcal{P}E_{\frac{n-1}{2}}$ is a red path
 of length $n-1$. Let
$\mathcal{E}_{\frac{n-1}{2}}=h_1h_2\ldots
h_{n-1}$, $\overline{u}\in
(f_{\frac{n-1}{2}}\setminus\{f_{\mathcal{C}_2,f_{\frac{n-1}{2}}}\})\cap(h_{n-1}\setminus
h_{n-2})$ and $\hat{u}=u_{(k-1)(\frac{n+1}{2})}$. If $z\in h_1$, then set $\overline{u'}=u_1$. Otherwise, let $\overline{u'}\in
(f_1\setminus \{l_{\mathcal{C}_2,f_1}\})\cap(h_1\setminus
h_2)$. Now set
\begin{eqnarray*} g=(f_{\frac{n+1}{2}}\setminus\{f_{\mathcal{C}_2,f_{\frac{n+1}{2}}},l_{\mathcal{C}_2,f_{\frac{n+1}{2}}},\hat{u}\})\cup\{\overline{u},\overline{u'},z\}.
\end{eqnarray*}
Clearly, $\mathcal{E}_{\frac{n-1}{2}}g$ is a red copy of
$\mathcal{C}^k_n$, a contradiction to our assumption. This contradiction completes the proof of Claim
\ref{333}. $\square$

\begin{emp}\label{11}
Let $e_i$ and $f_j$ be two arbitrary edges of $\mathcal{C}_3$ and
$\mathcal{C}_2$, respectively.
For every vertices  $\overline{u'}\in
f_j\setminus\{f_{\mathcal{C}_2,f_j}\}$ and $x\in f_{j-1}\setminus\{f_{\mathcal{C}_2,f_{j-1}}\}$ {\rm \Big(}resp. $\overline{u'}\in
f_j\setminus\{l_{\mathcal{C}_2,f_j}\}$ and $x\in f_{j+1}\setminus\{l_{\mathcal{C}_2,f_{j+1}}\}${\rm \Big)},
there are vertices
 $\overline{v}\in e_i\setminus
\{f_{\mathcal{C}_3,e_i},l_{\mathcal{C}_3,e_i}\}$
and  $\overline{u}\in f_{j-1}\setminus
\{f_{\mathcal{C}_2,f_{j-1}},x\}$ {\rm \Big(}resp. $\overline{u}\in
f_{j+1}\setminus \{l_{\mathcal{C}_2,f_{j+1}},x\}${\rm \Big)} so
that the edge $(V(f_j)\setminus
\{f_{\mathcal{C}_2,f_j},\overline{u'}\})\cup\{\overline{v},\overline{u}\}$ {\rm \Big(}resp. the edge $(V(f_j)\setminus
\{l_{\mathcal{C}_2,f_j},\overline{u'}\})\cup\{\overline{v},\overline{u}\}${\rm \Big)}
 is blue.

\end{emp}

\noindent{\bf Proof of Claim \ref{11}.}
 By symmetry  it only suffices to show that for every vertices  $\overline{u'}\in
f_j\setminus\{f_{\mathcal{C}_2,f_j}\}$ and $x\in f_{j-1}\setminus\{f_{\mathcal{C}_2,f_{j-1}}\}$, there are vertices  $\overline{v}\in e_i\setminus
\{f_{\mathcal{C}_3,e_i},l_{\mathcal{C}_3,e_i}\}$
and  $\overline{u}\in f_{j-1}\setminus
\{f_{\mathcal{C}_2,f_{j-1}},x\}$ so that the edge
\begin{eqnarray*}
(V(f_j)\setminus
\{f_{\mathcal{C}_2,f_j},\overline{u'}\})\cup\{\overline{v},\overline{u}\}
\end{eqnarray*}
 is blue. With no loss of generality we may assume that
  $e_i=e_1$ and $f_j=f_{\frac{n+1}{2}}$. Suppose to the
contrary that there are vertices  $\overline{u'}\in f_{\frac{n+1}{2}}\setminus\{f_{\mathcal{C}_2,f_{\frac{n+1}{2}}}\}$ and $x\in f_{\frac{n-1}{2}}\setminus\{f_{\mathcal{C}_2,f_{\frac{n-1}{2}}}\}$ such that
 for every  vertices $\overline{v}\in e_1\setminus
\{v_1,v_k\}$ and $\overline{u}\in f_{\frac{n-1}{2}}\setminus
\{f_{\mathcal{C}_2,f_\frac{n-1}{2}},x\}$ the edge
\begin{eqnarray*}(V(f_{\frac{n+1}{2}})\setminus\{f_{\mathcal{C}_2,f_{\frac{n+1}{2}}},\overline{u'}\})\cup\{\overline{v},\overline{u}\}
\end{eqnarray*}
is red.  We can find a red  copy of $\mathcal{C}^k_n$ as follows.\\

Use Lemma \ref{there is a P2} for $e=e_1$
(resp. $e=f_1$)
 to obtain a red path
 $E_1=g_1g'_1$ (resp. $F_1=\overline{g}_1\overline{g'}_1$)   with the mentioned properties of Lemma
 \ref{there is a P2} (by putting  $i=j=1$, $v'=v_1$, $v''=v_k$, $u'=\overline{u'}$,
$u''=u_k$, $C=\{v_{k-1}\}$ and  $B=\{w_1,w_2\}\subseteq W$). Let $\mathcal{E}_1=E_1$ and $\mathcal{F}_1=F_1$.  Use
Lemma \ref{there is a P}, $\frac{n-5}{2}$ times to obtain  two red
paths $\mathcal{E}_{\frac{n-3}{2}}$ and
$\mathcal{F}_{\frac{n-3}{2}}$ of length $n-3$. \\

 Let $i=\frac{n-1}{2}$, $v=v_{k-1}$,
 $u=f_{\mathcal{C}_2,f_{\frac{n+1}{2}}}$ and use Lemma
 \ref{there is a P:2} to obtain a  red path
 $E_{\frac{n-1}{2}}=g_{\frac{n-1}{2}}g'_{\frac{n-1}{2}}$
so that  for some
$\mathcal{P}\in\{\mathcal{E}_{\frac{n-3}{2}},\mathcal{F}_{\frac{n-3}{2}}\}$,
$\mathcal{E}_{\frac{n-1}{2}}=\mathcal{P}E_{\frac{n-1}{2}}$ is a
red path of length $n-1$.
Let
$\mathcal{E}_{\frac{n-1}{2}}=h_1h_2\ldots
h_{n-1}$, $\overline{v}\in (e_1\setminus \{v_1,v_k\})\cap
(h_1\setminus h_2)$, $\overline{u}\in
(f_{\frac{n-1}{2}}\setminus
\{f_{\mathcal{C}_2,f_{\frac{n-1}{2}}},x\})\cap (h_{{n-1}}\setminus
h_{n-2})$ and
\begin{eqnarray*}g=(V(f_{\frac{n+1}{2}})\setminus
\{f_{\mathcal{C}_2,f_{\frac{n+1}{2}}},\overline{u'}\})\cup\{\overline{v},\overline{u}\}.
\end{eqnarray*}
 Clearly, $\mathcal{E}_{\frac{n-1}{2}}g$ is a red copy of
$\mathcal{C}^k_n$, a contradiction to our assumption. This
contradiction completes the proof of Claim \ref{11}.
$\square$\\

\begin{emp}\label{12}
Let $e_i$ and $f_j$ be two arbitrary edges of $\mathcal{C}_3$ and
$\mathcal{C}_2$, respectively. If $n\geq 7$, then for every vertices $z\in f_j$ and
$v\in
(e_{i+1}\setminus\{l_{\mathcal{C}_3,e_{i+1}}\})\cup(e_{i-1}\setminus\{f_{\mathcal{C}_3,e_{i-1}}\})$
the edge
$(e_i\setminus\{f_{\mathcal{C}_3,e_i},l_{\mathcal{C}_3,e_i}\})\cup\{z,v\}$
is red. If $n=5$, then for every $z\in f_j$ and
$v\in
\{f_{\mathcal{C}_3,e_{i}},l_{\mathcal{C}_3,e_{i}}\}$
the edge
\begin{eqnarray*}(e_i\setminus\{f_{\mathcal{C}_3,e_i},l_{\mathcal{C}_3,e_i}\})\cup\{z,v\}
\end{eqnarray*}
is red.
\end{emp}

\noindent{\bf Proof of Claim \ref{12}.}
We give only a proof for  $n\geq 7$. The proof for $n=5$ is similar.
Suppose for a contradiction that there are
vertices $z\in f_j$ and $v\in
(e_{i+1}\setminus\{l_{\mathcal{C}_3,e_{i+1}}\})\cup(e_{i-1}\setminus\{f_{\mathcal{C}_3,e_{i-1}}\})$
so that the edge
\begin{eqnarray*}
g=(V(e_i)\setminus\{f_{\mathcal{C}_3,e_i},l_{\mathcal{C}_3,e_i}\})\cup\{z,v\}
\end{eqnarray*}
is blue. With no loss of generality assume that $v\in
e_{i+1}\setminus\{l_{\mathcal{C}_3,e_{i+1}}\}$ (by symmetry
the case $v\in  e_{i-1}\setminus\{f_{\mathcal{C}_3,e_{i-1}}\}$
is similar). If  $z=l_{\mathcal{C}_2,f_{j}}$, then put $\overline{u'}=z$ and  $x=l_{\mathcal{C}_2,f_{j-1}}$ and use  Claim
\ref{11} to obtain a blue edge
\begin{eqnarray*}
g'=(V(f_j)\setminus\{f_{\mathcal{C}_2,f_{j}},\overline{u'}\})\cup\{\overline{v},\overline{u}\}
\end{eqnarray*}
 for some
$\overline{v}\in e_{i-1}\setminus \{f_{\mathcal{C}_3,e_{i-1}},l_{\mathcal{C}_3,e_{i-1}}\}$ and
$\overline{u}\in f_{j-1}\setminus \{f_{\mathcal{C}_2,f_{j-1}},x\}$. Clearly
\begin{eqnarray*}
g'f_{j-1}f_{j-2}\ldots f_{1} f_{\frac{n+1}{2}}\ldots f_{j+1} g e_{i+1}e_{i+2}\ldots e_{\frac{n-1}{2}} e_1\ldots e_{i-1}\end{eqnarray*}
 is a blue $\mathcal{C}_n^k$, a contradiction to our assumption.
If  $z\in
f_j\setminus\{l_{\mathcal{C}_2,f_{j}}\}$,  then put $\overline{u'}=l_{\mathcal{C}_2,f_{j-1}}$ and $x=l_{\mathcal{C}_2,f_{j-2}}$ and use Claim \ref{11} to obtain a blue edge
\begin{eqnarray*}
g''=(V(f_{j-1})\setminus\{f_{\mathcal{C}_2,f_{j-1}},\overline{u'}\})\cup\{\overline{v},\overline{u}\}
\end{eqnarray*}
 for some
$\overline{v}\in e_{i-1}\setminus \{f_{\mathcal{C}_3,e_{i-1}},l_{\mathcal{C}_3,e_{i-1}}\}$ and
$\overline{u}\in f_{j-2}\setminus \{f_{\mathcal{C}_2,f_{j-2}},x\}$.
Again
\begin{eqnarray*}
g''f_{j-2}f_{j-3}\ldots f_{1} f_{\frac{n+1}{2}}\ldots f_{j} g e_{i+1}e_{i+2}\ldots e_{\frac{n-1}{2}} e_1\ldots e_{i-1}
\end{eqnarray*}
 is a
blue copy  of $\mathcal{C}^k_n$, a contradiction to our assumption. $\square$\\

\begin{emp}\label{13}
Let $e_i$  be an arbitrary edge of $\mathcal{C}_3$. If $n\geq 7$, then for every
vertices  $z\in W$ and $v\in
(e_{i-1}\setminus\{f_{\mathcal{C}_3,e_{i-1}}\})\cup(e_{i+1}\setminus\{l_{\mathcal{C}_3,e_{i+1}}\})$
the edge
$(e_i\setminus\{f_{\mathcal{C}_3,e_i},l_{\mathcal{C}_3,e_i}\})\cup\{z,v\}$
is red. If $n=5$, then for every vertices  $z\in W$ and
$v\in
\{f_{\mathcal{C}_3,e_{i}},l_{\mathcal{C}_3,e_{i}}\}$
the edge
\begin{eqnarray*}
(e_i\setminus\{f_{\mathcal{C}_3,e_i},l_{\mathcal{C}_3,e_i}\})\cup\{z,v\}
\end{eqnarray*}
is red.
\end{emp}

\noindent{\bf Proof of Claim \ref{13}.} We give only a proof for   $n\geq 7$.
The proof for $n=5$ is similar.  Suppose indirectly  that there are vertices
$z\in W$ and $v\in
(e_{i-1}\setminus\{f_{\mathcal{C}_3,e_{i-1}}\})\cup
(e_{i+1}\setminus\{l_{\mathcal{C}_3,e_{i+1}}\})$ so that the
edge
\begin{eqnarray*}
g=(e_i\setminus\{f_{\mathcal{C}_3,e_i},l_{\mathcal{C}_3,e_i}\})\cup\{z,v\}
\end{eqnarray*}
 is blue. With no loss of generality assume that $v\in
e_{i+1}\setminus\{l_{\mathcal{C}_3,e_{i+1}}\}$. Using
 Claim \ref{333} with $f_j=f_1$, there are vertices $\overline{u}\in
f_{\frac{n+1}{2}}\setminus\{f_{\mathcal{C}_2,f_{\frac{n+1}{2}}}\}$,
$\overline{u'}\in f_{2}\setminus\{l_{\mathcal{C}_2,f_{2}}\}$
and $\hat{u}\in
f_{1}\setminus\{f_{\mathcal{C}_2,f_{1}},l_{\mathcal{C}_2,f_{1}}\}$
so that  the edge
\begin{eqnarray*}
g'=(f_{1}\setminus\{f_{\mathcal{C}_2,f_{1}},l_{\mathcal{C}_2,f_{1}},\hat{u}\})\cup\{\overline{u},\overline{u'},z\}
\end{eqnarray*}
is blue. Use
Claim \ref{11}
to obtain a blue edge
\begin{eqnarray*}
g''=(f_{\frac{n+1}{2}}\setminus\{f_{\mathcal{C}_2,f_{\frac{n+1}{2}}},\overline{u}\})\cup\{\overline{x},\overline{y}\}
\end{eqnarray*}
 for some
 $\overline{x}\in e_{i-1}\setminus\{f_{\mathcal{C}_3,e_{i-1}},l_{\mathcal{C}_3,e_{i-1}}\}$ and $\overline{y}\in f_{\frac{n-1}{2}}\setminus\{f_{\mathcal{C}_2,f_{\frac{n-1}{2}}}\}$. Clearly
 \begin{eqnarray*}
 e_1e_2 \ldots e_{i-1}g'' f_{\frac{n-1}{2}}f_{\frac{n-3}{2}} \ldots f_2 g' g e_{i+1} e_{i+2} \ldots e_{\frac{n-1}{2}}
 \end{eqnarray*}
  is a
  blue copy of $\mathcal{C}^k_n$, a contradiction to our
assumption. This contradiction completes the proof of Claim
\ref{13}. $\square$\\

Now, set

\begin{eqnarray*}
h_1= \left\lbrace
\begin{array}{ll}
(e_1\setminus \{v_k\})\cup\{z_1\}  & n\equiv^41,\vspace{.5 cm}\\
(e_1\setminus\{v_1,v_k\})\cup\{z_1,v_{(k-1)(\frac{n-1}{2})}\} &
 n\equiv^43,
\end{array}
\right.\vspace{.2 cm}
\end{eqnarray*}
\noindent and for $2\leq i\leq \frac{n-1}{2}$,

\begin{eqnarray*}
h_i= \left\lbrace
\begin{array}{ll}
(e_i\setminus \{l_{\mathcal{C}_3,e_i}\})\cup\{z_{\frac{i+1}{2}}\}  &\mbox{if~}  i~\mbox{is~odd},\vspace{.5 cm}\\
(e_i\setminus\{f_{\mathcal{C}_3,e_i}\})\cup\{z_{\frac{i}{2}}\}
&\mbox{if~} i~\mbox{is~even},
\end{array}
\right.\vspace{.2 cm}
\end{eqnarray*}
\noindent where $\{z_1,z_2,\ldots,z_{\frac{n-1}{4}}\}\subseteq
V(\mathcal{H})\setminus V(\mathcal{C}_1\cup \mathcal{C}_3)$. Using
Claims \ref{12} and \ref{13}, $\mathcal{C}_4=h_1h_2\ldots
h_{\frac{n-1}{2}}$ is a red copy of
$\mathcal{C}^k_{\frac{n-1}{2}}$ disjoint from $\mathcal{C}_1$.
 The proof of the
following claim is similar to the proof of  Claim \ref{11}. So we
omit it here.





\begin{emp}\label{14}
Let $d_i$ and $h_j$ be two arbitrary edges of $\mathcal{C}_1$ and
$\mathcal{C}_4$, respectively. For every vertex $\overline{u'}\in
d_i\setminus\{f_{\mathcal{C}_1,d_i}\}$ {\rm (}resp. $\overline{u'}\in
d_i\setminus\{l_{\mathcal{C}_1,d_i}\}${\rm )},
there are vertices  $\overline{v}\in h_j\setminus
\{f_{\mathcal{C}_4,h_j},l_{\mathcal{C}_4,h_j}\}$
and  $\overline{u}\in d_{i-1}\setminus
\{f_{\mathcal{C}_1,d_{i-1}}\}$ {\rm (}resp. $\overline{u}\in
d_{i+1}\setminus \{l_{\mathcal{C}_1,d_{i+1}}\}${\rm )} so
that the edge
$V(d_i)\setminus
\{f_{\mathcal{C}_1,d_i},\overline{u'}\})\cup\{\overline{v},\overline{u}\}$
 {\rm (}resp. the edge $(V(d_i)\setminus
\{l_{\mathcal{C}_1,d_i},\overline{u'}\})\cup\{\overline{v},\overline{u}\}${\rm )}
 is red.
\end{emp}

 Now, using Claim \ref{14}
there are vertices $\overline{v}\in
h_1\setminus\{f_{\mathcal{C}_4,h_1},l_{\mathcal{C}_4,h_1}\}$
and $\overline{u}\in
d_{\frac{n+1}{2}}\setminus\{f_{\mathcal{C}_1,d_{\frac{n+1}{2}}}\}$
so that the edge
\begin{eqnarray*}
g=(V(d_1)\setminus\{f_{\mathcal{C}_1,d_{1}},l_{\mathcal{C}_1,d_{1}}\})\cup\{\overline{u},\overline{v}\}
\end{eqnarray*}
is red (by putting  $i=j=1$ and $\overline{u'}=l_{\mathcal{C}_1,d_{1}}$). Now let $z=l_{\mathcal{C}_1,d_{1}}$. If $n\neq 5$, then set  $y\in
(h_3\setminus\{f_{\mathcal{C}_4,h_3},l_{\mathcal{C}_4,h_3}\})\cap(e_3\setminus\{f_{\mathcal{C}_3,e_3},l_{\mathcal{C}_3,e_3}\})$
 and if $n=5$, then set  $y=f_{\mathcal{C}_3,e_1}$. When $z\in \mathcal{C}_2$, Claim \ref{12} and when $z\in
W$, Claim \ref{13} implies that the edge
\begin{eqnarray*}
g'=(h_2\setminus\{f_{\mathcal{C}_4,h_2},l_{\mathcal{C}_4,h_2}\})\cup\{z,y\}
\end{eqnarray*}
is red. Now, clearly
\begin{eqnarray*}
gd_{\frac{n+1}{2}}d_{\frac{n-1}{2}}\ldots
d_2 g' h_3h_4\ldots h_{\frac{n-1}{2}}h_1
\end{eqnarray*}
for $n\neq 5$ and $gd_3d_2g'h_1$ for $n=5$
 is a red copy of
$\mathcal{C}^k_n$, a contradiction.\\\\


\noindent{\bf Case 2: $n\equiv 0$ (mod 4)}\\

 In this case,  we show that
 there are two disjoint isochromatic
$\mathcal{C}^k_{\frac{n}{2}}$. Since
\begin{eqnarray*}
R(\mathcal{C}^k_{\frac{n}{2}},\mathcal{C}^k_{\frac{n}{2}})=f(\frac{n}{2})<f(n),
\end{eqnarray*}
there is a monochromatic
$\mathcal{C}_1=\mathcal{C}^k_{\frac{n}{2}}$. By symmetry we may
assume that $\mathcal{C}_1\subseteq \mathcal{H}_{blue}$. Since
$|V(\mathcal{C}_1)|=(k-1)(\frac{n}{2})$ and
\begin{eqnarray*}
f(\frac{n}{2})<f(n)-(k-1)(\frac{n}{2}),
\end{eqnarray*}
there is a monochromatic
$\mathcal{C}_2=\mathcal{C}^k_{\frac{n}{2}}$ in
$V(\mathcal{H})\setminus V(\mathcal{C}_1)$. If $\mathcal{C}_2$ is
blue, we are done. So suppose that $\mathcal{C}_2$ is red. Among
all red-blue copies of $\mathcal{C}^k_{\frac{n}{2}}$'s,  choose
red- blue copies with maximum intersection, say $\mathcal{C'}_1$
and $\mathcal{C'}_2$. Similar to the proof of Claim \ref{red-blue C} we have
\begin{eqnarray*}
|V(\mathcal{C'}_1)\cup
V(\mathcal{C'}_2)|\leq
R(\mathcal{C}^k_{\frac{n}{2}},\mathcal{C}^k_{\frac{n}{2}})+1=f(\frac{n}{2})+1.
\end{eqnarray*}
Since
\begin{eqnarray*}
f(n)-|V(\mathcal{C'}_1)\cup
V(\mathcal{C'}_2)|\geq f(\frac{n}{2})=R(\mathcal{C}^k_{\frac{n}{2}},\mathcal{C}^k_{\frac{n}{2}}),
\end{eqnarray*}
 there is a monochromatic $\mathcal{C}^k_{\frac{n}{2}}$, say
 $\mathcal{C}$,  disjoint from $\mathcal{C'}_1$ and
 $\mathcal{C'}_2$. Clearly $\mathcal{C}$ with one of $\mathcal{C'}_1$ and
 $\mathcal{C'}_2$ form our favorable cycles.\\

With no loss of generality assume that $\mathcal{C}_1$ and
$\mathcal{C}_2$ are two blue $\mathcal{C}^k_{\frac{n}{2}}$. Let
$\mathcal{C}_1=e_1e_2\ldots e_{\frac{n}{2}}$,
    $\mathcal{C}_2=f_1f_2\ldots
  f_{\frac{n}{2}}$  and
$W=V(\mathcal{H})\setminus V(\mathcal{C}_1\cup \mathcal{C}_2)$
where
\begin{eqnarray*}e_i=\{v_1,v_2,\ldots,v_k\}+(k-1)(i-1) ({\rm mod}\ (k-1)\frac{n}{2}),\hspace{0.5 cm} i=1,2,\ldots,\frac{n}{2}
\end{eqnarray*}
\noindent and
\begin{eqnarray*}
f_i=\{u_1,u_2,\ldots,u_k\}+(k-1)(i-1) ({\rm mod}\  (k-1)\frac{n}{2}),\hspace{0.5 cm} i=1,2,\ldots,\frac{n}{2}.
\end{eqnarray*}
Since $n\geq 8$, we have $|W|\geq 3$. Use Lemma \ref{there is a P2:3} for
$e_i=e_1$, $f_j=f_1$ and $B=\{w_1,w_2,w_3\}\subseteq W$ to obtain
two red paths $E_1$ and $F_1$ with desired properties in Lemma \ref{there is a P2:3}. As
mentioned in Lemma \ref{there is a P2:3}, there are distinct
vertices $\overline{v}\in e_1\setminus(V(E_1)\cup V(F_1))$, $v\in
e_1\setminus(V(E_1)\cup\{v_1,v_k,\overline{v}\})$,
$\overline{u}\in f_1\setminus(V(E_1)\cup V(F_1))$ and $u\in
f_1\setminus(V(F_1)\cup\{u_1,u_k,\overline{u}\})$. If
$\overline{v}\neq v_1$, then set $g_i=e_i$ for $1\leq i\leq
\frac{n}{2}$. If $\overline{v}=v_1$ then set $g_1=e_1$ and
$g_i=e_{\frac{n}{2}-i+2}$ for $2\leq i \leq \frac{n}{2}$. If
$\overline{u}\neq u_1$, then set $h_i=f_i$ for $1\leq i\leq
\frac{n}{2}$. If $\overline{u}=u_1$ then set $h_1=f_1$ and
$h_i=f_{\frac{n}{2}-i+2}$ for $2\leq i \leq \frac{n}{2}$. Clearly $\mathcal{C}_1=g_1g_2\ldots g_{\frac{n}{2}}$ and
$\mathcal{C}_2=h_1h_2\ldots h_{\frac{n}{2}}$. Assume that
\begin{eqnarray*}
&&g_i=\{x_1,x_2,\ldots,x_k\}+(k-1)(i-1) ({\rm mod}\ (k-1)\frac{n}{2}),\hspace{0.5 cm} i=1,2,\ldots,\frac{n}{2},\\
&&h_i=\{y_1,y_2,\ldots,y_k\}+(k-1)(i-1) ({\rm mod}\  (k-1)\frac{n}{2}),\hspace{0.5 cm} i=1,2,\ldots,\frac{n}{2}.
\end{eqnarray*}
 Now, use
Lemma \ref{there is a P:3} for $e_i=g_2$, $f_j=h_2$ and
$\mathcal{E}_1=E_1$ (resp. $\mathcal{E}_1=F_1$) to obtain two red
paths $\mathcal{E}_2$ and $\mathcal{F}_2$ (resp.
$\overline{\mathcal{E}}_2$ and $\overline{\mathcal{F}}_2$) of
length $4$ with desired properties of Lemma \ref{there is a
P:3}.
 Apply
Lemma \ref{there is a P} for $\mathcal{E}_l$, $\mathcal{F}_l$
(resp.  for $\overline{\mathcal{E}}_l$ and
$\overline{\mathcal{F}}_l$), $e_i=g_{l+1}$ and $f_j=h_{l+1}$ where
$2\leq l \leq \frac{n}{2}-2$ to obtain two red paths
$\mathcal{E}_{l+1}$ and $\mathcal{F}_{l+1}$ (resp.
$\overline{\mathcal{E}}_{l+1}$ and $\overline{\mathcal{F}}_{l+1}$)
with properties of Lemma \ref{there is a P}. Let
\begin{eqnarray*}
\mathcal{E}_{\frac{n}{2}-1}=p_1p_2\ldots p_{n-2},
\mathcal{F}_{\frac{n}{2}-1}=p_1p_2p'_3\ldots p'_{n-2},
\overline{\mathcal{E}}_{\frac{n}{2}-1}=p_1\overline{p}_2\ldots
\overline{p}_{n-2}
\end{eqnarray*}
 and
$\overline{\mathcal{F}}_{\frac{n}{2}-1}=p_1\overline{p}_2\overline{p'}_3\ldots
\overline{p'}_{n-2}$. Also, let $x\in
g_{\frac{n}{2}-1}\setminus(\mathcal{E}_{\frac{n}{2}-1}\cup\{f_{\mathcal{C}_1,g_{\frac{n}{2}-1}}\})$,
$y\in
h_{\frac{n}{2}-1}\setminus(\mathcal{F}_{\frac{n}{2}-1}\cup\{f_{\mathcal{C}_2,h_{\frac{n}{2}-1}}\})$,
$\overline{x}\in
g_{\frac{n}{2}-1}\setminus(\overline{\mathcal{E}}_{\frac{n}{2}-1}\cup\{f_{\mathcal{C}_1,g_{\frac{n}{2}-1}}\})$,
$\overline{y}\in
h_{\frac{n}{2}-1}\setminus(\overline{\mathcal{F}}_{\frac{n}{2}-1}\cup\{f_{\mathcal{C}_2,h_{\frac{n}{2}-1}}\})$,
$y''\in
 (h_{\frac{n}{2}-1}\setminus\{f_{\mathcal{C}_2,h_{\frac{n}{2}-1}}\})\cap(p_{n-2}\setminus
 p_{n-3})$ and $x''\in (g_{\frac{n}{2}-1}\setminus\{f_{\mathcal{C}_1,g_{\frac{n}{2}-1}}\})\cap(\overline{p'}_{n-2}\setminus \overline{p'}_{n-3})$.
  Consider an edge $q=E\cup F$ in
$\mathcal{A}_{\frac{n}{2}\frac{n}{2}}$ so that
\begin{eqnarray*}
&&E=\{x,x_{(k-1)(\frac{n}{2}-1)+2},\ldots,x_{(k-1)(\frac{n}{2}-1)+\lfloor\frac{k}{2}\rfloor}\},\\
&&F=\{y_{(k-1)\frac{n}{2}-\lceil\frac{k}{2}\rceil+2},\ldots,y_{(k-1)\frac{n}{2}},y'\},
\end{eqnarray*}
where $y'\in (h_1\setminus\{y_k\})\cap(p_1\setminus p_2)$.

\begin{emp}\label{g is red:6}
The edge $q$ is red.
\end{emp}
{\bf Proof of Claim \ref{g is red:6}.} Suppose indirectly that the
edge $q_1=q$ is blue. Since there is no blue copy of
$\mathcal{C}^k_{n}$, using Remark \ref{complementary edges are not
blue}, every  edge in $\mathcal{B}_{\frac{n}{2}\frac{n}{2}}$ that
is disjoint from $q_1$ is red.
 For $2\leq l\leq
\lfloor\frac{k}{2}\rfloor$, let
 $q_l=(q_{l-1}\setminus\{x_{(k-1)(\frac{n}{2}-1)+l}\})\cup\{x_{(k-1)\frac{n}{2}-l+2}\}$.
 Assume that
 $l'$ is the maximum $l\in[1,\lfloor\frac{k}{2}\rfloor]$ for
 which $q_l$ is blue. If $l'<\lfloor\frac{k}{2}\rfloor$, then $q_{l'+1}$ is red. Set
 \begin{eqnarray*}
 q'_{l'+1}=\Big((g_{\frac{n}{2}}\cup h_{\frac{n}{2}})\setminus(q_{l'}\cup
\{f_{\mathcal{C}_1,g_{\frac{n}{2}}},l_{\mathcal{C}_1,g_{\frac{n}{2}}},f_{\mathcal{C}_2,h_{\frac{n}{2}}},l_{\mathcal{C}_2,h_{\frac{n}{2}}}\})\Big)
\cup\{v,y''\}.
\end{eqnarray*}
  Clearly, $q'_{l'+1}$ is an edge in $\mathcal{B}_{\frac{n}{2}\frac{n}{2}}$ disjoint from $q_{l'}$. Since there is no blue copy of $\mathcal{C}_n^k$, the edge $q'_{l'+1}$ is red. Therefore,
   $\mathcal{E}_{\frac{n}{2}-1}q'_{l'+1}q_{l'+1}$ is a red
 copy of $\mathcal{C}^k_{n}$, a contradiction.
 Therefore, we may assume that $l'=\lfloor\frac{k}{2}\rfloor$ and
hence
the  edge $q_{\lfloor\frac{k}{2}\rfloor}=E' \cup F$ is
 blue, where
 \begin{eqnarray*}E'=\{x,x_{(k-1)\frac{n}{2}-\lfloor\frac{k}{2}\rfloor+2},\ldots,x_{(k-1)\frac{n}{2}-1},x_{(k-1)\frac{n}{2}}\}.
 \end{eqnarray*}

Now, set $m=\lfloor\frac{k}{2}\rfloor-1$.
For $1\leq l\leq m$, let
 \begin{eqnarray*}q_{\lfloor\frac{k}{2}\rfloor+l}=(q_{\lfloor\frac{k}{2}\rfloor+l-1}\setminus\{y_{(k-1)\frac{n}{2}-l+1}\})\cup\{y_{(k-1)(\frac{n}{2}-1)+l+1}\}.
 \end{eqnarray*} Now, let
 $l'$ be the maximum $l\in[0,m]$ for
 which $q_{\lfloor\frac{k}{2}\rfloor+l}$ is blue. If $l'<m$, then $q_{\lfloor\frac{k}{2}\rfloor+l'+1}$ is red. Set
 \begin{eqnarray*}q'_{\lfloor\frac{k}{2}\rfloor+l'+1}=\Big((g_{\frac{n}{2}}\cup h_{\frac{n}{2}})\setminus(q_{\lfloor\frac{k}{2}\rfloor+l'}\cup
\{f_{\mathcal{C}_1,g_{\frac{n}{2}}},l_{\mathcal{C}_1,g_{\frac{n}{2}}},f_{\mathcal{C}_2,h_{\frac{n}{2}}},l_{\mathcal{C}_2,h_{\frac{n}{2}}}\})\Big)\cup
\{v,y''\}.
\end{eqnarray*}
 Since there is no blue copy of $\mathcal{C}^k_{n}$, the edge
 $q'_{\lfloor\frac{k}{2}\rfloor+l'+1}$ is red. So
 \begin{eqnarray*}
 \mathcal{E}_{\frac{n}{2}-1}q'_{\lfloor\frac{k}{2}\rfloor+l'+1}q_{\lfloor\frac{k}{2}\rfloor+l'+1}
 \end{eqnarray*}
  is a red
 copy of $\mathcal{C}^k_{n}$, a contradiction.
  So we may
assume that $l'=m$ and
 hence the
 edge $q_{\lfloor\frac{k}{2}\rfloor+m}=E'\cup F'$ is
 blue, where
\begin{eqnarray*}
F'=\left\lbrace\begin{array}{ll}
\{y_{(k-1)(\frac{n}{2}-1)+2},\ldots,y_{(k-1)(\frac{n}{2}-1)+m+1},y'\}
&\mbox{if $k$ is even},\vspace{.5
cm}\\
\{y_{(k-1)(\frac{n}{2}-1)+2},\ldots,y_{(k-1)(\frac{n}{2}-1)+m+1},y_{(k-1)(\frac{n}{2}-1)+m+2},y'\}
&\mbox{if $k$ is odd}. \end{array} \right.\vspace{.2
cm}\end{eqnarray*}\\

Let
\begin{eqnarray*}
q_{\lfloor\frac{k}{2}\rfloor+m+1}=(q_{\lfloor\frac{k}{2}\rfloor+m}\setminus\{x\})\cup\{v\}.
\end{eqnarray*}
 If $q_{\lfloor\frac{k}{2}\rfloor+m+1}$ is red, then set
\begin{eqnarray*}
&&\hspace{-0.5
cm}q'_{\lfloor\frac{k}{2}\rfloor+m+1}=\Big((g_{\frac{n}{2}}\cup
h_{\frac{n}{2}})\setminus(q_{\lfloor\frac{k}{2}\rfloor+m}\cup
\{f_{\mathcal{C}_1,g_{\frac{n}{2}}},l_{\mathcal{C}_1,g_{\frac{n}{2}}},f_{\mathcal{C}_2,h_{\frac{n}{2}}},
l_{\mathcal{C}_2,h_{\frac{n}{2}}}\})\Big)\cup\{v,y''\}.
\end{eqnarray*} Since there is no blue copy of $\mathcal{C}_n^k$, the edge  $q'_{\lfloor\frac{k}{2}\rfloor+m+1}$ is red
 and
 \begin{eqnarray*}\mathcal{E}_{\frac{n}{2}-1}q'_{\lfloor\frac{k}{2}\rfloor+m+1}q_{\lfloor\frac{k}{2}\rfloor+m+1}
 \end{eqnarray*}
  is a red
 copy of $\mathcal{C}^k_{n}$, a contradiction.
 So we may assume that the edge
$q_{\lfloor\frac{k}{2}\rfloor+m+1}$ is blue.\\

 Now, let
 \begin{eqnarray*}
q_{\lfloor\frac{k}{2}\rfloor+m+2}=(q_{\lfloor\frac{k}{2}\rfloor+m+1}\setminus\{y',v\})\cup\{\overline{y},x'\}  \end{eqnarray*} where  $x'\in (g_1\setminus\{x_k\})\cap(p_1\setminus
\overline{p}_2).$

If $q_{\lfloor\frac{k}{2}\rfloor+m+2}$ is red, then set
\begin{eqnarray*}
q'_{\lfloor\frac{k}{2}\rfloor+m+2}=\Big((g_{\frac{n}{2}}\cup
h_{\frac{n}{2}})\setminus(q_{\lfloor\frac{k}{2}\rfloor+m+1}\cup
\{f_{\mathcal{C}_1,g_{\frac{n}{2}}},l_{\mathcal{C}_1,g_{\frac{n}{2}}},f_{\mathcal{C}_2,h_{\frac{n}{2}}},l_{\mathcal{C}_2,h_{\frac{n}{2}}}\})\Big)
\cup\{\overline{y},x''\}.
\end{eqnarray*}
\noindent
  Since there is no blue copy of $\mathcal{C}^k_{n}$, the edge
 $q'_{\lfloor\frac{k}{2}\rfloor+m+2}$ is red
 and
 \begin{eqnarray*} \overline{\mathcal{F}}_{\frac{n}{2}-1}q'_{\lfloor\frac{k}{2}\rfloor+m+2}q_{\lfloor\frac{k}{2}\rfloor+m+2}
   \end{eqnarray*}
is a red
 copy of $\mathcal{C}^k_{n}$, a contradiction.
   So we may assume that the edge
 $q_{\lfloor\frac{k}{2}\rfloor+m+2}$ is blue.\\

If $k$ is even, then clearly $q_{\lfloor\frac{k}{2}\rfloor+m+2}$
is an edge in $\mathcal{B}_{\frac{n}{2}\frac{n}{2}}$ disjoint from
$q_1$. This is impossible, by Remark \ref{complementary edges are
not blue}.
 Now we may assume that $k$ is odd.
One can easily see that
$x_{(k-1)(\frac{n}{2}-1)+\frac{k+1}{2}}\notin q_1\cup
q_{\lfloor\frac{k}{2}\rfloor+m+2}$ and
$y_{(k-1)(\frac{n}{2}-1)+\frac{k+1}{2}}\in q_1\cap
q_{\lfloor\frac{k}{2}\rfloor+m+2}$. Similarly,
 we can show
that the edge
$$q_{\lfloor\frac{k}{2}\rfloor+m+3}=q_{\lfloor\frac{k}{2}\rfloor+m+2}\setminus
\{y_{(k-1)(\frac{n}{2}-1)+\frac{k+1}{2}}\})\cup\{x_{(k-1)(\frac{n}{2}-1)+\frac{k+1}{2}}\}$$
is blue. That is a contradiction to Remark \ref{complementary
edges are not blue}. This contradiction completes the proof of
Claim \ref{g is red:6}. $\square$\\\\

Now, consider an edge $q'=E'\cup F'$ in
$\mathcal{B}_{\frac{n}{2}\frac{n}{2}}$ so that
\begin{eqnarray*}
&&E'=\{x_{(k-1)(\frac{n}{2}-1)+\lfloor\frac{k}{2}\rfloor},x_{(k-1)(\frac{n}{2}-1)+\lfloor\frac{k}{2}\rfloor+2},x_{(k-1)(\frac{n}{2}-1)+\lfloor\frac{k}{2}\rfloor+3}
 ,\ldots,x_{(k-1)\frac{n}{2}},v\},\\
&&F'=\{y'',y_{(k-1)(\frac{n}{2}-1)+2},\ldots,y_{(k-1)(\frac{n}{2}-1)+\lfloor\frac{k}{2}\rfloor}\}.
\end{eqnarray*}
 By an argument similar to the proof of Claim  \ref{g is red:6} we can show
 that $q'$ is red. Clearly $\mathcal{E}_{\frac{n}{2}-1}q'q$ is
 a red copy of $\mathcal{C}^k_{n}$, a contradiction to our
 assumption.\\\\

\noindent{\bf Case 3: $n\equiv 2$ (mod 4)}\\

By an argument similar to the proof of Claim \ref{red-blue C} we have the following claim.
\begin{emp}\label{red-blue C_{(n/2)+1}}
There are two monochromatic copies of
$\mathcal{C}^k_{\frac{n}{2}+1}$ of colors red and blue.
\end{emp}

 Among all red-blue copies of
$\mathcal{C}^k_{\frac{n}{2}+1}$'s choose red-blue copies with
maximum intersection. Let $\mathcal{C}_1=d_1d_2\ldots
d_{\frac{n}{2}+1}\subseteq \mathcal{H}_{\rm red}$ and
$\mathcal{C}_2\subseteq \mathcal{H}_{\rm blue}$ be such copies. It
is easy to see that $|V(\mathcal{C}_1\cup \mathcal{C}_2)|\leq
R(\mathcal{C}^k_{\frac{n}{2}+1},\mathcal{C}^k_{\frac{n}{2}+1})+1$.
 Since

   \begin{eqnarray*}|V(\mathcal{H})\setminus V(\mathcal{C}_1\cup
\mathcal{C}_2)|\geq f(\frac{n}{2}-1),
   \end{eqnarray*} using induction hypothesis, there is a monochromatic  $\mathcal{C}^k_{\frac{n}{2}-1}$  in $V(\mathcal{H})\setminus V(\mathcal{C}_1\cup
\mathcal{C}_2)$, say $\mathcal{C}_3$.
  By symmetry we may assume that $\mathcal{C}_3\subseteq \mathcal{H}_{\rm
  blue}$. Let
    $\mathcal{C}_2=f_1f_2\ldots
  f_{\frac{n}{2}+1}$, $\mathcal{C}_3=e_1e_2\ldots e_{\frac{n}{2}-1}$ and
$W=V(\mathcal{H})\setminus V(\mathcal{C}_2\cup \mathcal{C}_3)$
where
$$f_i=\{u_1,u_2,\ldots,u_k\}+(k-1)(i-1) \Big({\rm mod}\  (k-1)(\frac{n}{2}+1)\Big),\hspace{0.5 cm} i=1,2,\ldots,\frac{n}{2}+1$$
\noindent and
$$e_i=\{v_1,v_2,\ldots,v_k\}+(k-1)(i-1) \Big({\rm mod}\ (k-1)(\frac{n}{2}-1)\Big),\hspace{0.5 cm} i=1,2,\ldots,\frac{n}{2}-1.$$
\noindent  We have the following Claims.

\begin{emp}\label{22}
Let $e_i$ and $f_j$ be two arbitrary edges of $\mathcal{C}_3$ and
$\mathcal{C}_2$, respectively and $y\in e_i$. For every vertices $x\in f_j\setminus \{f_{\mathcal{C}_2,f_j}\}$ and $x'\in
f_{j+2}\setminus \{l_{\mathcal{C}_2,f_{j+2}}\}$, there are  distinct vertices
$v,v'\in e_i\setminus
\{f_{\mathcal{C}_3,e_i},l_{\mathcal{C}_3,e_i},y\}$,  $u\in
f_{j-1}\setminus \{f_{\mathcal{C}_2,f_{j-1}}\}$ and $u'\in
f_{j+3}\setminus \{l_{\mathcal{C}_2,f_{j+3}}\}$ so that at least one of
 the edges
$g=(V(f_j)\setminus \{f_{\mathcal{C}_2,f_j},x\})\cup\{v,u\}$ or
  \begin{eqnarray*}
g'=(V(f_{j+2})\setminus
\{x',l_{\mathcal{C}_2,f_{j+2}}\})\cup\{v',u'\}
  \end{eqnarray*} is blue.
\end{emp}

\noindent{\bf Proof of Claim \ref{22}.} Suppose to the contrary that there are vertices $x\in f_j\setminus \{f_{\mathcal{C}_2,f_j}\}$ and $x'\in
f_{j+2}\setminus \{l_{\mathcal{C}_2,f_{j+2}}\}$ so that
for every distinct vertices
$v,v'\in e_i\setminus
\{f_{\mathcal{C}_3,e_i},l_{\mathcal{C}_3,e_i},y\}$,  $u\in
f_{j-1}\setminus \{f_{\mathcal{C}_2,f_{j-1}}\}$ and $u'\in
f_{j+3}\setminus \{l_{\mathcal{C}_2,f_{j+3}}\}$
 the edges
   \begin{eqnarray*}
g=(V(f_j)\setminus \{f_{\mathcal{C}_2,f_j},x\})\cup\{v,u\}
  \end{eqnarray*} and
  \begin{eqnarray*}
g'=(V(f_{j+2})\setminus
\{x',l_{\mathcal{C}_2,f_{j+2}}\})\cup\{v',u'\}
  \end{eqnarray*} are red.
By symmetry we may assume that
$e_i=e_{\frac{n}{2}-1}$ and $f_j=f_{\frac{n}{2}-1}$.

 Use Lemma \ref{there is a P2} for $e=e_1$  (resp.
$e=f_1$) to obtain a red path $\mathcal{E}_1=g_1g'_1$ (resp.
$\mathcal{F}_1=\overline{g}_1\overline{g'}_1$) with the mentioned properties
of Lemma \ref{there is a P2} (by putting $i=j=1$, $v'=v_1$,
$v''=v_k$, $u'=u_1$, $u''=u_k$, $C=\{v_{k-1}\}$ and
$B=\{w_1,w_2\}\subseteq W$). Use Lemma \ref{there is a P},
$\frac{n}{2}-3$ times, to obtain two red paths
$\mathcal{E}_{\frac{n}{2}-2}$ and $\mathcal{F}_{\frac{n}{2}-2}$ of
length $n-4$ with the properties of Lemma \ref{there is a P}.

With no loss of generality assume that $y'\in
e_{\frac{n}{2}-2}\setminus(\mathcal{E}_{\frac{n}{2}-2}\cup\{f_{\mathcal{C}_3,e_{\frac{n}{2}-2}},v_{k-1}\})$.
Use Lemma \ref{there is a P2:2} for $i=\frac{n}{2}-1$,
$j=\frac{n}{2}$, $v'=y'$, $v''=v_{k-1}$, $u'=x$, $u''=x'$ to obtain
a red path $\mathcal{P}=g_{\frac{n}{2}-1}g'_{\frac{n}{2}-1}$ with
the properties of Lemma \ref{there is a P2:2}.

Let $\mathcal{E}_{\frac{n}{2}-2}=h_1 h_2 \ldots h_{n-4}$,
$\mathcal{P}=h_{n-3}h_{n-2}$,  $v\in
(e_{\frac{n}{2}-1}\setminus\{f_{\mathcal{C}_3,e_{\frac{n}{2}-1}},v_1,y\})\cap
(h_{n-3}\setminus h_{n-2})$, $v'\in
(e_{\frac{n}{2}-1}\setminus\{f_{\mathcal{C}_3,e_{\frac{n}{2}-1}},v_1,y\})\cap
(h_{n-2}\setminus h_{n-3})$, $u\in
(f_{\frac{n}{2}-2}\setminus\{f_{\mathcal{C}_2,f_{\frac{n}{2}-2}}\})\cap
(h_{n-4}\setminus h_{n-5})$ and $u'\in
(f_{1}\setminus\{l_{\mathcal{C}_2,f_{1}}\})\cap (h_1\setminus
h_2)$. Clearly,
$\mathcal{E}_{\frac{n}{2}-2}g\mathcal{P}g'$ is a
red copy of $\mathcal{C}_n^k$. This contradiction completes the proof of Claim \ref{22}.
$\square$



\begin{emp}\label{21}
Let $e_i$ and $f_j$ be two arbitrary edges of $\mathcal{C}_3$ and
$\mathcal{C}_2$, respectively. Also, let $A=(e_{i+1}\setminus\{l_{\mathcal{C}_3,e_{i+1}}\})\cup(e_{i-1}\setminus\{f_{\mathcal{C}_3,e_{i-1}}\})$.
For every vertices $z\in
f_j\setminus
\{f_{\mathcal{C}_2,f_j},l_{\mathcal{C}_2,f_j}\}$ and $\overline{v}\in
A$ the edge
$(e_i\setminus\{f_{\mathcal{C}_3,e_i},l_{\mathcal{C}_3,e_i}\})\cup\{z,\overline{v}\}$
is red.
\end{emp}

\noindent{\bf Proof of Claim \ref{21}.}  By symmetry it only  suffices to show that for every vertices $z\in
f_j\setminus
\{f_{\mathcal{C}_2,f_j},l_{\mathcal{C}_2,f_j}\}$ and $\overline{v}\in
e_{i+1}\setminus\{l_{\mathcal{C}_3,e_{i+1}}\}$,  the edge
  \begin{eqnarray*}
(e_i\setminus\{f_{\mathcal{C}_3,e_i},l_{\mathcal{C}_3,e_i}\})\cup\{z,\overline{v}\}
  \end{eqnarray*}is red.
With no loss of generality  we may assume that $e_i=e_{1}$
and $f_j=f_{\frac{n}{2}}$. Suppose for the sake of contradiction that there are
vertices  $z\in f_{\frac{n}{2}}\setminus
\{f_{\mathcal{C}_2,f_{\frac{n}{2}}},l_{\mathcal{C}_2,f_{\frac{n}{2}}}\}$
and  $\overline{v}\in
e_{2}\setminus\{l_{\mathcal{C}_3,e_{2}}\}$ so that the edge
 \begin{eqnarray*}
h=(e_1\setminus\{v_1,v_k\})\cup\{z,\overline{v}\}
  \end{eqnarray*} is blue. Since there is no blue copy of $\mathcal{C}^k_n$, then  for every
distinct vertices $v,v'\in
e_{\frac{n}{2}-1}\setminus\{f_{\mathcal{C}_3,e_{\frac{n}{2}-1}},l_{\mathcal{C}_3,e_{\frac{n}{2}-1}}, \overline{v}\}$,
$u\in
f_{\frac{n}{2}-2}\setminus\{f_{\mathcal{C}_2,f_{\frac{n}{2}-2}}\}$
and $u'\in f_{1}\setminus\{l_{\mathcal{C}_2,f_{1}}\}$ the
edges
\begin{eqnarray*}
g=(f_{\frac{n}{2}-1}\setminus\{f_{\mathcal{C}_2,f_{\frac{n}{2}-1}},l_{\mathcal{C}_2,f_{\frac{n}{2}-1}}\})\cup\{v,u\}
\end{eqnarray*} and
\begin{eqnarray*}
g'=(f_{\frac{n}{2}+1}\setminus\{f_{\mathcal{C}_2,f_{\frac{n}{2}+1}},l_{\mathcal{C}_2,f_{\frac{n}{2}+1}}\})\cup\{v',u'\}\end{eqnarray*}
 are red. That is a contradiction to  Claim \ref{22}. $\square$

 \begin{emp}\label{211}
Let   $e_i$ and $f_j$ be two arbitrary edges of $\mathcal{C}_3$ and
$\mathcal{C}_2$, respectively and $\overline{u}\in \{f_{\mathcal{C}_2,f_j},l_{\mathcal{C}_2,f_j}\}$. For $n\geq 10$ and $\overline{v}\in
e_{i+1}\setminus\{l_{\mathcal{C}_3,e_{i+1}}\}$ {\rm (}also for $n=6$ and $\overline{v}=l_{\mathcal{C}_3,e_{i}})$, if the edge
 \begin{eqnarray*}
(e_i\setminus\{f_{\mathcal{C}_3,e_i},l_{\mathcal{C}_3,e_i}\})\cup\{\overline{u},\overline{v}\}
\end{eqnarray*}
 is blue, then for every vertices $\hat{v}\in e_{i+1}\setminus\{l_{\mathcal{C}_3,e_{i+1}}\}$ and  $\hat{u}\in \{f_{\mathcal{C}_2,f_j},l_{\mathcal{C}_2,f_j}\}\setminus\{\overline{u}\}$ the edge
$
(e_i\setminus\{f_{\mathcal{C}_3,e_i},l_{\mathcal{C}_3,e_i}\})\cup\{\hat{u},\hat{v}\}$
 is red. For $n\geq 10$ and $\overline{v}\in
e_{i-1}\setminus\{f_{\mathcal{C}_3,e_{i-1}}\}$ {\rm (}also for $n=6$ and $\overline{v}=f_{\mathcal{C}_3,e_{i}})$, if the edge
\begin{eqnarray*}
(e_i\setminus\{f_{\mathcal{C}_3,e_i},l_{\mathcal{C}_3,e_i}\})\cup\{\overline{u},\overline{v}\}
\end{eqnarray*}
  is blue, then for every vertices $\hat{v}\in e_{i-1}\setminus\{f_{\mathcal{C}_3,e_{i-1}}\}$ and  $\hat{u}\in \{f_{\mathcal{C}_2,f_j},l_{\mathcal{C}_2,f_j}\}\setminus\{\overline{u}\}$ the edge
$(e_i\setminus\{f_{\mathcal{C}_3,e_i},l_{\mathcal{C}_3,e_i}\})\cup\{\hat{u},\hat{v}\}$
 is red.
\end{emp}

\noindent{\bf Proof of Claim \ref{211}.} We give only a proof for  $n\geq 10$. The proof for $n=6$ is similar. By symmetry we may assume that
 $\overline{v}\in
e_{i+1}\setminus\{l_{\mathcal{C}_3,e_{i+1}}\}$ and $\overline{u}=f_{\mathcal{C}_2,f_j}$.
With no loss of generality  we may assume that $e_i=e_{1}$
and $f_j=f_{\frac{n}{2}}$. Assume for the sake of contradiction that  there is a
 vertex $\hat{v}\in
e_{2}\setminus\{l_{\mathcal{C}_3,e_{2}}\}$ so that the edge
\begin{eqnarray*} (e_1\setminus\{v_1,v_k\})\cup\{l_{\mathcal{C}_2,f_{\frac{n}{2}}},\hat{v}\}
\end{eqnarray*}
 is blue. Since there is no blue copy of $\mathcal{C}^k_n$, then  for every
distinct vertices $v,v'\in
e_{\frac{n}{2}-1}\setminus\{f_{\mathcal{C}_3,e_{\frac{n}{2}-1}},l_{\mathcal{C}_3,e_{\frac{n}{2}-1}}\}$,
$u\in
f_{\frac{n}{2}-2}\setminus\{f_{\mathcal{C}_2,f_{\frac{n}{2}-2}}\}$
and $u'\in f_{1}\setminus\{l_{\mathcal{C}_2,f_{1}}\}$ the
edges
\begin{eqnarray*}
g=(f_{\frac{n}{2}-1}\setminus\{f_{\mathcal{C}_2,f_{\frac{n}{2}-1}},l_{\mathcal{C}_2,f_{\frac{n}{2}-1}}\})\cup\{v,u\}
\end{eqnarray*}
and
 \begin{eqnarray*}
g'=(f_{\frac{n}{2}+1}\setminus\{f_{\mathcal{C}_2,f_{\frac{n}{2}+1}},l_{\mathcal{C}_2,f_{\frac{n}{2}+1}}\})\cup\{v',u'\} \end{eqnarray*}
 are red,  a contradiction to  Claim \ref{22}. $\square$

\begin{emp}\label{23}
Let $e_i$  be an arbitrary edge of $\mathcal{C}_3$, $z\in W$ and
 $\overline{v}\in
(e_{i+1}\setminus\{l_{\mathcal{C}_3,e_{i+1}}\})\cup(e_{i-1}\setminus\{f_{\mathcal{C}_3,e_{i-1}}\})$.  If  the
edge
\begin{eqnarray*}
h=(e_i\setminus\{f_{\mathcal{C}_3,e_i},l_{\mathcal{C}_3,e_i}\})\cup\{z,\overline{v}\}
\end{eqnarray*}
is blue, then for every edge $f_j\in \mathcal{C}_2$ and every
distinct  vertices $x\in
f_{j-1}\setminus\{f_{\mathcal{C}_2,f_{j-1}}\}$, $x'\in
f_{j+1}\setminus\{l_{\mathcal{C}_2,f_{j+1}}\}$ and
$\hat{u}\in
f_{j}\setminus\{f_{\mathcal{C}_2,f_{j}},l_{\mathcal{C}_2,f_{j}}\}$
the edge
\begin{eqnarray*}
h'=(f_{j}\setminus\{f_{\mathcal{C}_2,f_{j}},l_{\mathcal{C}_2,f_{j}},\hat{u}\})\cup\{x,x',z\}
\end{eqnarray*} is red.
\end{emp}

\noindent{\bf Proof of Claim \ref{23}.}
By symmetry we may assume that $\overline{v}\in e_{i+1}\setminus\{l_{\mathcal{C}_3,e_{i+1}}\}$.
Suppose to the contrary that there is an
edge $f_j$ and there are distinct vertices $x\in
f_{j-1}\setminus\{f_{\mathcal{C}_2,f_{j-1}}\}$, $x'\in
f_{j+1}\setminus\{l_{\mathcal{C}_2,f_{j+1}}\}$ and
$\hat{u}\in
f_{j}\setminus\{f_{\mathcal{C}_2,f_{j}},l_{\mathcal{C}_2,f_{j}}\}$
so that  the edge
\begin{eqnarray*} h'=(f_{j}\setminus\{f_{\mathcal{C}_2,f_{j}},l_{\mathcal{C}_2,f_{j}},\hat{u}\})\cup\{x,x',z\}
\end{eqnarray*}
is blue. Since there is no blue copy of $\mathcal{C}_n^k$, then   for every distinct vertices $v,v'\in
e_{i-1}\setminus\{f_{\mathcal{C}_3,e_{i-1}},\overline{v}\}$, $u\in
f_{j-2}\setminus\{f_{\mathcal{C}_2,f_{j-2}}\}$ and $u'\in
f_{j+2}\setminus\{l_{\mathcal{C}_2,f_{j+2}}\}$ the edges
$(f_{j-1}\setminus\{f_{\mathcal{C}_2,f_{j-1}},x\})\cup\{v,u\}$ and
$(f_{j+1}\setminus\{x',l_{\mathcal{C}_2,f_{j+1}}\})\cup\{v',u'\}$
 are red, a contradiction to  Claim \ref{22}. So we are done. $\square$

\begin{emp}\label{24}
Let $n=6$, $W=\{z_1,z_2\}$ and $e_i$  be an arbitrary edge of $\mathcal{C}_3$ and $v\in
\{f_{\mathcal{C}_3,e_{i}},l_{\mathcal{C}_3,e_{i}}\}$. If the edge
\begin{eqnarray*}
(e_i\setminus\{f_{\mathcal{C}_3,e_i},l_{\mathcal{C}_3,e_i}\})\cup\{z_1,v\}
\end{eqnarray*}
is blue, then for every vertex $\overline{v}\in
(e_{i+1}\setminus\{f_{\mathcal{C}_3,e_{i+1}},l_{\mathcal{C}_3,e_{i+1}}\})\cup\{v\}$
the edge
\begin{eqnarray*}
(e_i\setminus\{f_{\mathcal{C}_3,e_i},l_{\mathcal{C}_3,e_i}\})\cup\{z_2,\overline{v}\}
\end{eqnarray*}
is red.
\end{emp}

\noindent{\bf Proof of Claim \ref{24}.} By symmetry we may assume that $e_i=e_1$ and $v=v_k$.
Suppose for the sake of contradiction that there is a vertex $\overline{v}\in e_2\setminus \{v_1\}$ so that
 the edges
$(e_1\setminus\{v_1\})\cup\{z_1\}$
and
$(e_1\setminus\{v_1,v_k\})\cup\{z_2,\overline{v}\}$
are blue. Let $u\in
f_{4}\setminus\{f_{\mathcal{C}_2,f_{4}},l_{\mathcal{C}_2,f_{4}}\}$
and  $\overline{u}\in f_{2}\setminus\{f_{\mathcal{C}_2,f_{2}},l_{\mathcal{C}_2,f_{2}}\}$.
 Use Lemma \ref{there is a P2} for $e=e_1$ to obtain a red path $\mathcal{P}=g_1g'_1$ with the mentioned properties of Lemma \ref{there is a P2} (by putting $i=j=1$, $v'=v_1$, $v''=v_k$, $u'=u$, $u''=\overline{u}$,  $C=\{v_{k-1}\}$ and $B=W=\{z_1,z_2\}$). By Lemma \ref{there is a P2}, there is a vertex $w\in e_1\setminus\{v_1,v_{k-1}\}$ so that $w\notin \mathcal{P}$. Now, use Lemma \ref{there is a P2:2} for $e_i=e_2$, $f_j=f_3$, $v'=w$, $v''=v_{k-1}$, $u'=f_{\mathcal{C}_2,f_{3}}$ and $u''=l_{\mathcal{C}_2,f_{3}}$ to obtain a red path $\mathcal{P}'=g_2g'_2$ with the mentioned properties of Lemma \ref{there is a P2:2}.
  Let $x\in
(f_{3}\setminus\{l_{\mathcal{C}_2,f_{3}}\})\cap
(g_{2}\setminus g'_{2})$, $x'\in
(f_{3}\setminus\{f_{\mathcal{C}_2,f_{3}}\})\cap
(g'_{2}\setminus g_{2})$  and $y\in
(f_{1}\setminus\{u_1,u_k\})\cap (g'_{1}\setminus g_{1})$.
 As mentioned in Lemma \ref{there is a P2}, we  may assume that $z_1\in g_1$.
 Now, let
  \begin{eqnarray*}
 h=(f_{4}\setminus\{u,f_{\mathcal{C}_2,f_{4}}\})\cup\{z_1,x'\}.
  \end{eqnarray*} Set $h'=(f_2\setminus\{\overline{u},l_{\mathcal{C}_2,f_2}\})\cup\{z_2,x\}$ if $z_2\in g'_1$ and $h'=(f_2\setminus\{\overline{u},f_{\mathcal{C}_2,f_2},l_{\mathcal{C}_2,f_2}\})\cup\{y,z_2,x\}$, otherwise. By Claim \ref{23}, the edges $h$ and $h'$ are red. So
$\mathcal{P}h'\mathcal{P}' h$
is a red copy of $\mathcal{C}^k_6$, a contradiction.
 $\square$

\begin{emp}\label{25}
Let $e_i$  be an arbitrary edge of $\mathcal{C}_3$ and $n\geq 10$. For every
vertices  $z\in W$ and $\hat{v}\in
(e_{i+1}\setminus\{l_{\mathcal{C}_3,e_{i+1}}\})\cup(e_{i-1}\setminus\{f_{\mathcal{C}_3,e_{i-1}}\})$
the edge
\begin{eqnarray*}
(e_i\setminus\{f_{\mathcal{C}_3,e_i},l_{\mathcal{C}_3,e_i}\})\cup\{z,\hat{v}\}
 \end{eqnarray*}
is red.
\end{emp}

\noindent{\bf Proof of Claim \ref{25}.}   Suppose indirectly  that there are vertices
$z\in W$ and $\hat{v}\in
(e_{i+1}\setminus\{l_{\mathcal{C}_3,e_{i+1}}\})\cup(e_{i-1}\setminus\{f_{\mathcal{C}_3,e_{i-1}}\})$ so that the
edge
\begin{eqnarray*}
g=(e_i\setminus\{f_{\mathcal{C}_3,e_i},l_{\mathcal{C}_3,e_i}\})\cup\{z,\hat{v}\}
\end{eqnarray*}
 is blue. We may assume that $i=1$ and    $\hat{v}\in e_{2}\setminus\{l_{\mathcal{C}_3,e_{2}}\}$. In the rest of the proof, we consider $W\setminus\{z\}$ instead of $W$ when we use Lemmas \ref{there is a P2:3}, \ref{there is a P:3} and \ref{there is a P}.

 Since $n\geq 10$, we have $|W|\geq 3$. Use Lemma \ref{there is a P2:3} for $e_i=e_1$,  $f_j=f_1$ and  $B=\{w_1,w_2,w_3\}\subseteq W$ to obtain a red path $E_1=g_1g'_1$ with the mentioned properties of Lemma \ref{there is a P2:3}. Let $B'=B\cap E_1$. As mentioned in  Lemma \ref{there is a P2:3}, there are distinct vertices $\overline{v}\in e_1\setminus E_1$, $\overline{u}\in f_1\setminus E_1$ and $\overline{v'}\in e_1\setminus(E_1\cup\{v_1,v_k,\overline{v}\})$. Set $\tilde{v}=\overline{v'}$, if $\overline{v}=v_1$ and $\tilde{v}=\overline{v}$, otherwise. Clearly $\tilde{v}\in e_1\setminus(E_1\cup \{v_1\})$. If $\overline{u}\neq u_1$, then set $g_i=f_i$ for $1\leq i \leq \frac{n}{2}+1$. If $\overline{u}=u_1$, then set $g_1=f_1$ and $g_t=f_{\frac{n}{2}-t+3}$ for $2\leq t \leq \frac{n}{2}+1$. Clearly $\mathcal{C}_2=g_1g_2\ldots g_{\frac{n}{2}+1}$. With no loss of generality assume that
 $$g_i=\{w_1,w_2,\ldots,w_k\}+(k-1)(i-1) \Big({\rm mod}\  (k-1)(\frac{n}{2}+1)\Big),\hspace{0.5 cm} i=1,2,\ldots,\frac{n}{2}+1.$$
 Use Lemma \ref{there is a P:3} for $e_i=e_2$ and $f_j=g_2$ and $\mathcal{E}_1=E_1$ to obtain
  two red paths $\mathcal{E}_2$ and $\mathcal{F}_2$ of length $4$ with the mentioned properties of Lemma \ref{there is a P:3}.
Now, use Lemma \ref{there is a P}, $\frac{n}{2}-4$ times (for $e_i=e_l$ and $f_j=g_l$ where $3\leq l \leq \frac{n}{2}-2$), to obtain two red paths $\mathcal{E}_{\frac{n}{2}-2}$ and $\mathcal{F}_{\frac{n}{2}-2}$ of length $n-4$ with the mentioned properties of Lemma \ref{there is a P}. Now, use Lemma \ref{there is a P:2} for $v\in\{\overline{v},\overline{v'}\}\setminus\{\tilde{v}\}$ and $u=l_{\mathcal{C}_2,g_{\frac{n}{2}-1}}$ to obtain a red path $\mathcal{E}_{\frac{n}{2}-1}$ of length $n-2$. With no loss of generality assume that
  $\mathcal{E}_{\frac{n}{2}-1}=h_1h_2 \ldots h_{n-2}$.
  Let $x\in (g_1\setminus\{l_{\mathcal{C}_2,g_1}\})\cap(h_1\setminus h_2)$ and $y\in  (g_{\frac{n}{2}-1}\setminus\{f_{\mathcal{C}_2,g_{\frac{n}{2}-1}}\})\cap(h_{n-2}\setminus h_{n-3})$. By Claim \ref{23} the edges
  \begin{eqnarray*}&&g=(g_{\frac{n}{2}}\setminus\{f_{\mathcal{C}_2,g_{\frac{n}{2}}},l_{\mathcal{C}_2,g_{\frac{n}{2}}},w_{(k-1)\frac{n}{2}}\})
  \{y,z,w_{(k-1)(\frac{n}{2}+1)}\},\\
&& g'=(g_{\frac{n}{2}+1}\setminus\{w_{(k-1)(\frac{n}{2}+1)},l_{\mathcal{C}_2,g_{\frac{n}{2}+1}}\})\cup\{x,z\}
  \end{eqnarray*} are red. So $\mathcal{E}_{\frac{n}{2}-1} g g'$ is a red copy of $\mathcal{C}_n^k$, a contradiction to our assumption.
 $\square$\\

\begin{emp}\label{30}
Let  $e_i$ and $f_j$ be two arbitrary edges of $\mathcal{C}_3$ and
$\mathcal{C}_2$, respectively. Choose $x\in\{f_{\mathcal{C}_2,f_{j-1}},l_{\mathcal{C}_2,f_j}\}$. For $n=6$, assume that
\begin{eqnarray*}(\overline{v},\hat{v})\in\Big(\{l_{\mathcal{C}_3,e_i}\}\times(e_{i+1}\setminus\{l_{\mathcal{C}_3,e_{i+1}}\})\Big)\cup
\Big(\{f_{\mathcal{C}_3,e_i}\}\times(e_{i-1}\setminus\{f_{\mathcal{C}_3,e_{i-1}}\})\Big).\end{eqnarray*}
 For $n> 6$, let $\overline{v},\hat{v}\in
e_{i+1}\setminus\{l_{\mathcal{C}_3,e_{i+1}}\}$ or $\overline{v},\hat{v}\in
e_{i-1}\setminus\{f_{\mathcal{C}_3,e_{i-1}}\})$. Then, at least one of the edges $(e_i\setminus\{f_{\mathcal{C}_3,e_i},l_{\mathcal{C}_3,e_i}\})\cup\{x,\overline{v}\}$ or $\Big((e_i\cup\{f_{\mathcal{C}_2,f_{j-1}},l_{\mathcal{C}_2,f_j}\})\setminus\{f_{\mathcal{C}_3,e_{i}},l_{\mathcal{C}_3,e_{i}},x\}\Big)\cup\{\hat{v}\}$ is red.
\end{emp}
\noindent{\bf Proof of Claim \ref{30}.} By symmetry we may assume that $e_i=e_1$, $f_j=f_{\frac{n}{2}+1}$ and $x=l_{\mathcal{C}_2,f_{\frac{n}{2}+1}}$.
Suppose indirectly that there
 are vertices $\overline{v},\hat{v}$ with the mentioned  properties so that the edges
  $(e_1\setminus\{f_{\mathcal{C}_3,e_1},
l_{\mathcal{C}_3,e_1}\})\cup\{l_{\mathcal{C}_2,f_{\frac{n}{2}+1}},\overline{v}\}$ and
 $(e_1\setminus \{f_{\mathcal{C}_3,e_{1}},l_{\mathcal{C}_3,e_{1}}\})\cup\{\hat{v},f_{\mathcal{C}_2,f_{\frac{n}{2}}}\}$ are blue.
With no loss of generality assume that $\overline{v}=l_{\mathcal{C}_3,e_1}$ and $\hat{v}\in
e_{2}\setminus\{l_{\mathcal{C}_3,e_{2}}\}$ for $n=6$ and $\overline{v},\hat{v}\in
e_{2}\setminus\{l_{\mathcal{C}_3,e_2}\}$, otherwise.
Use Lemma \ref{there is a P2} for $e=e_1$  (resp.
$e=f_1$)
 to obtain a red path $\mathcal{E}_1=g_1g'_1$ (resp. $\mathcal{F}_1=\overline{g}_1\overline{g'}_1$) with the mentioned properties
of Lemma \ref{there is a P2} (by putting $i=j=1$, $v'=v_1$,
$v''=v_k$, $u'=u_1$, $u''=u_k$, $C=\{v_{k-1}\}$ and
$B=\{w_1,w_2\}\subseteq W$). Use Lemma \ref{there is a P},
$\frac{n}{2}-3$ times, to obtain two red paths
$\mathcal{E}_{\frac{n}{2}-2}$ and $\mathcal{F}_{\frac{n}{2}-2}$ of
length $n-4$ with the properties of Lemma \ref{there is a P}.

Set  $y\in
e_{\frac{n}{2}-2}\setminus(\mathcal{E}_{\frac{n}{2}-2}\cup\{f_{\mathcal{C}_3,e_{\frac{n}{2}-2}},v_{k-1}\})$.
Use Lemma \ref{there is a P2:2} for $i=\frac{n}{2}-1$,
$j=\frac{n}{2}$, $v'=y$, $v''=v_{k-1}$, $u'=f_{\mathcal{C}_2,f_{\frac{n}{2}}}$, $u''=u_{\frac{n}{2}(k-1)+2}$ to obtain
a red path $\mathcal{P}=g_{\frac{n}{2}-1}g'_{\frac{n}{2}-1}$ with
the properties of Lemma \ref{there is a P2:2}.

Let $\mathcal{E}_{\frac{n}{2}-2}=h_1 h_2 \ldots h_{n-4}$,
$\mathcal{P}=h_{n-3}h_{n-2}$,
  $x'\in
(e_{\frac{n}{2}-1}\setminus\{f_{\mathcal{C}_3,e_{\frac{n}{2}-1}},l_{\mathcal{C}_3,e_{\frac{n}{2}-1}},\hat{v}\})\cap
(h_{n-2}\setminus h_{n-3})$, $x''\in
(e_{\frac{n}{2}-1}\setminus\{f_{\mathcal{C}_3,e_{\frac{n}{2}-1}},l_{\mathcal{C}_3,e_{\frac{n}{2}-1}},\hat{v}\})\cap
(h_{n-3}\setminus h_{n-2})$ and $y'\in
(f_{\frac{n}{2}-2}\setminus\{f_{\mathcal{C}_2,f_{\frac{n}{2}-2}}\})\cap
(h_{n-4}\setminus h_{n-5})$. By Lemma \ref{there is a P2}, we may
assume that $w_1\in (h_1\setminus h_2)\cap W$. Since there is no
blue copy of $\mathcal{C}_n^k$, the edges
\begin{eqnarray*}
h=(f_{\frac{n}{2}+1}\setminus\{l_{\mathcal{C}_2,f_{\frac{n}{2}+1}},u_{\frac{n}{2}(k-1)+2}\})\cup\{x',w_1\}
\end{eqnarray*}
and
\begin{eqnarray*}
h'=(f_{\frac{n}{2}-1}\setminus\{f_{\mathcal{C}_2,f_{\frac{n}{2}-1}},l_{\mathcal{C}_2,f_{\frac{n}{2}-1}}\})\cup\{x'',y'\}
\end{eqnarray*}
are red. Therefore, $\mathcal{E}_{\frac{n}{2}-2}h'\mathcal{P}h$ is
a red copy of $\mathcal{C}_n^k$. This contradiction completes the
proof of Claim \ref{30}.
$\square$\\

\begin{emp}\label{31}
Let $n=6$, $W=\{z_1,z_2\}$ and   $e_i$ and $f_j$ be two arbitrary
edges of $\mathcal{C}_3$ and $\mathcal{C}_2$, respectively. Let $A=\{a,b\}$ with $f_{\mathcal{C}_2,f_j}\in A$ and $|A\cap W|=1$. Assume that
\begin{eqnarray*}(\overline{v},\hat{v})\in \Big(\{l_{\mathcal{C}_3,e_i}\}\times (e_{i+1}\setminus\{l_{\mathcal{C}_3,e_{i+1}}\})\Big)\cup \Big(\{f_{\mathcal{C}_3,e_i}\}\times(e_{i-1}\setminus\{f_{\mathcal{C}_3,e_{i-1}}\})\Big).
\end{eqnarray*}
 At least one of the edges $(e_i\setminus\{f_{\mathcal{C}_3,e_i},l_{\mathcal{C}_3,e_i}\})\cup\{a,\overline{v}\}$
or
 $(e_i\setminus \{f_{\mathcal{C}_3,e_{i}},l_{\mathcal{C}_3,e_{i}}\})\cup\{b,\hat{v}\}$ is red.
\end{emp}
\noindent{\bf Proof of Claim \ref{31}.}  By symmetry we may assume
that $e_i=e_1$,  $f_j=f_{4}$, $a=f_{\mathcal{C}_2,f_4}$, $b=z_1$. Suppose indirectly that there are vertices $\overline{v}$  and $\hat{v}$ with the desired properties  so that
the edges
 $(e_1\setminus\{f_{\mathcal{C}_3,e_1},l_{\mathcal{C}_3,e_1}\})\cup\{a,\overline{v}\}$ and
 $(e_1\setminus \{f_{\mathcal{C}_3,e_{1}},l_{\mathcal{C}_3,e_{1}}\})\cup\{\hat{v},b\}$ are
 blue. With no loss of generality we may assume that $\overline{v}=l_{\mathcal{C}_3,e_1}=v_k$ and $\hat{v}\in e_{2}\setminus\{l_{\mathcal{C}_3,e_{2}}\}$.
Use Lemma \ref{there is a P2} for $e=e_1$
 to obtain a red path $\mathcal{E}_1=g_1g'_1$  with the mentioned properties
of Lemma \ref{there is a P2} (by putting $i=j=1$, $v'=v_1$,
$v''=v_k$, $u'=u_1$, $u''=u_{k+1}$, $C=\{v_{k-1}\}$ and
$B=W=\{z_1,z_2\}$). By symmetry we may assume that either  $z_1\in
g_1$ or $z_1\notin \mathcal{E}_1$.
As mentioned in  Lemma \ref{there is a P2}, there is a vertex   $y\in
e_{1}\setminus(\mathcal{E}_{1}\cup\{v_1,v_{k-1}\})$.

Use Lemma
\ref{there is a P2:2} for $i=2$, $j=3$, $v'=y$, $v''=v_{k-1}$,
$u'=f_{\mathcal{C}_2,f_{3}}$, $u''=l_{\mathcal{C}_2,f_{3}}$ to
obtain a red path $\mathcal{P}=g_{2}g'_{2}$ with the properties of
Lemma \ref{there is a P2:2}.

Let
  $x\in
(e_{2}\setminus\{v_k,v_1,\hat{v}\})\cap
(g'_2\setminus g_2)$,   $\overline{y}\in (f_{1}\setminus\{u_1,u_k\})\cap
(g_1\setminus g'_1)$, $y'\in (f_{1}\setminus\{u_1,u_k\})\cap
(g'_1\setminus g_1)$ and $y''\in
(f_3\setminus\{f_{\mathcal{C}_2,f_{3}},l_{\mathcal{C}_2,f_{3}}\})\cap(g_2\setminus
g'_2)$. Since there is no blue copy of $\mathcal{C}_n^k$, the edge
$h=(f_{4}\setminus\{f_{\mathcal{C}_2,f_{4}},l_{\mathcal{C}_2,f_{4}}\})\cup\{x,y'\}$
is red. If $z_1\in g_1$, then set \begin{eqnarray*}h'=(f_{2}\setminus\{u_{k+1},l_{\mathcal{C}_2,f_{2}}\})\cup\{z_1,y''\}.
\end{eqnarray*}
Otherwise, set
\begin{eqnarray*}
h'=(f_{2}\setminus\{f_{\mathcal{C}_2,f_{2}},l_{\mathcal{C}_2,f_{2}},u_{k+1}\})\cup\{z_1,\overline{y},y''\}.
\end{eqnarray*}
By Claim \ref{23}, the edge $h'$ is red. So, $h
\mathcal{E}_1h'\mathcal{P}$ is a red copy of  $\mathcal{C}_6^k$.
This contradiction completes the proof of Claim \ref{31}.
$\square$\\

\begin{emp}\label{28} Let  $n\geq 10$.
There is a red copy of $\mathcal{C}_{\frac{n}{2}-1}^k$, say
$\mathcal{C}_4=h_1h_2\ldots h_{\frac{n}{2}-1}$, disjoint from
$\mathcal{C}_1$ so that for each $1\leq i\leq \frac{n}{2}-1$,
$ k-2\leq|h_i\cap e_{i}|\leq k-1$. Moreover, $|h_{\frac{n}{2}-1}\cap
e_{\frac{n}{2}-1}|=k-1$ and
  \begin{eqnarray*}h_{\frac{n}{2}-1}\setminus\{f_{\mathcal{C}_4,h_{\frac{n}{2}-1}}\}=e_{\frac{n}{2}-1}\setminus\{f_{\mathcal{C}_3,e_{\frac{n}{2}-1}}\}.
      \end{eqnarray*}
\end{emp}

\noindent{\bf Proof of Claim \ref{28}.} Let $n=4l+2$ and
$W'=V(\mathcal{H})\setminus V(\mathcal{C}_1\cup \mathcal{C}_2\cup
\mathcal{C}_3)$.
  Since $|V(\mathcal{C}_1\cup \mathcal{C}_2)|\leq (k-1)(\frac{n}{2}+1)+l+1$, clearly $|W'|\geq l-1$. First let $|W'|\geq l$.
 For $1\leq i\leq \frac{n}{2}-1$, set
\begin{eqnarray*}
h_i= \left\lbrace
\begin{array}{ll}
(e_i\setminus \{l_{\mathcal{C}_3,e_i}\})\cup\{z_{\frac{i+1}{2}}\}  &\mbox{if~}  i~\mbox{is~odd},\vspace{.5 cm}\\
(e_i\setminus\{f_{\mathcal{C}_3,e_i}\})\cup\{z_{\frac{i}{2}}\}
&\mbox{if~} i~\mbox{is~even},
\end{array}
\right.\vspace{.2 cm}
\end{eqnarray*}
\noindent where $\{z_1,z_2,\ldots,z_{l}\}\subseteq W'$. By Claim
\ref{25}, $\mathcal{C}_4=h_1h_2\ldots h_{\frac{n}{2}-1}$ is a red
copy of $\mathcal{C}^k_{\frac{n}{2}-1}$. Now, let
$W'=\{z_1,z_2,\ldots,z_{l-1}\}$. This is the case only when
$|V(\mathcal{C}_1\cup \mathcal{C}_2)|=
(k-1)(\frac{n}{2}+1)+l+1$. Clearly, $|V(\mathcal{C}_2)\setminus
V(\mathcal{C}_1)|=l+1$. If there is a vertex  $x\in
f_j\setminus(\{f_{\mathcal{C}_2,f_j},l_{\mathcal{C}_2,f_j}\}\cup V(\mathcal{C}_1))$, for some edge $f_j$ of $\mathcal{C}_2$, then
the same argument yields a
 red $\mathcal{C}_4=h_1h_2\ldots
h_{\frac{n}{2}-1}$ (consider $W'\cup \{x\}$ instead of $W'$ in the
above argument and   use Claims
 \ref{21} and \ref{25}). Therefore, we may assume that each vertex of $V(\mathcal{C}_2)\setminus V(\mathcal{C}_1)$ is a first vertex
  of $f_j$ for some edge  $f_j$ of $\mathcal{C}_2$.

If there is an edge $f_j$ of $\mathcal{C}_2$ so that $\{f_{\mathcal{C}_2,f_j},l_{\mathcal{C}_2,f_j}\}\subseteq V(\mathcal{C}_2)\setminus V(\mathcal{C}_1)$,
 then do the following process.
Using Claim \ref{211}, there is a vertex $z\in\{f_{\mathcal{C}_2,f_j},l_{\mathcal{C}_2,f_j}\}$ so that   the edge $h_1=(e_1\setminus\{v_k\})\cup\{z\}$
 is red. If the edge $(e_2\setminus\{l_{\mathcal{C}_3,e_2}\})\cup\{z\}$ is red, then set $h_2=(e_2\setminus\{l_{\mathcal{C}_3,e_2}\})\cup\{z\}$.
  Otherwise, set $h_2=(e_2\setminus\{f_{\mathcal{C}_3,e_2},l_{\mathcal{C}_3,e_2}\})\cup\{v_{k-1},z'\}$,
  where $z'\in \{f_{\mathcal{C}_2,f_j},l_{\mathcal{C}_2,f_j}\}\setminus \{z\}$. By Claim \ref{211}, $h_2$ is red.
 For $3\leq i\leq \frac{n}{2}-1$, set
\begin{eqnarray*}
h_i= \left\lbrace
\begin{array}{ll}
(e_i\setminus \{f_{\mathcal{C}_3,e_i},l_{\mathcal{C}_3,e_i}\})\cup\{v_{(i-1)(k-1)},z_{\frac{i-1}{2}}\}  &\mbox{if~}  i~\mbox{is~odd},\vspace{.5 cm}\\
(e_i\setminus\{f_{\mathcal{C}_3,e_i}\})\cup\{z_{\frac{i}{2}-1}\}
&\mbox{if~} i~\mbox{is~even}.
\end{array}
\right.\vspace{.2 cm}
\end{eqnarray*} \noindent  By Claim \ref{25}, $h_i$'s, $3\leq i\leq
\frac{n}{2}-1$, are red and clearly $\mathcal{C}_4=h_1h_2\ldots
h_{\frac{n}{2}-1}$ is a red copy of
$\mathcal{C}^k_{\frac{n}{2}-1}$. If each the above cases does not
occur, then there is an edge $f_j$ so that
$\{f_{\mathcal{C}_2,f_{j-1}},l_{\mathcal{C}_2,f_{j}}\}\subseteq
V(\mathcal{C}_2)\setminus V(\mathcal{C}_1)$. Similar to the above
cases, by Claim \ref{30}, there is a vertex
$z\in \{f_{\mathcal{C}_2,f_{j-1}},l_{\mathcal{C}_2,f_{j}}\}$
so that the edge $h_1=(e_1\setminus\{v_k\})\cup\{z\}$
 is red. If the edge $(e_2\setminus\{l_{\mathcal{C}_3,e_2}\})\cup\{z\}$ is red, then set $h_2=(e_2\setminus\{l_{\mathcal{C}_3,e_2}\})\cup\{z\}$.
  Otherwise, set $h_2=(e_2\setminus\{f_{\mathcal{C}_3,e_2},l_{\mathcal{C}_3,e_2}\})\cup\{v_{k-1},z'\}$,
  where $z'\in \{f_{\mathcal{C}_2,f_{j-1}},l_{\mathcal{C}_2,f_{j}}\}\setminus \{z\}$. By Claim \ref{30}, $h_2$ is red.
 For $3\leq i\leq \frac{n}{2}-1$, set
\begin{eqnarray*}
h_i= \left\lbrace
\begin{array}{ll}
(e_i\setminus \{f_{\mathcal{C}_3,e_i},l_{\mathcal{C}_3,e_i}\})\cup\{v_{(i-1)(k-1)},z_{\frac{i-1}{2}}\}  &\mbox{if~}  i~\mbox{is~odd},\vspace{.5 cm}\\
(e_i\setminus\{f_{\mathcal{C}_3,e_i}\})\cup\{z_{\frac{i}{2}-1}\}
&\mbox{if~} i~\mbox{is~even}.
\end{array}
\right.\vspace{.2 cm}
\end{eqnarray*} \noindent  By Claim \ref{25}, $h_i$'s, $3\leq i\leq
\frac{n}{2}-1$, are red and clearly $\mathcal{C}_4=h_1h_2\ldots
h_{\frac{n}{2}-1}$ is a red copy of
$\mathcal{C}^k_{\frac{n}{2}-1}$.
 Clearly, in each the above cases, $\mathcal{C}_4=h_1h_2\ldots h_{\frac{n}{2}-1}$ is a red copy of $\mathcal{C}^k_{\frac{n}{2}-1}$
  disjoint form $\mathcal{C}_1$ with desired properties. $\square$\\

\begin{emp}\label{32} Let  $n=6$.
There is a red copy of $\mathcal{C}_{2}^k$, say
$\mathcal{C}_4=h_1h_2$, disjoint from $\mathcal{C}_1$ so that for $1\leq i \leq 2$, $ k-2\leq|h_i\cap e_i|\leq k-1$. Moreover,
 \begin{eqnarray*}h_1=(e_1\setminus\{v_1,v_k\})\cup\{v,w\}
 \end{eqnarray*} for some $v\in e_2\setminus\{f_{\mathcal{C}_3,e_2}\}$ and $w\in V(\mathcal{H})\setminus V(\mathcal{C}_1\cup \mathcal{C}_3)$.
\end{emp}

\noindent{\bf Proof of Claim \ref{32}.} Let
$W'=V(\mathcal{H})\setminus V(\mathcal{C}_1\cup \mathcal{C}_2\cup
\mathcal{C}_3)$. Since
\begin{eqnarray*}4(k-1)\leq |V(\mathcal{C}_1\cup \mathcal{C}_2)|\leq 4(k-1)+2,
 \end{eqnarray*} clearly $0 \leq|W'|\leq 2$. First let $|W'|=2$ and $W'=\{z_1,z_2\}$. If there is a vertex $z\in W'$, say $z_1$, so that the edge $(e_1\setminus\{v_k\})\cup\{z_1\}$ is blue, then using
 Claim \ref{24},  the edge
 \begin{eqnarray*}
h=(e_1\setminus\{v_1,v_k\})\cup\{v_{2k-2},z_2\}
\end{eqnarray*}
 is red. If the
edge $(e_2\setminus\{v_k\})\cup\{z_2\}$ is red, then set
$h_2=(e_2\setminus\{v_k\})\cup\{z_2\}$. Otherwise, set $h_2=(e_2\setminus\{v_1,v_k\})\cup\{v_2,z_1\}$.
Using Claim
\ref{24}, the edge $h_2$
is red. By choosing $h_1=h$,  $\mathcal{C}_4=h_1h_2$ is the desired
cycle. Therefore, we may assume that for each $z\in W'=\{z_1,z_2\}$, the edge $(e_1\setminus\{v_k\})\cup\{z\}$ is red.
Now, using Claim \ref{24}, there is a vertex $z'\in W'$ so that the edge $(e_2\setminus\{v_k\})\cup\{z'\}$ is red. With no loss of generality assume that $z'=z_1$. By choosing $h_1=(e_1\setminus\{v_k\})\cup\{z_1\}$ and $h_2=(e_2\setminus\{v_k\})\cup\{z_1\}$, $\mathcal{C}_4=h_1h_2$ is  the favorable $\mathcal{C}^k_2$. \\

Now, assume that  $|W'|=1$ and  $W'=\{z\}$. Clearly $|V(\mathcal{C}_2)\setminus V(\mathcal{C}_1)|=1$.  Let $x\in V(\mathcal{C}_2)\setminus V(\mathcal{C}_1)$.
  If   $x\in
f_j\setminus\{f_{\mathcal{C}_2,f_j},l_{\mathcal{C}_2,f_j}\}$, for some edge $f_j$ of $\mathcal{C}_2$, then  using Claim \ref{21}, the edges $h_1=(e_1\setminus\{v_k\})\cup\{x\}$ and $h_2=(e_2\setminus\{v_k\})\cup\{x\}$ are red and $\mathcal{C}_4=h_1h_2$ is
the desired cycle.
 Therefore, we may assume that $x$ is   a first vertex
  of $f_j$ for some edge  $f_j$ of $\mathcal{C}_2$. If there is a vertex $\overline{x}\in \{x,z\}$  so that the edge $(e_1\setminus\{v_k\})\cup\{\overline{x}\}$ is blue, then using
 Claim \ref{31},  the edge
 \begin{eqnarray*}
h'=(e_1\setminus\{v_1,v_k\})\cup\{v_{2k-2},\overline{y}\}
\end{eqnarray*} is red where   $\overline{y}\in\{x,z\}\setminus \{\overline{x}\}$. If the
edge $(e_2\setminus\{v_k\})\cup\{\overline{y}\}$ is red, then set
$h_2=(e_2\setminus\{v_k\})\cup\{\overline{y}\}$. Otherwise, set $h_2=(e_2\setminus\{v_1,v_k\})\cup\{v_2,\overline{x}\}$.
Using Claim
\ref{31}, the edge $h_2$
is red. By choosing $h_1=h'$,  $\mathcal{C}_4=h_1h_2$ is the desired
cycle. Therefore, we may assume that for each $\overline{x}\in \{x,z\}$, the edge $(e_1\setminus\{v_k\})\cup\{\overline{x}\}$ is red. Now, using Claim \ref{31}, there is a vertex $y\in \{x,z\}$ so that the edge $(e_2\setminus\{v_k\})\cup\{y\}$ is red. By choosing $h_1=(e_1\setminus\{v_k\})\cup\{y\}$ and $h_2=(e_2\setminus\{v_k\})\cup\{y\}$,  $\mathcal{C}_4=h_1h_2$  is the desired  $\mathcal{C}^k_2$.\\

 So we may assume that $|W'|=0$.
  This is the case only when
$|V(\mathcal{C}_1\cup \mathcal{C}_2)|=
(k-1)(\frac{n}{2}+1)+2$. Clearly, $|V(\mathcal{C}_2)\setminus
V(\mathcal{C}_1)|=2$. Similar to the
 discussion in the above paragraph,  we may assume that each vertex of $V(\mathcal{C}_2)\setminus V(\mathcal{C}_1)$ is a first vertex
  of $f_j$ for some edge  $f_j$ of $\mathcal{C}_2$. So we have two following cases.
\begin{itemize}
\item[(i)] There is an edge $f_j$ of $\mathcal{C}_2$ so that $V(\mathcal{C}_2)\setminus V(\mathcal{C}_1)=\{f_{\mathcal{C}_2,f_j},l_{\mathcal{C}_2,f_j}\}$,
 \item[(ii)] There is an edge $f_j$ so that
$V(\mathcal{C}_2)\setminus V(\mathcal{C}_1)=\{f_{\mathcal{C}_2,f_{j-1}},l_{\mathcal{C}_2,f_{j}}\}$.

 \end{itemize}
Note that, in  each  cases, we can find the favorable red $\mathcal{C}_4$, by an argument that used for case $|W'|=1$ (use Claim \ref{211} for case ${\rm (i)}$ and  Claim \ref{30} for case ${\rm (ii)}$). So we are done. $\square$\\\\

 The proof of the
following claim is similar to the proof of  Claim \ref{22}. So we
omit it here.


\begin{emp}\label{26}
Let $h_i$ and $d_j$ be two arbitrary edges of $\mathcal{C}_4$ and
$\mathcal{C}_1$, respectively and $y\in h_i$.
 For every vertices   $x\in d_j\setminus \{f_{\mathcal{C}_1,d_j}\}$ and $x'\in
d_{j+2}\setminus \{l_{\mathcal{C}_1,d_{j+2}}\}$, there are
distinct vertices
 $\overline{v},v'\in
h_i\setminus \{f_{\mathcal{C}_4,h_i},l_{\mathcal{C}_4,h_i},y\}$,
$u\in d_{j-1}\setminus \{f_{\mathcal{C}_1,d_{j-1}}\}$ and $u'\in
d_{j+3}\setminus \{l_{\mathcal{C}_1,d_{j+3}}\}$ so that  at least
one of
 the edges
$g=(V(d_j)\setminus \{f_{\mathcal{C}_1,d_j},x\})\cup\{\overline{v},u\}$ or
$g'=(V(d_{j+2})\setminus
\{x',l_{\mathcal{C}_1,d_{j+2}}\})\cup\{v',u'\}$ is red.\\

\end{emp}

In the rest of this section, we show that there is a red copy of $\mathcal{C}_n^k$.
First let $n=6$. By Claim \ref{32}, there is a red cycle $\mathcal{C}_4=h_1h_2$ disjoint from $\mathcal{C}_1$ where
 $h_1=(e_1\setminus\{v_1,v_k\})\cup\{v,w\}$
 for some $v\in e_2\setminus\{f_{\mathcal{C}_3,e_2}\}$ and $w\in V(\mathcal{H})\setminus (V(\mathcal{C}_1)\cup V(\mathcal{C}_3))$.
Since
  \begin{eqnarray*}
4(k-1)\leq |V(\mathcal{C}_1\cup \mathcal{C}_2)|\leq 4(k-1)+2,
 \end{eqnarray*} there is a vertex $z\in (d_l\setminus\{f_{\mathcal{C}_1,d_l},l_{\mathcal{C}_1,d_l}\})\cap(f_{l'}\setminus\{f_{\mathcal{C}_2,f_{l'}},l_{\mathcal{C}_2,f_{l'}}\})$ for some $1\leq l,l'\leq 4.$ Using Claim \ref{21} the edge
\begin{eqnarray*}
 g=(h_{1}\setminus\{v,w\})\cup\{z,v_1\}=(e_{1}\setminus\{v_k\})\cup\{z\} \end{eqnarray*} is red.
With no loss of generality assume that $d_l=d_1$.
Now, using Claim
 \ref{26} for $h_i=h_{2}$, $y=v_1$, $d_j=d_{4}$, $x=l_{\mathcal{C}_1,d_{4}}$ and $x'=f_{\mathcal{C}_1,d_{2}}$, there are distinct  vertices $\overline{v},v'\in h_{2}\setminus \{f_{\mathcal{C}_4,h_{2}},l_{\mathcal{C}_4,h_{2}},y\}$, $u\in d_{3}\setminus\{f_{\mathcal{C}_1,d_{3}}\}$ and $u'\in d_{3}\setminus\{l_{\mathcal{C}_1,d_{3}}\}$ so that  at least one of the edges $g'=(d_{4}\setminus\{x,f_{\mathcal{C}_1,d_{4}}\})\cup\{u,\overline{v}\}$
or $g''=(d_{2}\setminus\{x',l_{\mathcal{C}_1,d_{2}}\})\cup\{u',v'\}$  is red. If $g'$ is red, then $g'd_{3}d_{2}d_1 g h_2$ is a red copy of
$\mathcal{C}^k_6$,  a contradiction. If $g''$ is red, then $gd_1d_{4}d_{3}g''  h_2$ is a red copy of
$\mathcal{C}^k_n$,  a contradiction to our assumption. \\

Now, let $n> 6$. By Claim \ref{28}, there is a red cycle $\mathcal{C}_4=h_1h_2\ldots h_{\frac{n}{2}-1}$, disjoint from
$\mathcal{C}_1$ so that
\begin{eqnarray*}
 |h_{\frac{n}{2}-1}\cap e_{\frac{n}{2}-1}|=k-1,   h_{\frac{n}{2}-1}\setminus\{f_{\mathcal{C}_4,h_{\frac{n}{2}-1}}\}=e_{\frac{n}{2}-1}\setminus\{f_{\mathcal{C}_3,e_{\frac{n}{2}-1}}\}
\end{eqnarray*}
  and for each $1\leq i\leq \frac{n}{2}-1$,
$|h_i\cap e_{i}|\geq k-2$.  Since
\begin{eqnarray*}(k-1)(\frac{n}{2}+1)\leq|V(\mathcal{C}_1\cup \mathcal{C}_2)|\leq (k-1)(\frac{n}{2}+1)+\lfloor\frac{n}{4}\rfloor+1,
 \end{eqnarray*}
 there is a vertex $z\in (d_l\setminus\{f_{\mathcal{C}_1,d_l},l_{\mathcal{C}_1,d_l}\})\cap(f_{l'}\setminus\{f_{\mathcal{C}_2,f_{l'}},l_{\mathcal{C}_2,f_{l'}}\})$ for some $1\leq l,l'\leq \frac{n}{2}+1.$ Using Claim \ref{21} the edge \begin{eqnarray*} g=(h_{\frac{n}{2}-1}\setminus\{f_{\mathcal{C}_4,h_{\frac{n}{2}-1}}\})\cup\{z\}=
(e_{\frac{n}{2}-1}\setminus\{f_{\mathcal{C}_3,e_{\frac{n}{2}-1}}\})\cup\{z\}
\end{eqnarray*} is red. Now, using Claim
 \ref{26} for $h_i=h_{\frac{n}{2}-2}$, $d_j=d_{l-1}$, $x=l_{\mathcal{C}_1,d_{l-1}}$ and $x'=f_{\mathcal{C}_1,d_{l+1}}$, there are distinct  vertices $\overline{v},v'\in h_{\frac{n}{2}-2}\setminus \{f_{\mathcal{C}_4,h_{\frac{n}{2}-2}},l_{\mathcal{C}_4,h_{\frac{n}{2}-2}}\}$, $u\in d_{l-2}\setminus\{f_{\mathcal{C}_1,d_{l-2}}\}$ and $u'\in d_{l+2}\setminus\{l_{\mathcal{C}_1,d_{l+2}}\}$ so that  at least one of the edges $g'=(d_{l-1}\setminus\{x,f_{\mathcal{C}_1,d_{l-1}}\})\cup\{u,\overline{v}\}$
or $g''=(d_{l+1}\setminus\{x',l_{\mathcal{C}_1,d_{l+1}}\})\cup\{u',v'\}$  is red. If $g'$ is red, then
$g'd_{l-2}d_{l-3}\ldots d_1d_{\frac{n}{2}+1}d_{\frac{n}{2}}\ldots d_l g h_1h_2\ldots h_{\frac{n}{2}-2}$ is a red copy of
$\mathcal{C}^k_n$,  a contradiction. If $g''$ is red, then
 \begin{eqnarray*}
gd_ld_{l-1}d_{l-2}\ldots d_1d_{\frac{n}{2}+1}d_{\frac{n}{2}}\ldots d_{l+2} g'' h_{\frac{n}{2}-2} h_{\frac{n}{2}-3}\ldots h_2h_1
\end{eqnarray*} is a red copy of
$\mathcal{C}^k_n$,  a contradiction.

 $\hfill\blacksquare$\\\\



\section{Concluding remarks and open problems}

Throughout this paper, we consider $k\geq 8$ for simplicity. We believe that our approach can be used to prove Conjecture
  \ref{our conjecture2} for $n=m$ and  $k=7$, however much more details are required. Therefore,
it would be interesting to investigate Conjecture \ref{our conjecture2} for $n=m$ and $3< k< 7$.  \\

It is known   that  Conjecture \ref{our conjecture2} is true for a fixed $m\geq 3$ if it holds for every $m\leq
n\leq 2m$ (\cite{Ramsey numbers of loose cycles in uniform hypergraphs}).
So it would  be interesting to deduce whether  Conjecture \ref{our conjecture2} holds for every $m\leq n\leq 2m$ and $k\geq 4$.
 It seems that our method   can be used to prove Conjecture \ref{our conjecture2} for $m\leq n \leq 2m$ and sufficiently large $k$,  but  too much efforts and  details are needed.\\

Another interesting  question in this direction is to prove Conjecture \ref{our conjecture2}  for small values of $k$.
 As we noted in the introduction, the case $k=3$  is proved in \cite{Ramsey numbers of loose cycles in uniform hypergraphs} and \cite{The Ramsey number of loose paths  and loose cycles in 3-uniform hypergraphs}.


\section{Appendix A}

 \noindent\textbf{Proof the part $(i)$ of Theorem \ref{R(Pk3,Pk3)}.}
  Let $\mathcal{H}=\mathcal{K}^k_{2(k-1)}$ is two edge
 colored red and blue. If $\mathcal{H}$ is monochromatic, then clearly we are done. So we may assume $\mathcal{H}_{\rm red}$ and $\mathcal{H}_{\rm blue}$ are both non empty.
 It is easy to see that if a $k$-uniform complete hypergraph is $2$-colored
 and both colors are used at least once, then there are two edges
 of distinct colors intersecting in $k-1$ vertices (see Remark $3$ of \cite{subm}).
 Thereby,
 we can select $e=\{v_1,v_2,\ldots,v_k\}\in \mathcal{H}_{\rm
 red}$ and $f=\{v_2,v_3,\ldots,v_{k+1}\}\in \mathcal{H}_{\rm
 blue}.$ Let $W=V(\mathcal{H})\setminus
 \{v_1,v_2,\ldots,v_{k+1}\}.$ Clearly $|W|=k-3$.  Consider the edge $g=\{v_1,v_2,v_{k+1}\}\cup W$. If $g$ is red, then $eg$ is a red copy of $\mathcal{C}^k_2$. Otherwise, $fg$ is a blue $\mathcal{C}^k_2$. $\hfill\blacksquare$\\

 \noindent\textbf{Proof of Lemma \ref{No  red Cn-1 implies blue Cn}. }Let $\mathcal{C}=e_1e_2\ldots e_n$ be a
 copy of $\mathcal{C}_{n}^k$ in  $\mathcal{H}_{\rm red}$ with edges
 $$e_i=\{v_1,v_2,\ldots,v_k\}+(k-1)(i-1) \hspace{0.5 cm} ({\rm mod}\ \  (k-1)n), \hspace{0.5 cm}i=1,\ldots, n.$$ First assume that $n$ is odd.
 Let
 \begin{eqnarray*}
 f_i= \left\lbrace
 \begin{array}{ll}
 (e_{2i-1}\setminus\{l_{\mathcal{C},e_{2i-1}}\})\cup\{l_{\mathcal{C},e_{2i}}\} &  1\leq i \leq \frac{n+1}{2},\vspace{.5 cm}\\
 (e_{2i-n-1}\setminus\{l_{\mathcal{C},e_{2i-n-1}}\})\cup\{l_{\mathcal{C},e_{2i-n}}\}&
 \frac{n+3}{2}\leq i\leq n.
 \end{array}
 \right.\vspace{.2 cm}
 \end{eqnarray*}
  Since there is no red copy of
 $\mathcal{C}_{n-1}^k$, all $f_i$'s ,$1\leq i\leq  n$ are blue
 (otherwise, for some $i$, the edges  $f_ie_{2i+1}e_{2i+2}\ldots
 e_n e_1\ldots e_{2i-2}$ form a red $\mathcal{C}_{n-1}^k$). Clearly
 $f_1f_2\ldots f_n$ is
 a blue copy of  $\mathcal{C}_{n}^k$.\\

 Now assume that $n$ is even. Let


 \begin{eqnarray*} \hspace{-1 cm}f_i= \left\lbrace
 \begin{array}{ll}
 (e_{2i-1}\setminus \{v_{(2i-1)(k-1)},v_{(2i-1)(k-1)+1}\})\cup
 \{v_{2i(k-1)},v_{2i(k-1)+1}\} &  1\leq i \leq \frac{n}{2}-1,\vspace{.5 cm}\\
 (e_{n-1}\setminus \{v_{(n-1)(k-1)},v_{(n-1)(k-1)+1}\})\cup
 \{v_{n(k-1)},v_{k-1}\}&  i=\frac{n}{2},\vspace{.5 cm}\\
 \{v_{a(k-1),v_{a(k-1)+1}}\}\cup(e_{a-1}\setminus
 \{v_{(a-1)(k-1)},v_{(a-1)(k-1)+1}\})&  \frac{n}{2}+1\leq i\leq n-1,\vspace{.5 cm}\\
 \{v_{k-2},v_{3(k-1)},v_{3(k-1)+1}\}\cup (e_2\setminus
 \{v_k,v_{2(k-1)},v_{2(k-1)+1}\})&  i=n,

 \end{array}
 \right.\vspace{.3 cm}
 \end{eqnarray*}
 where $1\leq a\leq n$ and  $a=3-2i({\rm mod}\  n)$.\\
 It is obvious that all $f_i$'s are blue. So $f_1\ldots f_n$ is a
 copy of $\mathcal{C}_{n}^k$ in $\mathcal{H}_{\rm blue}.$
 $\hfill\blacksquare$\\

 \noindent\textbf{Proof of Lemma \ref{No  red Cn-2 implies blue Cn}. }Let $\mathcal{C}=e_1e_2\ldots e_n$ be a
 copy of $\mathcal{C}_{n}^k$ in  $\mathcal{H}_{\rm red}$ with edges
 $$e_i=\{v_1,v_2,\ldots,v_k\}+(k-1)(i-1)\hspace{0.5 cm} ({\rm mod}\ \
 (k-1)n),\hspace{0.5 cm} i=1,\ldots, n.$$  First assume that
 $n\equiv 1,2$ (mod 3). Let $f_i=e_{3i-2}\setminus
 \{l_{\mathcal{C},e_{3i-2}}\}\cup\{l_{\mathcal{C},e_{3i}}\}$,
 $1\leq i \leq n$ where indices the $e_i$'s are mod $n$. Since
 there is no red copy of $\mathcal{C}_{n-2}^k$, all $f_i$'s ,$1\leq
 i\leq n$ are blue (otherwise, for some $i$, the edges
 $f_ie_{3i+1}\ldots e_{3i-3}$ form a red $\mathcal{C}_{n-2}^k$).
 Clearly $f_1f_2\ldots f_n$ is
 a blue copy of $\mathcal{C}_{n}^k$.\\

 Now assume that $n\equiv 0$ (mod 3). Partition the vertices of
 $e_i$ into three  parts $A_i,B_i$ and $C_i$ with
 $|A_i|\geq|B_i|\geq|C_i|\geq|A_i|-1$ so that
 $f_{\mathcal{C},e_i}\in A_i$ and $l_{\mathcal{C},e_i}\in
 C_i$. Clearly $|C_i|\geq 2$. Let $v\in C_1\setminus \{v_k\},$
 $v'\in B_1$ and $v''\in B_2.$ Set

 \begin{eqnarray*}
 f_i= \left\lbrace
 \begin{array}{ll}
 A_{3(i-1)+1}\cup B_{3(i-1)+2}\cup C_{3(i-1)+3}  &  1\leq i\leq \frac{n}{3}-1,\vspace{.5 cm}\\
 A_{n-2}\cup B_{n-1}\cup (C_{n}\setminus\{v_1\})\cup \{v\}&  i=\frac{n}{3},\vspace{.5 cm}\\
 C_{2n+4-3i}\cup B_{2n+3-3i}\cup A_{2n+2-3i}&  \frac{n}{3}+1\leq i\leq \frac{2n}{3}-1,\vspace{.5 cm}\\
 C_{4}\cup B_{3}\cup (A_{2}\setminus\{v_k\})\cup \{v'\}&  i=\frac{2 n}{3},\vspace{.5 cm}\\
 C_{3n+5-3i}\cup B_{3n+4-3i}\cup A_{3n+3-3i}&  \frac{2n}{3}+1\leq i\leq n-1,\vspace{.5 cm}\\
 C_{5}\cup B_{4}\cup (A_{3}\setminus\{v_{2(k-1)+1}\})\cup \{v''\}&
 i=n,\vspace{.5 cm}
 \end{array}
 \right.\vspace{.2 cm}
 \end{eqnarray*}
 where the indices are mod $n$.\\
 It is obvious that all $f_i$'s are blue. So $f_1\ldots f_n$ is a
 copy of $\mathcal{C}_{n}^k$ in $\mathcal{H}_{\rm blue}.$
  $\hfill\blacksquare$\\

 \noindent{\bf Proof of Lemma \ref{there is a P2}.} Suppose for a contradiction that there is no red path
 $\mathcal{P}$ with the above  conditions.
  By symmetry we may assume that $i=j=1$,
 $e_1=\{v_1,v_2,\ldots,v_k\}$, $f_1=\{u_1,u_2,\ldots,u_k\}$ and
 $C\subseteq\{v_{k-1}\}$. Let  $|C|=l\in \{0,1\}$. Consider an edge
 $h=E_1\dot{\cup}W_1\dot{\cup}F_1$  in
 $\mathcal{A}_{11}\cup\mathcal{B}_{11}$ so that
 \begin{itemize}
 \item $W_1\subseteq \{w_1,w_2\}$, $|W_1|=l$. \item $E_1\subseteq
 \Big(V(e_1)\setminus (C\cup\{v_1,v_k\})\Big)\cup\{x\}$, $x\in
 \{v',v''\}$.
  \item $F_1\subseteq
 (V(f_1)\setminus \{u_1,u_k\})\cup\{y\}$, $y\in \{u',u''\}$. \item
 $\Big| |E_1|-|F_1|\Big|\leq 1$.
 \end{itemize}
 \begin{emp}\label{g is red}
 The  edge  $h$ is red.
 \end{emp}
 {\bf Proof of Claim \ref{g is red}.} By symmetry we may assume
 that $h\in \mathcal{A}_{11}$, $|F_1|\geq |E_1|$ and
 $W_1\subseteq\{w_1\}$. It is easy to see that
 $|E_1|=\lfloor\frac{k-l}{2}\rfloor$ and
 $|F_1|=\lceil\frac{k-l}{2}\rceil$. W.l.g. assume that
 \begin{eqnarray*}
 &&E_1=\{v',v_2,\ldots,v_{\lfloor\frac{k-l}{2}\rfloor}\},\\
 &&F_1=\{u_{k-\lceil\frac{k-l}{2}\rceil+1},\ldots,u_{k-1},u''\}.
 \end{eqnarray*}
 Suppose indirectly that the edge $h_1=h$ is blue. Since there is
 no blue copy of $\mathcal{C}_{l_1+l_2}$, using Remark
 \ref{complementary edges are not blue}, every edge in
 $\mathcal{B}_{11}$ that is disjoint from $h_1$ is red.
  Now let
 $h_2=(h_1\setminus\{v_2\})\cup \{v_{k-1-l}\}$.  Clearly
 $h_2\in\mathcal{A}_{11}$. If $h_2$ is red, then set
 $$h'_2=\Big((e_1\cup f_1\cup \{w_2\})\setminus(h_1\cup C\cup \{v_1,v_k,u_1,u_k,u\})\Big)\cup\{u',v''\}$$
  where $u\in
  f_1\setminus(h_2\cup\{u_1,u_k\})$. Since there is no blue copy of
  $\mathcal{C}_{l_1+l_2}$, $\mathcal{P}=g_1g_2$ is the desired path
  where $g_1=h_2$ and $g_2=h'_2$ (clearly $w\in e\setminus \mathcal{P}$ where $w=v_2$ for $e=e_1$ and $w=u$ for $e=f_1$).
  Therefore, we may assume that the edge $h_2$ is
  blue.
  For $2\leq i\leq \lfloor\frac{k-l}{2}\rfloor$,  let
  $h_i=(h_{i-1}\setminus\{v_i\})\cup\{v_{k-(i+l)+1}\}$. Assume that
  $j$ is  the maximum $i\in[1,\lfloor\frac{k-l}{2}\rfloor]$ for
  which $h_i$ is blue. If $j<\lfloor\frac{k-l}{2}\rfloor$, then $h_{j+1}$ is red. Set
  $$h'_{j+1}=\Big((e_1\cup f_1\cup \{w_2\})\setminus(h_{j}\cup C\cup \{v_1,v_k,u_1,u_k,u\})\Big)\cup\{u',v''\},$$
  where  $u\in
  f_1\setminus(h_{j+1}\cup\{u_1,u_k\})$. Clearly, $h'_{j+1}$ is
  red. Therefore, $\mathcal{P}=g_1g_2$ is the desired path where
  $g_1=h_{j+1}$and $g_2=h'_{j+1}$ (clearly $w\in e\setminus \mathcal{P}$ where $w=v_{j+1}$ for $e=e_1$ and $w=u$ for
  $e=f_1$). So we may assume that $j=\lfloor\frac{k-l}{2}\rfloor$ and hence
 the  edge
 $h_{\lfloor\frac{k-l}{2}\rfloor}=E'\dot{\cup}W_1\dot{\cup}F_1$ is
  blue, where
  $$E'=\{v',v_{k-l-(\lfloor\frac{k-l}{2}\rfloor-1)}, \ldots, v_{k-l-2},v_{k-l-1}\}.$$
 Now, set
 \begin{eqnarray*}
 m=\left\lbrace
 \begin{array}{ll}
 |F_1|-2 &\mbox{if~} \ l=0\ \mbox{and $k$ is odd},\vspace{.5
 cm}\\
 |F_1|-1  &\mbox{otherwise}.
 \end{array}
 \right.\vspace{.2 cm}
 \end{eqnarray*}
 For $1\leq i\leq m$, let
  $h_{\lfloor\frac{k-l}{2}\rfloor+i}=(h_{\lfloor\frac{k-l}{2}\rfloor+i-1}\setminus\{u_{k-i}\})\cup\{u_{i+1}\}$. Now, let
  $j$ be the maximum $i\in[0,m]$ for
  which $h_{\lfloor\frac{k-l}{2}\rfloor+i}$ is blue. If $j<m$, then $h_{\lfloor\frac{k-l}{2}\rfloor+j+1}$ is red. Set
  $$h'_{\lfloor\frac{k-l}{2}\rfloor+j+1}=\Big((e_1\cup f_1\cup \{w_2\})\setminus(h_{\lfloor\frac{k-l}{2}\rfloor+j}\cup C\cup
   \{v_1,v_k,u_1,u_k,v\})\Big)\cup\{u',v''\},$$
  where  $v\in
  e_1\setminus(h_{\lfloor\frac{k-l}{2}\rfloor+j+1}\cup\{v_1,v_k\}\cup C)$.
  Since there is no blue copy of $\mathcal{C}_{l_1+l_2}$, the edge
  $h'_{\lfloor\frac{k-l}{2}\rfloor+j+1}$ is red.
   Therefore, $\mathcal{P}=g_1g_2$ is the desired path where
  $g_1=h_{\lfloor\frac{k-l}{2}\rfloor+j+1}$and $g_2=h'_{\lfloor\frac{k-l}{2}\rfloor+j+1}$
   (clearly $w\in e\setminus \mathcal{P}$ where $w=v$ for $e=e_1$ and $w=u_{k-j-1}$ for
  $e=f_1$). So we may assume that $j=m$ and
  hence the
  edge $h_{\lfloor\frac{k-l}{2}\rfloor+m}=E'\dot{\cup}W_1\dot{\cup}F'$ is
  blue, where
 \begin{eqnarray*}
 F'=\left\lbrace\begin{array}{ll}
 \{u_2,\ldots,u_{m+1},u_{\frac{k+1}{2}},u''\} &\mbox{if~} \ l=0\
 \mbox{and $k$ is odd},\vspace{.5
 cm}\\
 \{u_2,\ldots,u_{m+1},u''\} &\mbox{otherwise}. \end{array}
 \right.\vspace{.2 cm}\end{eqnarray*}\\

 Set
 $h_{\lfloor\frac{k-l}{2}\rfloor+m+1}=(h_{\lfloor\frac{k-l}{2}\rfloor+m}\setminus\{v'\})\cup\{v''\}$.
  If $h_{\lfloor\frac{k-l}{2}\rfloor+m+1}$ is red, then set
 \begin{eqnarray*}
 &&\hspace{-1.5
 cm}h'_{\lfloor\frac{k-l}{2}\rfloor+m+1}=\Big((e_1\cup f_1\cup
 \{w_2\})\setminus(h_{\lfloor\frac{k-l}{2}\rfloor+m}\cup C\cup
 \{v_1,v_k,u_1,u_k,v\})\Big)\cup\{v'',u'\},\\
 \vspace*{0.4 cm} &&\hspace{-1.5
 cm}h''_{\lfloor\frac{k-l}{2}\rfloor+m+1}=\Big((e_1\cup f_1\cup
 \{w_2\})\setminus(h_{\lfloor\frac{k-l}{2}\rfloor+m}\cup C\cup
 \{v_1,v_k,u_1,u_k,u\})\Big)\cup\{v'',u'\}.
 \end{eqnarray*}
  where  $v\in
  e_1\setminus(h_{\lfloor\frac{k-l}{2}\rfloor+m+1}\cup\{v_1,v_k\}\cup C)$ and $u\in  f_1\setminus(h_{\lfloor\frac{k-l}{2}\rfloor+m+1}\cup\{u_1,u_k\})$.
 For $e=e_1$, set
 $\mathcal{P}=h_{\lfloor\frac{k-l}{2}\rfloor+m+1}h'_{\lfloor\frac{k-l}{2}\rfloor+m+1}$
 and for $e=f_1$ set
 $\mathcal{P}=h_{\lfloor\frac{k-l}{2}\rfloor+m+1}h''_{\lfloor\frac{k-l}{2}\rfloor+m+1}$.
 It is easy to see that $\mathcal{P}$ is the desired path. Hence,
 we may assume that the edge $h_{\lfloor\frac{k-l}{2}\rfloor+m+1}$
 is blue.\\

 Similarly, we may assume that the edge
 $h_{\lfloor\frac{k-l}{2}\rfloor+m+2}=(h_{\lfloor\frac{k-l}{2}\rfloor+m+1}\setminus\{u''\})\cup\{u'\}$
 is blue.  If $l=0$ and $k$ is even, then clearly
 $h_{\lfloor\frac{k-l}{2}\rfloor+m+2}$ is an edge in
 $\mathcal{B}_{11}$ disjoint from  $h_1$. This is impossible, by
 Remark \ref{complementary edges are not blue}.
  Now we
 have one of the following
 cases.\\\\
 {\bf Case 1: $l=0$ and $k$ is odd.}

  One can easily check that
 $u_{\frac{k+1}{2}}\in h_1\cap h_{\lfloor\frac{k-l}{2}\rfloor+m+2}$
 and $v_{\frac{k+1}{2}}\notin h_1\cup
 h_{\lfloor\frac{k-l}{2}\rfloor+m+2} $. Let
 $h_{\lfloor\frac{k-l}{2}\rfloor+m+3}=(h_{\lfloor\frac{k-l}{2}\rfloor+m+2}\setminus\{u_{\frac{k+1}{2}}\})\cup\{v_{\frac{k+1}{2}}\}$.
 If the edge $h_{\lfloor\frac{k-l}{2}\rfloor+m+3}$ is red, then set
 $$h'_{\lfloor\frac{k-l}{2}\rfloor+m+3}=\Big((e_1\cup f_1\cup
 \{w_2\})\setminus(h_{\lfloor\frac{k-l}{2}\rfloor+m+2}\cup
 \{v_1,v_k,u_1,u_k,v\})\Big)\cup\{v',u''\},$$  where  $v\in
  e_1\setminus(h_{\lfloor\frac{k-l}{2}\rfloor+m+3}\cup\{v_1,v_k\})$.
  It is easy to see that $\mathcal{P}=g_1g_2$ is the desired path
  where $g_1=h_{\lfloor\frac{k-l}{2}\rfloor+m+3}$ and
  $g_2=h'_{\lfloor\frac{k-l}{2}\rfloor+m+3}$. Therefore, we may
  assume that the edge $h_{\lfloor\frac{k-l}{2}\rfloor+m+3}$ is
  blue.
  That is a contradiction to Remark \ref{complementary
 edges are not blue}, since $h_{\lfloor\frac{k-l}{2}\rfloor+m+3}$
 is an edge in $\mathcal{B}_{11}$
 disjoint from $h_1$.\\\\
 {\bf Case 2: $l=1$.}

  First let $k$ is even. Clearly,
 $v_{\frac{k}{2}}\notin h_1\cup
 h_{\lfloor\frac{k-l}{2}\rfloor+m+2}$ and $W_1=\{w_1\}$. Let
 $h_{\lfloor\frac{k-l}{2}\rfloor+m+3}=(h_{\lfloor\frac{k-l}{2}\rfloor+m+2}\setminus\{w_1\})\cup\{v_{\frac{k}{2}}\}$. If $h_{\lfloor\frac{k-l}{2}\rfloor+m+3}$ is
 red, then set
 \begin{eqnarray*}
 &&\hspace{-1.5 cm}
 h'_{\lfloor\frac{k-l}{2}\rfloor+m+3}=\Big((e_1\cup f_1\cup
 \{w_2\})\setminus(h_{\lfloor\frac{k-l}{2}\rfloor+m+2}\cup C\cup
 \{v_1,v_k,u_1,u_k,v\})\Big)\cup\{v',u''\},\\
 \vspace*{0.4 cm} &&\hspace{-1.5
 cm}h''_{\lfloor\frac{k-l}{2}\rfloor+m+3}=\Big((e_1\cup f_1\cup
 \{w_2\})\setminus(h_{\lfloor\frac{k-l}{2}\rfloor+m+2}\cup C\cup
 \{v_1,v_k,u_1,u_k,u\})\Big)\cup\{v',u''\}.
 \end{eqnarray*}
  where  $v\in
  e_1\setminus(h_{\lfloor\frac{k-l}{2}\rfloor+m+3}\cup\{v_1,v_k,v_{k-1}\})$ and $u\in  f_1\setminus(h_{\lfloor\frac{k-l}{2}\rfloor+m+3}\cup\{u_1,u_k\})$.
  Let $g_2=h_{\lfloor\frac{k-l}{2}\rfloor+m+3}$. It is easy to
  check that $\mathcal{P}=g_1g_2$ is the desired path where
 $g_1=h'_{\lfloor\frac{k-l}{2}\rfloor+m+3}$ for $e=e_1$ and
  $g_1=h''_{\lfloor\frac{k-l}{2}\rfloor+m+3}$ for $e=f_1$.
  So, we may assume that
  the edge $h_{\lfloor\frac{k-l}{2}\rfloor+m+3}$ is blue, a
  contradiction to Remark \ref{complementary edges are not blue}.

  Now,
 let $k$ be odd. One can easily see that $u_{\frac{k+1}{2}}\notin
 h_1\cup h_{\lfloor\frac{k-l}{2}\rfloor+m+2}$. Similarly,
  we may assume
 that the edge
 $h_{\lfloor\frac{k-l}{2}\rfloor+m+3}=h_{\lfloor\frac{k-l}{2}\rfloor+m+2}\setminus
 \{w_1\})\cup\{u_{\frac{k+1}{2}}\}$  is blue. That is a
 contradiction to Remark \ref{complementary edges are not blue}.
 This contradiction
 completes the proof of Claim \ref{g is red}. $\square$\\\\

 \noindent By Claim \ref{g is red}
 we can find the favorable red $\mathcal{P}$, as follow.\\

 \noindent First let $l=1$. Set
 $$g_1=\{v',v_2,\ldots,v_{\lfloor\frac{k-1}{2}\rfloor}\}\cup\{w_1\}\cup\{u_{\lfloor\frac{k-1}{2}\rfloor+2},\ldots,u_{k-1},u''\}$$
 and
 $$g_2=\{v_{\lfloor\frac{k-1}{2}\rfloor+2},\ldots,v_{k-2},v''\}\cup\{w_2\}\cup\{u',u_2,u_3,\ldots,u_{\lfloor\frac{k-1}{2}\rfloor},u_{\lfloor\frac{k-1}{2}\rfloor+2}\}.$$
 Clearly, by Claim \ref{g is red}, $\mathcal{P}=g_1g_2$ is the
 desired path. Now let $l=0$. Set
 $$g_1=\{v',v_2,\ldots,v_{\lfloor\frac{k}{2}\rfloor}\}\cup\{u_{\lfloor\frac{k}{2}\rfloor+1},\ldots,u_{k-1},u''\}.$$
 For $e=e_1$ let
 $$g_2=\{v_{\lfloor\frac{k}{2}\rfloor},v_{\lfloor\frac{k}{2}\rfloor+2},\ldots,v_{k-1},v''\}\cup\{u',u_2,\ldots,u_{\lfloor\frac{k}{2}\rfloor}\}$$
 and for $e=f_1$ let
 $$g_2=\{v_{\lfloor\frac{k}{2}\rfloor+1},\ldots,v_{k-1},v''\}\cup\{u',u_2,\ldots,u_{\lfloor\frac{k}{2}\rfloor-1},u_{\lfloor\frac{k}{2}\rfloor+1}\}.$$
 It is easy to see that $\mathcal{P}=g_1g_2$ is the desired path.
 So we are done.
  $\hfill\blacksquare$\\

 \noindent{\bf Proof of Lemma \ref{there is a P2:3}.} Suppose not.
  By symmetry we may assume that $i=j=1$,
 $e_1=\{v_1,v_2,\ldots,v_k\}$ and $f_1=\{u_1,u_2,\ldots,u_k\}$.  Consider an edge
 $h=E\dot{\cup}W'\dot{\cup}F$  in
 $\mathcal{A}_{11}\cup\mathcal{B}_{11}$ so that
 \begin{itemize}
 \item $W'\subseteq B$, $|W'|=1$. \item $E\subseteq
 V(e_1)$,  $F\subseteq V(f_1)$ and $\Big| |E|-|F|\Big|\leq 1$.
 \end{itemize}
 \begin{emp}\label{g is red:jadid}
 The  edge  $h$ is red.
 \end{emp}
 {\bf Proof of Claim \ref{g is red:jadid}.} By symmetry we may assume
 that $h\in \mathcal{A}_{11}$, $|F|\geq |E|$ and
 $W'=\{w_1\}$. It is easy to see that
 $|E|=\lfloor\frac{k-1}{2}\rfloor$ and
 $|F|=\lceil\frac{k-1}{2}\rceil$. W.l.g. assume that
 \begin{eqnarray*}
 &&E=\{v_1,v_2,\ldots,v_{\lfloor\frac{k-1}{2}\rfloor}\},\\
 &&F=\{u_{k-\lceil\frac{k-1}{2}\rceil+1},\ldots,u_{k-1},u_k\}.
 \end{eqnarray*}
 Suppose indirectly that the edge $h_1=h$ is blue. Since there is
 no blue copy of $\mathcal{C}_{l_1+l_2}$, using Remark
 \ref{complementary edges are not blue}, every edge in
 $\mathcal{B}_{11}$ that is disjoint from $h_1$ is red.
 For $2\leq i\leq \lfloor\frac{k-1}{2}\rfloor$,  let
  $h_i=(h_{i-1}\setminus\{v_i\})\cup\{v_{k-i+1}\}$. Now, let
  $j$ be the maximum $i\in[1,\lfloor\frac{k-1}{2}\rfloor]$ for
  which $h_i$ is blue. If $j<\lfloor\frac{k-1}{2}\rfloor$, then $h_{j+1}$ is red. Set
  \begin{eqnarray*}
  &&h'_{j+1}=\Big((e_1\cup f_1\cup \{w_2\})\setminus(h_{j}\cup \{u,v\})\Big),\\
  &&h''_{j+1}=\Big((e_1\cup f_1\cup \{w_2\})\setminus(h_{j}\cup \{u,u'\})\Big),
  \end{eqnarray*}
  where $u,u',v$ are distinct vertices so that $u,u'\in
  f_1\setminus(h_{j+1}\cup\{u_1,u_k\})$ and $v\in e_1\setminus(h_{j+1}\cup\{v_1,v_k,v_{j+1}\})$ . Clearly, $h'_{j+1}$ and $h''_{j+1}$ are
  red. Set $g_1=h_{j+1}$, $g'_1=h'_{j+1}$ and $\overline{g'}_1=h''_{j+1}$. Then $E_1=g_1g'_1$ and $F_1=g_1\overline{g'}_1$ are  desired paths where $\overline{v}=v_{j+1}$ and  $\overline{u}=u$.
   So we may assume that $j=\lfloor\frac{k-1}{2}\rfloor$ and hence
 the  edge
 $h_{\lfloor\frac{k-1}{2}\rfloor}=E'\dot{\cup}W'\dot{\cup}F$ is
  blue, where
  $$E'=\{v_1,v_{k-\lfloor\frac{k-1}{2}\rfloor+1}, \ldots, v_{k-2},v_{k-1}\}.$$
  Now,  for $1\leq i \leq \lceil\frac{k-1}{2}\rceil-1$ let $h_{\lfloor\frac{k-1}{2}\rfloor+i}=(h_{\lfloor\frac{k-1}{2}\rfloor+i-1}\setminus\{u_{k-i}\})\cup\{u_{i+1}\}$.
 In a similar way we can show that the
  edge $h_{k-2}=E'\dot{\cup}W'\dot{\cup}F'$ is
  blue, where

 $$F'=\{u_2,u_3,\ldots,u_{\lceil\frac{k-1}{2}\rceil},u_k\}.$$

 Now, let $h_{k-1}=(h_{k-2}\setminus\{v_1,w_1\})\cup\{v_k,w_3\}$.
  If $h_{k-1}$ is red, then set
 \begin{eqnarray*}
 &&h'_{k-1}=\Big((e_1\cup f_1\cup
 \{w_2\})\setminus(h_{k-2}\cup
 \{u,v\})\Big),\\
 \vspace*{0.4 cm} &&h''_{k-1}=\Big((e_1\cup f_1\cup
 \{w_2\})\setminus(h_{k-2}\cup
 \{u,u'\})\Big).
 \end{eqnarray*}
 where $u,u',v$ are distinct vertices so that $u,u'\in
  f_1\setminus(h_{k-1}\cup\{u_1,u_k\})$ and $v\in e_1\setminus(h_{k-1}\cup\{v_1,v_k\})$. Clearly, $h'_{k-1}$ and $h''_{k-1}$ are red. Set $g_1=h_{k-1}$, $g'_1=h'_{k-1}$ and $\overline{g'}_1=h''_{k-1}$. Then $E_1=g_1g'_1$ and $F_1=g_1\overline{g'}_1$ are  desired
   paths where $\overline{v}=v_1$ and  $\overline{u}=u$.
  Hence,
 we may assume that the edge $h_{k-1}$
 is blue. Similarly, we may assume that the edge $h_k=(h_{k-1}\setminus\{u_k\})\cup\{u_1\}$ is blue. This is a contradiction to
  Remark \ref{complementary edges are not blue}. This contradiction completes the proof of Claim \ref{g is red:jadid}. $\square$\\

 Now, we can find favorable paths as follows:\\
 Let
 \begin{eqnarray*}
 &&g_1=\{v_1,v_2,\ldots,v_{\lfloor\frac{k-1}{2}\rfloor}\}\cup\{w_1\}\cup\{u_{k-\lceil\frac{k-1}{2}\rceil+1},\ldots,u_{k-1},u_k\},\\
 &&g'_1=\{v_{\lfloor\frac{k-1}{2}\rfloor},v_{\lfloor\frac{k-1}{2}\rfloor+3},\ldots,v_{k}\}\cup\{w_2\}\cup\{u_{1},\ldots,u_{\lfloor\frac{k-1}{2}\rfloor}\},\\
 &&\overline{g'}_1=\{v_{\lfloor\frac{k-1}{2}\rfloor+2},\ldots,v_{k}\}\cup\{w_2\}
 \cup\{u_{1},\ldots,u_{\lfloor\frac{k-1}{2}\rfloor-1},u_{k-\lceil\frac{k-1}{2}\rceil+1}\}.
 \end{eqnarray*}
 Using Claim \ref{g is red:jadid}, the edges $g_1$, $g'_1$ and $\overline{g'}_1$ are red and so  $E_1=g_1g'_1$ and $F_1=g_1\overline{g'}_1$ are desired paths  where $\overline{v}=v_{\lfloor\frac{k-1}{2}\rfloor+1}$,  $\overline{u}=u_{\lfloor\frac{k-1}{2}\rfloor+1}$ and $B'=\{w_1,w_2\}$.$\hfill\blacksquare$\\

 \medskip
 \noindent{\bf Proof of Lemma \ref{there is a P:3}.} By symmetry we may assume that $i=j=2$. Suppose for a
 contradiction that there is no red paths $E_2$ and $F_2$  with
 desired properties. Let
  $x\in (e_{1}\setminus\{v_1\})\cap({g'}_{1}\setminus
  {g}_{1})$. Assume that
 $h=E\cup F$ is an edge in $\mathcal{A}_{22}$ so that
 \begin{eqnarray*}
 &&E=\{x,v_{k+1},\ldots,v_{(k-1)+\lfloor\frac{k}{2}\rfloor}\},\\
 &&F=\{u_{2k-\lceil\frac{k}{2}\rceil},\ldots,u_{2k-2},u_{2k-1}\}.
 \end{eqnarray*}

 \begin{emp}\label{g is red:4}
 The edge $h$ is red.
 \end{emp}
 {\bf Proof of Claim \ref{g is red:4}.} Suppose indirectly that the
 edge $h_1=h$ is blue. Since there is no blue copy of
 $\mathcal{C}_{l_1+l_2}$, using Remark \ref{complementary edges are
 not blue}, every  edge in $\mathcal{B}_{22}$ that is disjoint from
 $h_1$ is red.
  Now let
 $h_2=(h_1\setminus\{x\})\cup\{v_{2k-1}\}$. If $h_2$ is red, then
 set \begin{eqnarray*}
 &&h'_2=\Big((e_2\cup f_2\cup
 \{w\})\setminus(h_1\cup
 \{f_{\mathcal{C}_1,e_{2}},f_{\mathcal{C}_2,f_{2}},\tilde{v}\})\Big)\cup\{\overline{u}\},\\
 &&h''_2=\Big((e_2\cup f_2\cup \{w\})\setminus(h_1\cup
 \{f_{\mathcal{C}_1,e_{2}},f_{\mathcal{C}_2,f_{2}},\tilde{u}\})\Big)\cup\{\overline{u}\},
 \end{eqnarray*}
 where $\tilde{v}\in
 e_2\setminus(h_2\cup\{f_{\mathcal{C}_1,e_{2}},l_{\mathcal{C}_1,e_{2}}\})$
 and $\tilde{u}\in
 f_2\setminus(h_2\cup\{f_{\mathcal{C}_2,f_{2}},l_{\mathcal{C}_2,f_{2}}\})$.
 Set $g_2=h'_2$, $\overline{g}_2=h''_2$ and
 $g'_2=\overline{g'}_2=h_2$.  Since there is no blue copy of
 $\mathcal{C}_{l_1+l_2}$, $E_2=g_2g'_2$ and
 $F_2=\overline{g}_2\overline{g'}_2$ are desired paths.  Therefore, we
 may assume that the edge $h_2$ is blue. For $3\leq l\leq
 \lfloor\frac{k}{2}\rfloor+1$, let
  $h_l=(h_{l-1}\setminus\{v_{k-2+l}\})\cup\{v_{2k-l+1}\}$.
  Assume that
  $l'$ is the maximum $l\in[2,\lfloor\frac{k}{2}\rfloor+1]$ for
  which $h_l$ is blue. If $l'<\lfloor\frac{k}{2}\rfloor+1$, then $h_{l'+1}$ is red. Set
  $$h'_{l'+1}=\Big((e_2\cup f_2\cup
 \{w\})\setminus(h_{l'}\cup
 \{f_{\mathcal{C}_1,e_{2}},f_{\mathcal{C}_2,f_{2}},\tilde{u}\})\Big)\cup\{\overline{u},\overline{v}\},$$
  where  $\tilde{u}\in
  f_2\setminus(h_{l'+1}\cup\{f_{\mathcal{C}_2,f_{2}},l_{\mathcal{C}_2,f_{2}}\})$. Clearly, $h'_{l'+1}$ is
  red. Set $g_2=\overline{g}_2=h'_{l'+1}$ and $g'_2=\overline{g'}_2=h_{l'+1}$.
 Since there is no blue copy of $\mathcal{C}_{l_1+l_2}$,
 $E_2=g_2g'_2$ and $F_2=\overline{g}_2\overline{g'}_2$ are desired
 paths. Also,
  $\mathcal{E}_2=\mathcal{E}_1{E}_2$ and
 $\mathcal{F}_2=\mathcal{E}_1F_2$ are two red paths of length $4$. Therefore, we
 may assume that $l'=\lfloor\frac{k}{2}\rfloor+1$ and hence the
 edge  $h_{\lfloor\frac{k}{2}\rfloor+1}=E' \cup F$ is
  blue, where
  $$E'=\{v_{2k-\lfloor\frac{k}{2}\rfloor},\ldots,v_{2k-2},v_{2k-1}\}.$$

 Now, set $m=\lfloor\frac{k}{2}\rfloor-1$.
 For $1\leq l\leq m$, let
  $h_{\lfloor\frac{k}{2}\rfloor+l+1}=(h_{\lfloor\frac{k}{2}\rfloor+l}\setminus\{u_{2(k-1)-l+1}\})\cup\{u_{k+l}\}$. Let
  $l'$ be the maximum $l\in[0,m]$ for
  which $h_{\lfloor\frac{k}{2}\rfloor+l+1}$ is blue. Similar to the
  above argument we can show that $l'=m$ and
  hence the
  edge $h_{\lfloor\frac{k}{2}\rfloor+m+1}=E'\cup F'$ is
  blue, where
 \begin{eqnarray*}
 F'=\left\lbrace\begin{array}{ll}
 \{u_{k+1},\ldots,u_{k+m},u_{2k-1}\} &\mbox{if $k$ is
 even},\vspace{.5
 cm}\\
 \{u_{k+1},\ldots,u_{k+m},u_{k+m+1},u_{2k-1}\} &\mbox{if $k$ is
 odd}.
 \end{array} \right.\vspace{.2
 cm}\end{eqnarray*}\\

 Let
 $h_{\lfloor\frac{k}{2}\rfloor+m+2}=(h_{\lfloor\frac{k}{2}\rfloor+m+1}\setminus\{u_{2k-1}\})\cup\{\bar{u}\}$.
  If $h_{\lfloor\frac{k}{2}\rfloor+m+2}$ is red, then set
 \begin{eqnarray*}
 h'_{\lfloor\frac{k}{2}\rfloor+m+2}=\Big((e_2\cup f_2\cup
 \{w\})\setminus(h_{\lfloor\frac{k}{2}\rfloor+m+1}\cup
 \{f_{\mathcal{C}_1,e_{2}},f_{\mathcal{C}_2,f_{2}},l_{\mathcal{C}_2,f_{2}},\tilde{v}\})\Big)\cup\{\overline{u},\overline{v}\},
 \end{eqnarray*}
  where  $\tilde{v}\in
  e_2\setminus(h_{\lfloor\frac{k}{2}\rfloor+m+2}\cup\{f_{\mathcal{C}_1,e_{2}},l_{\mathcal{C}_1,e_{2}}\})$.
   Since there is no blue copy of $\mathcal{C}_{l_1+l_2}$, the edge
  $h'_{\lfloor\frac{k}{2}\rfloor+m+2}$ is red.
  Set $g_2=\overline{g}_2=h'_{\lfloor\frac{k}{2}\rfloor+m+2}$ and
  $g'_2=\overline{g'}_2=h_{\lfloor\frac{k}{2}\rfloor+m+2}$.
   Therefore, $E_2=g_2g'_2$ and $F_2=\overline{g}_2\overline{g'}_2$ are  desired paths.
    So we may assume that the edge
  $h_{\lfloor\frac{k}{2}\rfloor+m+2}$ is blue.\\

 If $k$ is even, then clearly $h_{\lfloor\frac{k}{2}\rfloor+m+2}$
 is an edge in $\mathcal{B}_{22}$ disjoint from  $h_1$. This is
 impossible, by Remark \ref{complementary edges are not blue}.
  Now we may assume that $k$ is odd.
 One can easily see that $v_{k+\frac{k-1}{2}}\notin h_1\cup
 h_{\lfloor\frac{k}{2}\rfloor+m+2}$ and $u_{k+\frac{k-1}{2}}\in
 h_1\cap h_{\lfloor\frac{k}{2}\rfloor+m+2}$. Similarly,
  we can show
 that the edge
 $h_{\lfloor\frac{k}{2}\rfloor+m+3}=h_{\lfloor\frac{k}{2}\rfloor+m+2}\setminus
 \{u_{k+\frac{k-1}{2}}\})\cup\{v_{k+\frac{k-1}{2}}\}$ is blue. That
 is a contradiction to Remark \ref{complementary edges are not
 blue}. This contradiction completes the proof of
 Claim \ref{g is red:4}. $\square$\\\\

 Now, let $h'=\overline{E}\cup \overline{F}$ be an edge in $\mathcal{B}_{22}$ so that
 \begin{eqnarray*}
 && \overline{E}=\overline{E'}\cup\{l_{\mathcal{C}_1,e_{2}}\},\ \
 |\overline{E}|=\lceil\frac{k}{2}\rceil,\ \
  \overline{E'}\subseteq
 e_2\setminus\{f_{\mathcal{C}_1,e_{2}},l_{\mathcal{C}_1,e_{2}}\},\\
 && \overline{F}=\overline{F'}\cup\{\overline{u}\},\ \ |\overline{F}|=
 \lfloor\frac{k}{2}\rfloor,\ \ \overline{F'}\subseteq
 f_2\setminus\{f_{\mathcal{C}_2,f_{2}},l_{\mathcal{C}_2,f_{2}}\}.
 \end{eqnarray*}
 By an argument similar  to the proof of Claim \ref{g is red:4} we can show
 the following.
 \begin{emp}\label{g' is red:4}
 The edge  $h'$ is red.
 \end{emp}
  Now, by choosing edges $h$ and $h'$ appropriately
 as follows, we can find red paths $E_2$ and $F_2$ with desired
 properties. Let
 $$g_2=\overline{g}_2=h=\{x,v_{k+1},\ldots,v_{(k-1)+\lfloor\frac{k}{2}\rfloor}\}\cup
 \{u_{2k-\lceil\frac{k}{2}\rceil},\ldots,u_{2k-2},u_{2k-1}\}.$$ Set
 $E_2=g_2g'_2$ and $F_2=\overline{g}_2\overline{g'}_2$ where
 \begin{eqnarray*}
 g'_2&=&\{v_{k-1+\lfloor\frac{k}{2}\rfloor},v_{k+\lfloor\frac{k}{2}\rfloor+1},\ldots,v_{2k-2},l_{\mathcal{C}_1,e_{2}}\}
 \cup\{\overline{u},u_{k+1},\ldots,u_{k+\lfloor\frac{k}{2}\rfloor-1}\},\\
 \overline{g'}_2&=&\{v_{k+\lfloor\frac{k}{2}\rfloor},\ldots,v_{2k-2},l_{\mathcal{C}_1,e_{2}}\}
 \cup
 \{\overline{u},u_{k+1},\ldots,u_{k+\lfloor\frac{k}{2}\rfloor-2},u_{k+\lfloor\frac{k}{2}\rfloor}\}.
 \end{eqnarray*}
 Using Claims \ref{g is red:4} and \ref{g' is red:4}, $E_2=g_2g'_2$ and
 $F_2=\overline{g}_2\overline{g'}_2$ are  desired paths.
 Note that,
 $\mathcal{E}_2=\mathcal{E}_1E_2$ and
 $\mathcal{F}_2=\mathcal{E}_1F_2$ are two red paths of length $4$.
  This is a contradiction  to
 our
 assumption and so we are done. $\hfill\blacksquare$\\
 \medskip

 \medskip
 \noindent{\bf Proof of Lemma \ref{there is a P}.} Suppose for a contradiction that there is no red
 paths $E_i$ and $F_i$  with  desired  properties.
 Let $u'\in
 f_{i-1}\setminus(F_{i-1}\cup\{f_{\mathcal{C}_2,f_{i-1}}\})$,
  $v\in (e_{i-1}\setminus\{f_{\mathcal{C}_1,e_{i-1}}\})\cap(\overline{g'}_{i-1}\setminus
  \overline{g}_{i-1})$ and $w\in W\setminus
  (\bigcup_{j=1}^{i-1}B_j)$. Assume that
 $h=E\cup F$ is an edge in $\mathcal{A}_{ii}$ so that
 \begin{eqnarray*}
 && E=\overline{E}\cup\{v\},\ \ |E|=\lfloor\frac{k}{2}\rfloor,\ \
 \overline{E}\subseteq
 e_i\setminus\{f_{\mathcal{C}_1,e_{i}},l_{\mathcal{C}_1,e_{i}}\},\\
 && F=\overline{F}\cup\{l_{\mathcal{C}_2,f_{i}}\},\ \ |F|=
 \lceil\frac{k}{2}\rceil,\ \ \overline{F}\subseteq
 f_i\setminus\{f_{\mathcal{C}_2,f_{i}},l_{\mathcal{C}_2,f_{i}}\}.
 \end{eqnarray*}

 \begin{emp}\label{g is red:2}
 The edge $h$ is red.
 \end{emp}
 {\bf Proof of Claim \ref{g is red:2}.} By changing the indices  we
 may assume that
 \begin{eqnarray*}
 &&E=\{v,v_{(k-1)(i-1)+2},\ldots,v_{(k-1)(i-1)+\lfloor\frac{k}{2}\rfloor}\},\\
 &&F=\{u_{(k-1)i-\lceil\frac{k}{2}\rceil+2},\ldots,u_{(k-1)i},l_{\mathcal{C}_2,f_{i}}\}.
 \end{eqnarray*}
 Suppose indirectly that the edge $h_1=h$ is blue. Since there is
 no blue copy of   $\mathcal{C}_{l_1+l_2}$, using Remark
 \ref{complementary edges are not blue}, every  edge in
 $\mathcal{B}_{ii}$ that is disjoint from $h_1$ is red.
  Now let
 $h_2=(h_1\setminus\{v_{(k-1)(i-1)+2}\})\cup\{v_{(k-1)i}\}$. If
 $h_2$ is red, then set
 $$h'_2=\Big((e_i\cup f_i\cup
 \{w\})\setminus(h_1\cup
 \{f_{\mathcal{C}_1,e_{i}},f_{\mathcal{C}_2,f_{i}},l_{\mathcal{C}_2,f_{i}},\bar{u}\})\Big)\cup\{u'\},$$
 where $\bar{u}\in
 f_i\setminus(h_2\cup\{f_{\mathcal{C}_2,f_{i}},l_{\mathcal{C}_2,f_{i}}\})$.
 Set $g_i=\overline{g}_i=h_2$ and $g'_i=\overline{g'}_i=h'_2$.
 Since there is no blue copy of $\mathcal{C}_{l_1+l_2}$,
 $E_i=g_ig'_i$ and $F_i=\overline{g}_i\overline{g'}_i$ are desired
 paths (clearly, $v_{(k-1)(i-1)+2}\in e_i\setminus
 (E_i\cup\{f_{\mathcal{C}_1,e_{i}}\})$, $\bar{u}\in
 f_i\setminus(F_i\cup\{f_{\mathcal{C}_2,f_{i}}\})$. So for
 $\mathcal{P}=\mathcal{F}_{i-1}$, $\mathcal{P}{E}_i$ and
 $\mathcal{P}F_i$ are two red paths of length $2i$. Therefore, we
 may assume that the edge $h_2$ is blue. For $2\leq l\leq
 \lfloor\frac{k}{2}\rfloor$, let
  $h_l=(h_{l-1}\setminus\{v_{(k-1)(i-1)+l}\})\cup\{v_{(k-1)i-l+2}\}$.
  Assume that
  $l'$ is the maximum $l\in[1,\lfloor\frac{k}{2}\rfloor]$ for
  which $h_l$ is blue. If $l'<\lfloor\frac{k}{2}\rfloor$, then $h_{l'+1}$ is red. Set
  $$h'_{l'+1}=\Big((e_i\cup f_i\cup
 \{w\})\setminus(h_{l'}\cup
 \{f_{\mathcal{C}_1,e_{i}},f_{\mathcal{C}_2,f_{i}},l_{\mathcal{C}_2,f_{i}},\bar{u}\})\Big)\cup\{u'\},$$
  where  $\bar{u}\in
  f_i\setminus(h_{l'+1}\cup\{f_{\mathcal{C}_2,f_{i}},l_{\mathcal{C}_2,f_{i}}\})$. Clearly, $h'_{l'+1}$ is
  red. Set $g_i=\overline{g}_i=h_{l'+1}$ and $g'_i=\overline{g'}_i=h'_{l'+1}$.
 Since there is no blue copy of $\mathcal{C}_{l_1+l_2}$,
 $E_i=g_ig'_i$ and $F_i=\overline{g}_i\overline{g'}_i$ are desired
 paths
  (clearly $v_{(k-1)(i-1)+l'+1}\in e_i\setminus
 (E_i\cup\{f_{\mathcal{C}_1,e_{i}}\})$
  and $\bar{u}\in
  f_i\setminus(F_i\cup\{f_{\mathcal{C}_2,f_{i}}\})$). So for
 $\mathcal{P}=\mathcal{F}_{i-1}$, $\mathcal{P}{E}_i$ and
 $\mathcal{P}F_i$ are two red paths of length $2i$. Therefore, we
 may assume that $l'=\lfloor\frac{k}{2}\rfloor$ and hence
 the  edge $h_{\lfloor\frac{k}{2}\rfloor}=E'' \cup F$ is
  blue, where
  $$E''=\{v,v_{(k-1)i-\lfloor\frac{k}{2}\rfloor+2},\ldots,v_{(k-1)i-1},v_{(k-1)i}\}.$$

 Now, set $m=\lfloor\frac{k}{2}\rfloor-1$.
 For $1\leq l\leq m$, let
  $h_{\lfloor\frac{k}{2}\rfloor+l}=(h_{\lfloor\frac{k}{2}\rfloor+l-1}\setminus\{u_{(k-1)i-l+1}\})\cup\{u_{(k-1)(i-1)+l+1}\}$. Now, let
  $l'$ be the maximum $l\in[0,m]$ for
  which $h_{\lfloor\frac{k}{2}\rfloor+l}$ is blue. If $l'<m$, then $h_{\lfloor\frac{k}{2}\rfloor+l'+1}$ is red. Set
  $$h'_{\lfloor\frac{k}{2}\rfloor+l'+1}=\Big((e_i\cup f_i\cup
 \{w\})\setminus(h_{\lfloor\frac{k}{2}\rfloor+l'}\cup
 \{f_{\mathcal{C}_1,e_{i}},f_{\mathcal{C}_2,f_{i}},l_{\mathcal{C}_2,f_{i}},\bar{v}\})\Big)\cup\{u'\},$$
  where  $\bar{v}\in
  e_i\setminus(h_{\lfloor\frac{k}{2}\rfloor+l'+1}\cup\{f_{\mathcal{C}_1,e_{i}},l_{\mathcal{C}_1,e_{i}}\})$.
  Since there is no blue copy of $\mathcal{C}_{l_1+l_2}$, the edge
  $h'_{\lfloor\frac{k}{2}\rfloor+l'+1}$ is red. Set
  $g_i=\overline{g}_i=h_{\lfloor\frac{k}{2}\rfloor+l'+1}$and
  $g'_i=\overline{g'}_i=h'_{\lfloor\frac{k}{2}\rfloor+l'+1}$.
   Therefore, $E_i=g_ig'_i$ and $F_i=\overline{g}_i\overline{g'}_i$ are desired paths
   (clearly $\bar{v}\in e_i\setminus
 (E_i\cup\{f_{\mathcal{C}_1,e_{i}}\})$ and $u_{(k-1)i-l'}\in
 f_i\setminus(F_i\cup\{f_{\mathcal{C}_2,f_{i}}\})$) and  for
 $\mathcal{P}=\mathcal{F}_{i-1}$, $\mathcal{P}E_i$ and
 $\mathcal{P}F_i$ are two red paths of length $2i$. So we may
 assume that $l'=m$ and
  hence the
  edge $h_{\lfloor\frac{k}{2}\rfloor+m}=E''\cup F''$ is
  blue, where
 \begin{eqnarray*}
 F''=\left\lbrace\begin{array}{ll}
 \{u_{(k-1)(i-1)+2},\ldots,u_{(k-1)(i-1)+m+1},l_{\mathcal{C}_2,f_{i}}\}
 &\mbox{if $k$ is even},\vspace{.5
 cm}\\
 \{u_{(k-1)(i-1)+2},\ldots,u_{(k-1)(i-1)+m+1},u_{(k-1)(i-1)+m+2},l_{\mathcal{C}_2,f_{i}}\}
 &\mbox{if $k$ is odd}. \end{array} \right.\vspace{.2
 cm}\end{eqnarray*}\\

 Let
 $h_{\lfloor\frac{k}{2}\rfloor+m+1}=(h_{\lfloor\frac{k}{2}\rfloor+m}\setminus\{v\})\cup\{l_{\mathcal{C}_1,e_{i}}\}$.
  If $h_{\lfloor\frac{k}{2}\rfloor+m+1}$ is red, then set
 \begin{eqnarray*}
 &&\hspace{-1.5 cm}h'_{\lfloor\frac{k}{2}\rfloor+m+1}=\Big((e_i\cup
 f_i\cup \{w\})\setminus(h_{\lfloor\frac{k}{2}\rfloor+m}\cup
 \{f_{\mathcal{C}_1,e_{i}},f_{\mathcal{C}_2,f_{i}},l_{\mathcal{C}_2,f_{i}},\bar{v}\})\Big)\cup\{u\},\\
 \vspace*{0.4 cm} &&\hspace{-1.5
 cm}h''_{\lfloor\frac{k}{2}\rfloor+m+1}=\Big((e_i\cup f_i\cup
 \{w\})\setminus(h_{\lfloor\frac{k}{2}\rfloor+m}\cup
 \{f_{\mathcal{C}_1,e_{i}},f_{\mathcal{C}_2,f_{i}},l_{\mathcal{C}_2,f_{i}},\bar{u}\})\Big)\cup\{u\}.
 \end{eqnarray*}
  where  $\bar{v}\in
  e_i\setminus(h_{\lfloor\frac{k}{2}\rfloor+m+1}\cup\{f_{\mathcal{C}_1,e_{i}},l_{\mathcal{C}_1,e_{i}}\})$,
   $\bar{u}\in  f_j\setminus(h_{\lfloor\frac{k}{2}\rfloor+m+1}\cup\{f_{\mathcal{C}_2,f_{j}},l_{\mathcal{C}_2,f_{j}}\})$ and
   $u\in (f_{j-1}\setminus\{f_{\mathcal{C}_2,f_{j-1}}\})\cap(\overline{g'}_{i-1}\setminus
   \overline{g}_{i-1})$. Set
   $g'_i=\overline{g'}_i=h_{\lfloor\frac{k}{2}\rfloor+m+1}$,
   $g_i=h'_{\lfloor\frac{k}{2}\rfloor+m+1}$and
   $\overline{g}_i=h''_{\lfloor\frac{k}{2}\rfloor+m+1}$.
 It is easy to see that $E_i=g_ig'_i$ and
 $F_i=\overline{g}_i\overline{g'}_i$ are the desired paths and
 hence for $\mathcal{P}=\mathcal{F}_{i-1}$,
 $\mathcal{E}_{i}=\mathcal{P}E_i$ and
 $\mathcal{F}_{i}=\mathcal{P}F_i$ are two favorable red paths of
 length $2i$. So, we may assume that the edge
 $h_{\lfloor\frac{k}{2}\rfloor+m+1}$ is blue.\\

  Now, let
 $h_{\lfloor\frac{k}{2}\rfloor+m+2}=(h_{\lfloor\frac{k}{2}\rfloor+m+1}\setminus\{l_{\mathcal{C}_2,f_{i}}\})\cup\{u'\}$.
 If $h_{\lfloor\frac{k}{2}\rfloor+m+2}$ is red, then set
 $$h'_{\lfloor\frac{k}{2}\rfloor+m+2}=\Big((e_i\cup
 f_i\cup \{w\})\setminus(h_{\lfloor\frac{k}{2}\rfloor+m+1}\cup
 \{f_{\mathcal{C}_1,e_{i}},l_{\mathcal{C}_1,e_{i}},f_{\mathcal{C}_2,f_{i}},l_{\mathcal{C}_2,f_{i}},\bar{v}\})\Big)\cup\{v,u'\},$$
  where  $\bar{v}\in
  e_i\setminus(h_{\lfloor\frac{k}{2}\rfloor+m+2}\cup\{f_{\mathcal{C}_1,e_{i}},l_{\mathcal{C}_1,e_{i}}\})$.
   Since there is no blue copy of $\mathcal{C}_{l_1+l_2}$, the edge
  $h'_{\lfloor\frac{k}{2}\rfloor+m+2}$ is red.
  Set $g_i=\overline{g}_i=h'_{\lfloor\frac{k}{2}\rfloor+m+2}$ and
  $g'_i=\overline{g'}_i=h_{\lfloor\frac{k}{2}\rfloor+m+2}$.
   Therefore, $E_i=g_ig'_i$ and $F_i=\overline{g}_i\overline{g'}_i$ are the desired path
   (clearly $\bar{v}\in e_i\setminus(E_i\cup\{f_{\mathcal{C}_1,e_{i}}\})$  and
    $l_{\mathcal{C}_2,f_{i}}\in
    f_i\setminus(F_i\cup\{f_{\mathcal{C}_2,f_{i}}\})$. So  for $\mathcal{P}=\mathcal{F}_{i-1}$,
 $\mathcal{E}_{i}=\mathcal{P}E_i$ and
 $\mathcal{F}_{i}=\mathcal{P}F_i$ are two red paths of length $2i$.
    So we may assume that the edge
  $h_{\lfloor\frac{k}{2}\rfloor+m+2}$ is blue.\\

 If $k$ is even, then clearly $h_{\lfloor\frac{k}{2}\rfloor+m+2}$
 is an edge in $\mathcal{B}_{ii}$ disjoint from  $h_1$. This is
 impossible, by Remark \ref{complementary edges are not blue}.
  Now we may assume that $k$ is odd.
 One can easily see that $v_{(k-1)(i-1)+\frac{k+1}{2}}\notin
 h_1\cup h_{\lfloor\frac{k}{2}\rfloor+m+2}$ and
 $u_{(k-1)(i-1)+\frac{k+1}{2}}\in h_1\cap
 h_{\lfloor\frac{k}{2}\rfloor+m+2}$. Similarly,
  we can show
 that the edge
 $h_{\lfloor\frac{k}{2}\rfloor+m+3}=h_{\lfloor\frac{k}{2}\rfloor+m+2}\setminus
 \{u_{(k-1)(i-1)+\frac{k+1}{2}}\})\cup\{v_{(k-1)(i-1)+\frac{k+1}{2}}\}$
 is blue. That is a contradiction to Remark \ref{complementary
 edges are not blue}. This contradiction completes the proof of
 Claim \ref{g is red:2}. $\square$\\\\

 Now, let $h'=E'\cup F'$ be an edge in $\mathcal{B}_{ii}$ so that
 \begin{eqnarray*}
 && E'=\overline{E'}\cup\{l_{\mathcal{C}_1,e_{i}}\},\ \
 |E'|=\lceil\frac{k}{2}\rceil,\ \
  \overline{E'}\subseteq
 e_i\setminus\{f_{\mathcal{C}_1,e_{i}},l_{\mathcal{C}_1,e_{i}}\},\\
 && F'=\overline{F'}\cup\{u'\},\ \ |F'|=
 \lfloor\frac{k}{2}\rfloor,\ \ \overline{F'}\subseteq
 f_i\setminus\{f_{\mathcal{C}_2,f_{i}},l_{\mathcal{C}_2,f_{i}}\}.
 \end{eqnarray*}

 \begin{emp}\label{g' is red:2}
 The edge  $h'$ is red.
 \end{emp}
 {\bf Proof of Claim \ref{g' is red:2}.} By changing  the indices
 we may assume that
 \begin{eqnarray*}
 &&E'=\{v_{(k-1)i-\lceil\frac{k}{2}\rceil+2},\ldots,v_{(k-1)i},l_{\mathcal{C}_1,e_{i}}\},\\
 &&F'=\{u',u_{(k-1)(i-1)+2},\ldots,u_{(k-1)(i-1)+\lfloor\frac{k}{2}\rfloor}\}.
 \end{eqnarray*}
 Suppose indirectly that the edge $h_1=h'$ is blue. Let
 $m=\lfloor\frac{k}{2}\rfloor$.
 \noindent  For $2\leq l\leq m$,  let
  $h_l=(h_{l-1}\setminus\{v_{(k-1)i-l+2}\})\cup\{v_{(k-1)(i-1)+l}\}$.
  Similar to the proof of Claim \ref{g is red:2}, we can show that
  the edge $h_m=E''\cup F'$ is blue where
  \begin{eqnarray*}
 \hspace*{-1 cm} E''=\left\lbrace \begin{array}{ll}
 \{v_{(k-1)(i-1)+2},\ldots,v_{(k-1)(i-1)+m},l_{\mathcal{C}_1,e_{i}}\}
 &\mbox{if $k$ is even},\vspace{.5 cm}\\
 \{v_{(k-1)(i-1)+2},\ldots,v_{(k-1)(i-1)+m},v_{(k-1)(i-1)+m+1},l_{\mathcal{C}_1,e_{i}}\}
 &\mbox{if $k$ is odd}. \end{array}\right. \end{eqnarray*}\\

  Now, set $m'=\lfloor\frac{k}{2}\rfloor-1$.
 For $1\leq l\leq m'$, let $h_{m+l}=(h_{m+l-1}\setminus
 \{u_{(k-1)(i-1)+l+1}\})\cup\{u_{(k-1)i+1-l}\}$. By an argument
 similar to the proof of Claim \ref{g is red:2}  we may assume that
 the
  edge $h_{m+m'}=E''\cup F''$ is
  blue, where
 $$F''=\{u',u_{(k-1)i-m'+1},\ldots,u_{(k-1)i}\}.$$

 Let
 $h_{m+m'+1}=(h_{m+m'}\setminus\{u'\})\cup\{l_{\mathcal{C}_2,f_{i}}\}$.
  If $h_{m+m'+1}$ is red, then set
 \begin{eqnarray*}
 &&\hspace{-1.5 cm}h'_{m+m'+1}=\Big((e_i\cup f_i\cup
 \{w\})\setminus(h_{m+m'}\cup
 \{f_{\mathcal{C}_1,e_{i}},l_{\mathcal{C}_1,e_{i}},f_{\mathcal{C}_2,f_{i}},\bar{v}\})\Big)\cup\{v\},\\
 \vspace*{0.4 cm} &&\hspace{-1.5 cm}h''_{m+m'+1}=\Big((e_i\cup
 f_j\cup \{w\})\setminus(h_{m+m'}\cup
 \{f_{\mathcal{C}_1,e_{i}},l_{\mathcal{C}_1,e_{i}},f_{\mathcal{C}_2,f_{i}},\bar{u}\})\Big)\cup\{v\}.
 \end{eqnarray*}
  where  $\bar{v}\in
  e_i\setminus(h_{m+m'+1}\cup\{f_{\mathcal{C}_1,e_{i}},l_{\mathcal{C}_1,e_{i}}\})$
  and
   $\bar{u}\in
   f_i\setminus(h_{m+m'+1}\cup\{f_{\mathcal{C}_2,f_{i}},l_{\mathcal{C}_2,f_{i}}\})$.
 Set $g'_i=\overline{g'}_i=h_{m+m'+1}$, $g_i=h'_{m+m'+1}$ and
 $\overline{g}_i=h''_{m+m'+1}$. It is easy to see that
 $E_i=g_ig'_i$ and $F_i=\overline{g}_i\overline{g'}_i$ are  desired
 paths. Clearly for $\mathcal{P}=\mathcal{F}_{i-1}$,
 $\mathcal{E}_i=\mathcal{P}E_i$ and $\mathcal{F}_i=\mathcal{P}F_i$
 are two red paths of length $2i$.
  Hence, we may assume that the edge $h_{m+m'+1}$ is blue.\\

 Let  $v'\in
  e_{i-1}\setminus(E_{i-1}\cup\{f_{\mathcal{C}_1,e_{i-1}}\})$ and
 $h_{m+m'+2}=(h_{m+m'+1}\setminus\{l_{\mathcal{C}_1,e_{i}}\})\cup\{v'\}$.
  If $h_{m+m'+2}$ is red, then set
 \begin{eqnarray*}
 &&\hspace{-1.5 cm}h'_{m+m'+2}=\Big((e_i\cup f_i\cup
 \{w\})\setminus(h_{m+m'+1}\cup
 \{f_{\mathcal{C}_1,e_{i}},l_{\mathcal{C}_1,e_{i}},f_{\mathcal{C}_2,f_{i}},l_{\mathcal{C}_2,f_{i}},\bar{v}\})\Big)\cup\{\hat{u},v'\},\\
 \vspace*{0.4 cm} &&\hspace{-1.5 cm}h''_{m+m'+2}=\Big((e_i\cup
 f_i\cup \{w\})\setminus(h_{m+m'+1}\cup
 \{f_{\mathcal{C}_1,e_{i}},l_{\mathcal{C}_1,e_{i}},f_{\mathcal{C}_2,f_{i}},l_{\mathcal{C}_2,f_{i}},\bar{u}\})\Big)\cup\{\hat{u},v'\}.
 \end{eqnarray*}
  where  $\bar{v}\in
  e_i\setminus(h_{m+m'+2}\cup\{f_{\mathcal{C}_1,e_{i}},l_{\mathcal{C}_1,e_{i}}\})$,
   $\bar{u}\in
   f_i\setminus(h_{m+m'+2}\cup\{f_{\mathcal{C}_2,f_{i}},l_{\mathcal{C}_2,f_{i}}\})$ and $\hat{u}\in (f_{i-1}\setminus\{f_{\mathcal{C}_2,f_{i-1}}\})
   \cap(g'_{i-1}\setminus {g}_{i-1})$. Set
   $g'_i=\overline{g'}_i=h_{m+m'+2}$, $g_i=h'_{m+m'+2}$ and
   $\overline{g}_i=h''_{m+m'+2}$. Then $E_i=g_ig'_i$ and
   $F_i=\overline{g}_i\overline{g'}_i$ are desired paths (note that $\bar{v}\in e_i\setminus
 (E_i\cup\{f_{\mathcal{C}_1,e_{i}}\})$ and $\bar{u}\in
 f_i\setminus(F_i\cup\{f_{\mathcal{C}_2,f_{i}}\})$) and so for
 $\mathcal{P}=\mathcal{E}_{i-1}$,
   $\mathcal{E}_i=\mathcal{P}E_i$ and
   $\mathcal{F}_i=\mathcal{P}F_i$ are two red paths of length
   $2i$, a contradiction to our assumption.
  Hence, we
 may assume that the edge $h_{m+m'+2}$ is blue. \\

 If $k$ is even, then clearly $h_{m+m'+2}$
 is an edge in $\mathcal{A}_{ii}$ disjoint from  $h_1$. This is
 impossible, by Remark \ref{complementary edges are not blue}.
  Now we may assume that $k$ is odd.
 One can easily see that $v_{(k-1)(i-1)+\frac{k+1}{2}}\in h_1\cap
 h_{m+m'+2}$ and $u_{(k-1)(i-1)+\frac{k+1}{2}}\notin h_1\cup
 h_{m+m'+2}$. Similarly,
  we can show
 that the edge $$h_{m+m'+3}=(h_{m+m'+2}\setminus
 \{v_{(k-1)(i-1)+\frac{k+1}{2}}\})\cup\{u_{(k-1)(i-1)+\frac{k+1}{2}}\}$$
 is blue. Hence, $h_{m+m'+3}$ is an edge in $\mathcal{A}_{ii}$
 disjoint from $h_1$,
  that is a contradiction to Remark \ref{complementary
 edges are not blue}. This contradiction
 completes the proof of Claim \ref{g' is red:2}. $\square$\\\\

  Now, by choosing edges $h$ and $h'$ appropriately
 as follows, we can find red paths $E_i$ and $F_i$ with desired
 properties. Let
 $$g_i=\{v,v_{(k-1)(i-1)+2},\ldots,v_{(k-1)(i-1)+\lfloor\frac{k}{2}\rfloor}\}
 \cup\{u_{(k-1)i-\lceil\frac{k}{2}\rceil+2},\ldots,u_{(k-1)i},l_{\mathcal{C}_2,f_{i}}\}.$$
 Set $E_i=g_ig'_i$ and $F_i=\overline{g}_i\overline{g'}_i$ where
 $\overline{g}_i=g_i$ and
 \begin{eqnarray*}
 g'_i&=&\{v_{(k-1)(i-1)+\lfloor\frac{k}{2}\rfloor},v_{(k-1)(i-1)+\lfloor\frac{k}{2}\rfloor+2},v_{(k-1)(i-1)+\lfloor\frac{k}{2}\rfloor+3},\ldots,v_{(k-1)i},l_{\mathcal{C}_1,e_{i}}\}\\
 &&\cup\{u',u_{(k-1)(i-1)+2},\ldots,u_{(k-1)(i-1)+\lfloor\frac{k}{2}\rfloor}\},\\
 \overline{g'}_i&=&\{v_{(k-1)(i-1)+\lfloor\frac{k}{2}\rfloor+1},\ldots,v_{(k-1)i},l_{\mathcal{C}_1,e_{i}}\}\\
 &&\cup
 \{u',u_{(k-1)(i-1)+2},\ldots,u_{(k-1)(i-1)+\lfloor\frac{k}{2}\rfloor-1},u_{(k-1)(i-1)+\lfloor\frac{k}{2}\rfloor+1}\}.
 \end{eqnarray*}
 Using Claims \ref{g is red:2} and \ref{g' is red:2}, $E_i$ and
 $F_i$ are  desired paths (clearly,
 $v_{(k-1)(i-1)+\lfloor\frac{k}{2}\rfloor+1}\in
 e_i\setminus(E_i\cup \{f_{\mathcal{C}_1,e_{i}}\})$ and
 $u_{(k-1)(i-1)+\lfloor\frac{k}{2}\rfloor}\in
 f_i\setminus(F_i\cup\{f_{\mathcal{C}_2,f_{i}}\})$). Note that,
 $\mathcal{E}_i=\mathcal{P}E_i$ and $\mathcal{F}_i=\mathcal{P}F_i$
 are two red paths of length $2i$ where
 $\mathcal{P}=\mathcal{F}_{i-1}$. This is a contradiction  to our
 assumption and so we are done. $\hfill\blacksquare$\\
 \medskip

 \noindent\textbf{Proof of Lemma \ref{No  blue Cn implies red Cn+1/2}. } Let $\mathcal{C}_1=e_1e_2\ldots
 e_{\frac{n-1}{2}}$ and $\mathcal{C}_2=f_1f_2\ldots
 f_{\frac{n+1}{2}}$ be copies of $\mathcal{C}_{\frac{n-1}{2}}^k$
 and $\mathcal{C}_{\frac{n+1}{2}}^k$ in $\mathcal{H}_{\rm blue}$
 with edges
 $$e_i=\{v_1,v_2,\ldots,v_k\}+(k-1)(i-1) \Big({\rm mod}\ (k-1)(\frac{n-1}{2})\Big),\hspace{0.5cm} i=1,2,\ldots,\frac{n-1}{2}$$
 and
 $$f_i=\{u_1,u_2,\ldots,u_k\}+(k-1)(i-1) \Big({\rm mod}\  (k-1)(\frac{n+1}{2})\Big),\hspace{0.5cm} i=1,2,\ldots,\frac{n+1}{2}.$$

 \noindent Let $W=V(\mathcal{H})\setminus V(\mathcal{C}_1\cup
 \mathcal{C}_2).$ First assume that $n=5.$ Since there is no blue
 copy of $\mathcal{C}^k_5$, consider $i=j=1$, $e=f_1$,
 $B=W=\{w_1,w_2\}$, $C=\{v_{k-1}\}$, $v'=v_1$, $u'=u_1$, $v''=v_k$,
 $u''=u_k$ and use Lemma
  \ref{there is a P2}, to obtain a red path
  $\mathcal{P}=g_1g_2$ with the mentioned properties in Lemma \ref{there is a
  P2}.
  Clearly,  $V(\mathcal{P})\subseteq \Big(V(e_1\cup f_1)\setminus\{v_{k-1}\}\Big)\cup W$ and $g_1\cap W\neq
  \emptyset$. With no loss of generality assume that $w_1\in
 W\cap g_1$. Set
 $$h=(e_2\setminus
 \{v_k,v_{2k-2},v_{1}\})\cup \{x_1,u_{3k-3},y_1\},$$ and
 $$h'=(f_3\setminus \{u_{3k-4},u_{3k-3},u_1\})\cup
 \{w_1,v_{2k-2},x_2\},$$

 \noindent where $x_1\in (g_1\setminus
 g_2)\cap(e_1\setminus\{v_1\}),$ $y_1\in (g_2\setminus
 g_1)\cap(f_1\setminus \{u_k\})$ and $x_2\in (g_2\setminus
 g_1)\cap(e_1\setminus \{v_k\}).$ Since there is no blue copy of
 $\mathcal{C}_{5}^k,$ at least one of $\mathcal{P}h$
 or $\mathcal{P}h'$ is a red $\mathcal{C}^k_3$ and so  we are done.\\

  Now, let $n\geq 7$. We find a red $\mathcal{C}^k_{\frac{n+1}{2}}$  as
  follows.\\\\
 {\bf Step 1:} Put $i=j=1$, $e=f_1$,  $C=\emptyset$, $v'=v_1$,
 $u'=u_1$, $v''=v_k$, $u''=u_k$ and use Lemma
  \ref{there is a P2}, to obtain a red path
  $\mathcal{P}=g_1g_2$ with the mentioned properties in Lemma \ref{there is a
  P2}.
 So  $V(\mathcal{P}_1)\subseteq V(e_1\cup
 f_1)\cup W$ and $|V(\mathcal{P}_1)\cap W|\leq 1$.\\\\
 %
 %
 %

 \noindent {\bf Step 2:} For $n=7$ go to step 3. Otherwise, do the
 following process $\frac{n-7}{2}$ times.
 For $1\leq i\leq \frac{n-7}{2}$, let $$h_{i+1}=(e_{i+1}\setminus
 \{f_{\mathcal{C}_1,e_{i+1}}, v_{(i+1)(k-1)},
 l_{\mathcal{C}_1,e_{i+1}}\})\cup
 $$$$\{x_{i+1}, u_{(i+1)(k-1)}, l_{\mathcal{C}_2,f_{i+1}}\},$$ and $$h'_{i+1}=
 (f_{i+1}\setminus \{f_{\mathcal{C}_2,f_{i+1}}, u_{(i+1)(k-1)},
 l_{\mathcal{C}_2,f_{i+1}}\}) \cup
 $$$$\{y_{i+1}, v_{(i+1)(k-1)}, l_{\mathcal{C}_1,e_{i+1}}\},$$ where $x_{i+1}\in
 (g_{i+1}\setminus g_{i})\cap (e_{i}\setminus
 \{f_{\mathcal{C}_1,e_{i}}\})$ and $y_{i+1}\in (g_{i+1}\setminus
 g_{i})\cap (f_{i}\setminus \{f_{\mathcal{C}_2,f_{i}}\})$ are
 vertices with maximum indices. One can easily check that at least
 one of $h_{i+1}$ or $h'_{i+1}$, say $g_{i+2}$, is red. Set
 $\mathcal{P}_{i+1}=\mathcal{P}_{i}g_{i+2}$.
  Clearly $\mathcal{P}_{i+1}$ is a red path of length
  $i+2$.\\

 \noindent{\bf Step 3:}
   Let $$h_{\frac{n-3}{2}}=\big(e_{\frac{n-1}{2}}\setminus
  \{f_{\mathcal{C}_1,e_{\frac{n-1}{2}}},v_{\frac{n-3}{2}(k-1)+2},v_{1}\}\big)\cup $$ $$
  \{x_{\frac{n-3}{2}},u_{\frac{n-3}{2}(k-1)},l_{\mathcal{C}_2,f_{\frac{n-3}{2}}}\},$$ and  $$h'_{\frac{n-3}{2}}=\big(f_{\frac{n-3}{2}}\setminus
 \{f_{\mathcal{C}_2,f_{\frac{n-3}{2}}}, u_{\frac{n-3}{2}(k-1)},
 l_{\mathcal{C}_2,f_{\frac{n-3}{2}}}\}\big)\cup
 $$ $$ \{y_{\frac{n-3}{2}}, f_{\mathcal{C}_1,e_{\frac{n-1}{2}}}, v_{\frac{n-3}{2}(k-1)+2}\},$$ where
 $x_{\frac{n-3}{2}}\in (g_1\setminus g_2)\cap (e_1\setminus
 \{v_k\})$ is the vertex with minimum indices and
 $y_{\frac{n-3}{2}}\in (g_{\frac{n-3}{2}}\setminus
 g_{\frac{n-5}{2}})\cap (f_{\frac{n-5}{2}}\setminus
 \{f_{\mathcal{C}_2,f_{\frac{n-5}{2}}}\})$ is the vertex with
 maximum indices. If $h_{\frac{n-3}{2}}$ is red, then set
 $$h_{\frac{n-1}{2}}=(e_{\frac{n-3}{2}}\setminus \{f_{\mathcal{C}_1,e_{\frac{n-3}{2}}}, v_{\frac{n-3}{2}(k-1)},
 l_{\mathcal{C}_1,e_{\frac{n-3}{2}}}\})\cup$$$$
 \{x_{\frac{n-1}{2}}, u_{\frac{n-5}{2}(k-1)+k-2},
 l_{\mathcal{C}_2,f_{\frac{n-3}{2}}}\},$$ and
 $$h'_{\frac{n-1}{2}}=(f_{\frac{n-3}{2}}\setminus \{f_{\mathcal{C}_2,f_{\frac{n-3}{2}}}, u_{\frac{n-5}{2}(k-1)+k-2},
 l_{\mathcal{C}_2,f_{\frac{n-3}{2}}}\})\cup$$
 $$ \{y_{\frac{n-1}{2}}, v_{\frac{n-3}{2}(k-1)}, l_{\mathcal{C}_1,e_{\frac{n-3}{2}}}\},$$
 where $x_{\frac{n-1}{2}}\in (g_{\frac{n-3}{2}}\setminus
 g_{\frac{n-5}{2}})\cap (e_{\frac{n-5}{2}}\setminus
 \{f_{\mathcal{C}_1,e_{\frac{n-5}{2}}}\})$ and
 $y_{\frac{n-1}{2}}\in (g_{\frac{n-3}{2}}\setminus
 g_{\frac{n-5}{2}})\cap (f_{\frac{n-5}{2}}\setminus
 \{f_{\mathcal{C}_2,f_{\frac{n-5}{2}}}\})$. So clearly at least
 one of
 $\mathcal{P}_{\frac{n-5}{2}}h_{\frac{n-1}{2}}h_{\frac{n-3}{2}}$ or
 $\mathcal{P}_{\frac{n-5}{2}}h'_{\frac{n-1}{2}}h_{\frac{n-3}{2}}$
 is a red $\mathcal{C}_{\frac{n+1}{2}}^k$. Now, assume that
 $h_{\frac{n-3}{2}}$ is blue. If $h'_{\frac{n-3}{2}}$ is blue, then
 $$e_1e_2\ldots
 e_{\frac{n-3}{2}}h'_{\frac{n-3}{2}}f_{\frac{n-5}{2}}f_{\frac{n-7}{2}}\ldots
 f_1 f_{\frac{n+1}{2}}f_{\frac{n-1}{2}}h_{\frac{n-3}{2}},$$ is a
 blue copy of  $\mathcal{C}_n^k$, a contradiction. So
 $h'_{\frac{n-3}{2}}$ is red. Therefore, set
 $$h_\frac{n-1}{2}=(e_\frac{n-1}{2}\setminus \{v_{\frac{n-3}{2}(k-1)+2}, v_{1}\})\cup
 \{u_{\frac{n+1}{2}(k-1)},y_{\frac{n-1}{2}}\},$$ and
 $$h'_{\frac{n-1}{2}}=(f_{\frac{n+1}{2}}\setminus \{u_{\frac{n+1}{2}(k-1)}, u_1\})\cup
 \{v_{\frac{n-3}{2}(k-1)+2},x_\frac{n-1}{2}\},$$ where
 $x_{\frac{n-1}{2}} \in (g_1\setminus g_2)\cap (e_1\setminus
 \{v_{k}\})$ and $y_{\frac{n-1}{2}} \in (g_1\setminus g_2)\cap
 (f_1\setminus \{u_k\})$ are vertices with minimum indices. Since
 at least one of $h_{\frac{n-1}{2}}$ and $h'_{\frac{n-1}{2}}$, say
 $\overline{h}_{\frac{n-1}{2}}$, is red, then
 $\mathcal{P}_{\frac{n-5}{2}}h'_{\frac{n-3}{2}}\overline{h}_{\frac{n-1}{2}}$
 is a red $\mathcal{C}_{\frac{n+1}{2}}^k$, which completes the
 proof.

 $\hfill\blacksquare$\\

 \noindent\textbf{Proof of Lemma \ref{No  red Cn implies blue Cn/2+1:2}. } Let $\mathcal{C}_1=e_1e_2\ldots
 e_{\frac{n}{2}-1}$ and $\mathcal{C}_2=f_1f_2\ldots
 f_{\frac{n}{2}+1}$ be copies of $\mathcal{C}_{\frac{n}{2}-1}^k$
 and $\mathcal{C}_{\frac{n}{2}+1}^k$ in $\mathcal{H}_{\rm blue}$
 with edges
 \begin{eqnarray*}
 &&\hspace{-1.5 cm}e_i=\{v_1,v_2,\ldots,v_k\}+(k-1)(i-1) \Big({\rm
 mod}\ (k-1)(\frac{n}{2}-1)\Big),\hspace{0.5 cm}
 i=1,2,\ldots,\frac{n}{2}-1,\\
 &&\hspace{-1.5 cm}f_i=\{u_1,u_2,\ldots,u_k\}+(k-1)(i-1) ({\rm
 mod}\ (k-1)(\frac{n}{2}+1)),\hspace{0.5 cm}
 i=1,2,\ldots,\frac{n}{2}+1,
 \end{eqnarray*}
 and  $W=V(\mathcal{H})\setminus V(\mathcal{C}_1\cup
 \mathcal{C}_2)$. First let $n\geq 8$. We can find a red copy of $\mathcal{C}_{\frac{n}{2}+1}^k$ as follows.\\\\
 {\bf Step 1:} Put $i=j=1$, $B=\{w_1,w_2\}\subseteq W$,
 $C=\{v_{k-1}\}$, $v'=v_1$, $u'=u_1$, $v''=v_k$ and $u''=u_k$. Then
  use Lemma
  \ref{there is a P2} for $e=e_1$ (resp. $e=f_1$), to obtain a red path
  $E_1$ (resp. $F_1$) of length $2$ with the mentioned properties in Lemma \ref{there is a
  P2}.
 Now, use Lemma \ref{there is a P} for $i=2$, to obtain two red
 paths $E_2$ and $F_2$ of length $2$ with the mentioned properties
 in Lemma \ref{there is a P}. So $\mathcal{E}_2=\mathcal{P}E_2$ is
 a red path of length $4$ for some $\mathcal{P}\in \{E_1,F_1\}$. Assume that $\mathcal{P}_2=\mathcal{E}_2=g_1g_2g_3g_4$.
 With no loss of generality we may assume that $w_1\in g_1\setminus g_2$.\\\\
 {\bf Step 2:} If $n=8$ go to step 3. Otherwise, do the following
 process $\frac{n}{2}-4$ times. For $1\leq i\leq \frac{n}{2}-4$,
 let

 $$h_{i+2}=(e_{i+2}\setminus \{f_{\mathcal{C}_1,e_{i+2}},
 v_{(i+2)(k-1)}, l_{\mathcal{C}_1,e_{i+2}}\})\cup
 $$$$\{x_{i+2}, u_{(i+2)(k-1)}, l_{\mathcal{C}_2,f_{i+2}}\}$$ and $$h'_{i+2}=
 (f_{i+2}\setminus \{f_{\mathcal{C}_2,f_{i+2}}, u_{(i+2)(k-1)},
 l_{\mathcal{C}_2,f_{i+2}}\}) \cup
 $$$$\{y_{i+2}, v_{(i+2)(k-1)}, l_{\mathcal{C}_1,e_{i+2}}\},$$ where $x_{i+2}\in
 (g_{i+3}\setminus g_{i+2})\cap (e_{i+1}\setminus
 \{f_{\mathcal{C}_1,e_{i+1}}\})$ and $y_{i+1}\in (g_{i+3}\setminus
 g_{i+2})\cap (f_{i+1}\setminus \{f_{\mathcal{C}_2,f_{i+1}}\})$ are
 vertices with maximum indices. Since at least one of $h_{i+2}$ and
 $h'_{i+2}$, say $g_{i+4}$, is red,
  clearly $\mathcal{P}_{i+2}=\mathcal{P}_{i+1}g_{i+4}$ is a red path of length
  $i+4$.\\\\
 {\bf Step 3:}
  Set
 $$h_{\frac{n}{2}-1}=\big(e_{\frac{n}{2}-1}\setminus
  \{f_{\mathcal{C}_1,e_{\frac{n}{2}-1}},v_{(\frac{n}{2}-1)(k-1)},v_{1}\}\big)\cup $$ $$
  \{x_{\frac{n}{2}-1},w_1,l_{\mathcal{C}_2,f_{\frac{n}{2}-1}}\}$$ and  $$h'_{\frac{n}{2}-1}=\big(f_{\frac{n}{2}-1}\setminus
 \{f_{\mathcal{C}_2,f_{\frac{n}{2}-1}},
 u_{(\frac{n}{2}-1)(k-1)},l_{\mathcal{C}_2,f_{\frac{n}{2}-1}}\}\big)\cup
 $$ $$ \{y_{\frac{n}{2}-1}, v_{(\frac{n}{2}-1)(k-1)}, x'_{\frac{n}{2}-1}\},$$ so that
 $x'_{\frac{n}{2}-1}\in (g_1\setminus g_2) \cap (e_1\setminus
 \{v_k\})$ is a vertex with minimum index and
 $x_{\frac{n}{2}-1}\in(g_{\frac{n}{2}}\setminus
 g_{\frac{n}{2}-1})\cap (e_{\frac{n}{2}-2}\setminus
 \{f_{\mathcal{C}_1,e_{\frac{n}{2}-2}}\})$, $y_{\frac{n}{2}-1}\in
 (g_{\frac{n}{2}}\setminus g_{\frac{n}{2}-1})\cap
 (f_{\frac{n}{2}-2}\setminus
 \{f_{\mathcal{C}_2,f_{\frac{n}{2}-2}}\})$ are vertices with
 maximum indices. Clearly, at least one of
 $\mathcal{P}_{\frac{n}{2}-2}h_{\frac{n}{2}-1}$ or
 $\mathcal{P}_{\frac{n}{2}-2}h'_{\frac{n}{2}-1}$ is a red
 $\mathcal{C}_{\frac{n}{2}+1}^k$.\\\\

 Now, let $n=6$. Suppose indirectly that there is no blue copy of
 $\mathcal{C}^k_{4}$. Put $i=j=1$, $e=f_1$, $B=W=\{w_1,w_2\}$,
 $C=\{v_{k-1}\}$, $v'=v_1$, $u'=u_1$, $v''=v_k$, $u''=u_k$ and use
 Lemma
  \ref{there is a P2}, to obtain a red path
  $\mathcal{P}=g_1g_2$ with
  so that $V(\mathcal{P})\subseteq \Big(V(e_1\cup f_1)\setminus\{v_{k-1}\}\Big)\cup B$, $(g_1\setminus g_2)\cap B\neq
  \emptyset$ and there is a vertex $y\in
 f_1\setminus\Big(V(\mathcal{P})\cup \{u_1\}\Big)$. With no loss of
 generality assume that
  $w_1\in g_1\setminus g_2$.
  Let
  $x\in(e_{1}\setminus\{v_1\})\cap(g_2\setminus
  g_1)$, $x'\in (e_1\setminus\{v_k\})\cap (g_1\setminus g_2)$ and
 $h=E\dot{\cup}F$ is an edge in $\mathcal{A}_{22}$ so that
 \begin{eqnarray*}
 &&E=\{x,v_{k+1},\ldots,v_{k+\lfloor\frac{k}{2}\rfloor-1}\},\\
 &&F=\{u_{2k-\lceil\frac{k}{2}\rceil},\ldots,u_{2k-1}\}.
 \end{eqnarray*}

 \begin{emp}\label{g is red:3}
 The edge $h$ is red.
 \end{emp}
 {\bf Proof of Claim \ref{g is red:3}.} Suppose indirectly that the
 edge $h_1=h$ is blue. Since there is no blue copy of
 $\mathcal{C}^k_{6}$, using Remark \ref{complementary edges are not
 blue}, every  edge in $\mathcal{B}_{22}$ that is disjoint from
 $h_1$ is red.
  For $2\leq l\leq
 \lfloor\frac{k}{2}\rfloor$, let
  $h_l=(h_{l-1}\setminus\{v_{k-1+l}\})\cup\{v_{2(k-1)-l+2}\}$.
  Assume that
  $l'$ is the maximum $l\in[1,\lfloor\frac{k}{2}\rfloor]$ for
  which $h_l$ is blue. If $l'<\lfloor\frac{k}{2}\rfloor$, then $h_{l'+1}$ is red. Set
  $$h'_{l'+1}=\Big((e_2\cup f_2)\setminus(h_{l'}\cup
 \{v_{k},v_1,u_{k},l_{\mathcal{C}_2,f_{2}}\})\Big)\cup\{y,x'\}.$$
    Since there is no blue copy of $\mathcal{C}^k_{6}$, the edge $h'_{l'+1}$ is
  red and  $\mathcal{P}h_{l'+1}h'_{l'+1}$ is a red copy of
  $\mathcal{C}_4^k$. So we may assume that $l'=\lfloor\frac{k}{2}\rfloor$ and  the edge
  $h_{\lfloor\frac{k}{2}\rfloor}=E' \cup F$ is
  blue, where
  $$E'=\{x,v_{2k-\lfloor\frac{k}{2}\rfloor},\ldots,v_{2k-3},v_{2k-2}\}.$$

 Now, set $m=\lfloor\frac{k}{2}\rfloor-1$. For $1\leq l\leq m$, let
  $h_{\lfloor\frac{k}{2}\rfloor+l}=(h_{\lfloor\frac{k}{2}\rfloor+l-1}\setminus\{u_{2(k-1)-l+1}\})\cup\{u_{k+l}\}$. Now, let
  $l'$ be the maximum $l\in[0,m]$ for
  which $h_{\lfloor\frac{k}{2}\rfloor+l}$ is blue. If $l'<m$, then $h_{\lfloor\frac{k}{2}\rfloor+l'+1}$ is red. Set
  $$h'_{\lfloor\frac{k}{2}\rfloor+l'+1}=\Big((e_2\cup f_2)\setminus(h_{\lfloor\frac{k}{2}\rfloor+l'}\cup
 \{v_1,v_k,u_{k},l_{\mathcal{C}_2,f_{2}}\})\Big)\cup\{y,x'\}.$$
  Since there is no blue copy of $\mathcal{C}^k_{6}$, the edge
  $h'_{\lfloor\frac{k}{2}\rfloor+l'+1}$ is red and $\mathcal{P}h_{\lfloor\frac{k}{2}\rfloor+l'+1}h'_{\lfloor\frac{k}{2}\rfloor+l'+1}$ is a red copy of
  $\mathcal{C}_4^k$. So we may assume that
  $l'=m$ and
  hence the
  edge $h_{\lfloor\frac{k}{2}\rfloor+m}=E'\cup F'$ is
  blue, where
  \begin{eqnarray*} F'=\left\lbrace\begin{array}{ll}
 \{u_{k+1},\ldots,u_{k+m},u_{2k-1}\}
 &\mbox{$k$ is even},\vspace{.5 cm}\\
 \{u_{k+1},\ldots,u_{k+m},u_{k+m+1},u_{2k-1}\} &\mbox{$k$ is odd}.
 \end{array}\right.
 \end{eqnarray*}

 Let $y'\in (g_2\setminus g_1)\cap(f_1\setminus\{u_1\})$. Now,
 consider the edge
 $h_{\lfloor\frac{k}{2}\rfloor+m+1}=(h_{\lfloor\frac{k}{2}\rfloor+m}\setminus\{x,u_{2k-1}\})\cup\{v_{k-1},y'\}$.
 If $h_{\lfloor\frac{k}{2}\rfloor+m+1}$ is red, then set
 $$h'_{\lfloor\frac{k}{2}\rfloor+m+1}=\Big((e_2\cup
 f_2)\setminus(h_{\lfloor\frac{k}{2}\rfloor+m}\cup
 \{v_{k},v_1,u_{k},l_{\mathcal{C}_2,f_{2}},\hat{u}\})\Big)\cup\{y,v_{k-1},w_1\},$$
 where $\hat{u}\in f_2\setminus
 (h_{\lfloor\frac{k}{2}\rfloor+m}\cup\{u_k,u_{2k-1}\})$. Clearly,
 $\mathcal{P}h_{\lfloor\frac{k}{2}\rfloor+m+1}h'_{\lfloor\frac{k}{2}\rfloor+m+1}$
 is a red copy of $\mathcal{C}^k_4$, a contradiction to our
 assumption. So we may assume that
 $h_{\lfloor\frac{k}{2}\rfloor+m+1}$ is blue.
  If $k$ is even, then clearly $h_{\lfloor\frac{k}{2}\rfloor+m+1}$
  is an edge in $\mathcal{B}_{22}$ disjoint from $h_1$, that is a
  contradiction to Remark \ref{complementary edges are not blue}.

  Now, we may assume that $k$ is odd.
 One can easily check that $v_{k+\lfloor\frac{k}{2}\rfloor}\notin
 h_1\cup h_{\lfloor\frac{k}{2}\rfloor+m+1}$ and
 $u_{k+\lfloor\frac{k}{2}\rfloor}\in h_1\cap
 h_{\lfloor\frac{k}{2}\rfloor+m+1}$. By an argument similar to the
 above, we can show that the edge
 $h_{\lfloor\frac{k}{2}\rfloor+m+2}=(h_{\lfloor\frac{k}{2}\rfloor+m+1}\setminus\{u_{k+\lfloor\frac{k}{2}\rfloor}\})\cup\{v_{k+\lfloor\frac{k}{2}\rfloor}\}$
 is blue. That is a contradiction to Remark \ref{complementary
 edges are not blue}, since $h_{\lfloor\frac{k}{2}\rfloor+m+2}$ is
 a blue edge in $\mathcal{B}_{22}$ disjoint from $h_1$. This
 contradiction completes the proof of Claim \ref{g is red:3}.
 $\square$\\\\

 Now, let $h'=\overline{E}\dot{\cup}\overline{F}$ be an edge in
 $\mathcal{B}_{22}$ so that
 \begin{eqnarray*}
 &&\overline{E}=\{v_{2k-\lceil\frac{k}{2}\rceil},\ldots,v_{2k-2},x'\},\\
 &&\overline{F}=\{y,u_{k+1},\ldots,u_{k+\lfloor\frac{k}{2}\rfloor-1}\}.
 \end{eqnarray*}

 \begin{emp}\label{g' is red:3}
 The edge  $h'$ is red.
 \end{emp}
 Suppose indirectly that the edge $h_1=h'$ is blue. Let
 $m=m'=\lfloor\frac{k}{2}\rfloor-1$.
  Similar to the proof of Claim \ref{g is red:3}, we can show that
  the edge $h_{m+m'+1}=\overline{E'}\cup \overline{F'}$ is blue where
 \begin{eqnarray*}
 \overline{F'}=\{y,u_{2k-\lfloor\frac{k}{2}\rfloor},\ldots,u_{2k-2}\}
 \end{eqnarray*}
 \noindent and
 \begin{eqnarray*}
  \overline{E'}=\left\lbrace \begin{array}{ll}
 \{v_{k+1},\ldots,v_{k+m},x'\}
 &\mbox{$k$ is even},\vspace{.5 cm}\\
 \{v_{k+1},\ldots,v_{k+m},v_{k+m+1},x'\} &\mbox{$k$ is odd}.
 \end{array}\right. \end{eqnarray*}

  Now, consider the edge
  $h_{m+m'+2}=(h_{m+m'+1}\setminus\{y\})\cup\{u_{2k-1}\}$. If
  $h_{m+m'+2}$ is red, then set
 \begin{eqnarray*}
 h'_{m+m'+2}=\Big((e_2\cup f_2)\setminus(h_{m+m'+1}\cup
 \{v_{k},v_1,u_{k}\})\Big)\cup\{x\}.
 \end{eqnarray*}
 It is easy to see that  $h'_{m+m'+2}$ is red  and
 $\mathcal{P}h'_{m+m'+2}h_{m+m'+2}$ is a red copy of
 $\mathcal{C}^k_4$, a contradiction to our assumptions. So we may assume that $h_{m+m'+2}$ is blue.\\

 Now, let
 $$h_{m+m'+3}=(h_{m+m'+2}\setminus\{x',v_{k+\lceil\frac{k}{2}\rceil-1}\})\cup\{v_{k-1},w_1\}.$$
 If $h_{m+m'+3}$ is red, then set
 $$h'_{m+m'+3}=\Big((e_2\cup f_2)\setminus(h_{m+m'+2}\cup
 \{v_{k},v_1,u_{k}\})\Big)\cup\{v_{k-1},y'\},$$ where
 $y'\in(f_1\setminus\{u_1\})\cap(g_2\setminus g_1).$ Since there is
 no blue copy of $\mathcal{C}^6_k$, then $h'_{m+m'+3}$ is red and
 hence $\mathcal{P}h'_{m+m'+3}h_{m+m'+3}$ is a red copy of
 $\mathcal{C}^k_4$, a contradiction to our assumption. So we may
 assume that the edge $h_{m+m'+3}$ is blue. Notice that
 $h_{m+m'+3}$ is a blue edge in $\mathcal{A}_{22}$ disjoint from
 $h_1$, a contradiction to Remark \ref{complementary edges are not
 blue}. This contradiction completes the proof of Claim \ref{g' is
 red:3}.$\square$\\\\

 Now, let
 \begin{eqnarray*}
 &&h=\{x,v_{k+1},\ldots,v_{k+\lfloor\frac{k}{2}\rfloor-1}\}\cup\{u_{2k-\lceil\frac{k}{2}\rceil},\ldots,u_{2k-1}\},\\
 &&h'=\{v_{k+\lfloor\frac{k}{2}\rfloor-1},v_{k+\lfloor\frac{k}{2}\rfloor+1},\ldots,v_{2k-2},x'\}\cup
 \{u',u_{k+1},\ldots,u_{k+\lfloor\frac{k}{2}\rfloor-1}\}.
 \end{eqnarray*}
 Using Claims \ref{g is red:3} and \ref{g' is red:3}, the edges $h$
 and $h'$ are red. So $\mathcal{P}hh'$ is a red copy of $\mathcal{C}^k_4$. This contradicts  our assumption. So we are done. $\hfill\blacksquare$\\ \end{document}